\documentclass[a4paper,reqno]{amsart}
\usepackage{latexsym}
\usepackage{amssymb,amsthm,amsmath}
\usepackage[dvips]{graphicx}
\usepackage{xcolor}

\calclayout

\theoremstyle{plain}
\newtheorem{theorem}{Theorem}
\newtheorem{lemma}[theorem]{Lemma}
\newtheorem{corollary}[theorem]{Corollary}
\newtheorem{proposition}[theorem]{Proposition}
\newtheorem{example}[theorem]{Example}

\theoremstyle{definition}
\newtheorem{definition}[theorem]{Definition}
\theoremstyle{remark}
\newtheorem{remark}[theorem]{Remark}

\DeclareMathOperator*{\esssup}{\rm ess\,sup}

\def\endmark{\hskip 2em\begin{picture}(8,10)
\put(0,0){$\Box$} \put(2,0){\rule{1.9mm}{0.3mm}}
\put(6.5,0){\rule{0.3mm}{1.9mm}}
\end{picture}
\par}

\def\endproof{\null\hfill\endmark\endtrivlist}

\def\d#1{{#1\kern-0.4em\char"16\kern-0.1em}}
\def\D#1{{\raise0.2ex\hbox{-}\kern-0.4em #1}}
\reversemarginpar

\newcounter{zd}

\newcounter{zdr}[subsection]

\newcommand{\eps}{\varepsilon}

\def\mff{{\mathfrak f}}

\def\Div{{\rm div}}

\def\pa{\partial}

\def\cal{\mathcal}
\let\mib=\boldsymbol

\def\R{\mathbb{R}}    
\def\N{\mathbb{N}}    

\def\pC#1#2{{{\rm C}^{#1}(#2)}}

\def\eps{\varepsilon}

\def\malpha{{\mib \alpha}}
\def\mbeta{{\mib \beta}}
\def\mgamma{{\mib \gamma}}
\def\mdelta{{\mib \delta}}
\def\mkappa{{\mib \kappa}}

\def\mx{{\bf x}}
\def\mxi{{\mib \xi}}
\def\meta{{\mib \eta}}

\def\my{{\bf y}}
\def\mz{{\bf z}}

\def\oi#1#2{(#1,#2)}   

\def\Rd{{\mathbb{R}^{d}}}    
\def\Sdmj{{\rm S}^{d-1}}

\def\supp{{\rm supp\,}}

\def\la{\lambda}

\begin{document}

\title[Degenerate parabolic equations]{Velocity averaging for diffusive transport equations with discontinuous flux}

\author{M.~Erceg}\address{Marko Erceg,
Department of Mathematics, Faculty of Science, University of Zagreb, Bijeni\v{c}ka cesta 30,
10000 Zagreb, Croatia}\email{maerceg@math.hr}

\author{M.~Mi\v{s}ur}\address{Marin Mi\v{s}ur,
Department of Mathematics, Faculty of Science, University of Zagreb, Bijeni\v{c}ka cesta 30,
10000 Zagreb, Croatia}\email{mmisur@math.hr}

\author{D.~Mitrovi\'{c}}\address{Darko Mitrovi\'{c},
University of Vienna, Faculty of Mathematics, Oscar Morgenstern platz 1,
1090 Vienna, Austria}\email{darko.mitrovic@univie.ac.at}



\begin{abstract}
We consider a diffusive transport equation with discontinuous flux and prove the velocity averaging result under non-degeneracy conditions. In order to achieve the result, we introduce a new variant of micro-local defect functionals which are able to ``recognise'' changes of the type of the equation. As a corollary, we show the existence of a weak solution for the Cauchy problem for nonlinear degenerate parabolic equation with discontinuous flux. We also show existence of strong traces at $t=0$ for so-called quasi-solutions to degenerate parabolic equations under non-degeneracy conditions on the diffusion term.
\end{abstract}

\maketitle

\section{Introduction}\label{sec:intro}

In \cite[Theorem C]{LPT} a result on velocity averaging for diffusive transport equations has been stated, but the proof of the theorem cannot be found neither in that paper or in later contributions (we shall provide a more detailed insight later in the introduction). The aim of the paper is to precisely proof \cite[Theorem C]{LPT} in the $\mathrm{L}^q$-setting, $q>2$, and to generalise the result on equations with discontinuous coefficients.

To be more precise, we aim to prove a velocity averaging result for a diffusive transport equation with discontinuous flux meaning that for the sequence $(u_n)$
of solutions to the sequence of equations of the form
\begin{equation}\label{app-1}
\begin{split}
\Div_\mx \bigl(f(\mx,\lambda) u_n(\mx,\lambda) \bigr) =& \,\Div_\mx \bigl(\Div_\mx \left( a(\lambda) u_n(\mx,\lambda) \right)\bigr) \\
&\qquad\quad + \pa_{\lambda} G_n(\mx,\lambda) + \Div_\mx P_n(\mx,\lambda) \qquad \hbox{in} \ {\cal D}'(\R^{d+1})\;,
\end{split}
\end{equation}
for every $\rho \in {\rm C}^1_c(\R)$, the sequence $\left(\int_{\R} \rho(\lambda) u_n(\mx,\lambda) \,d\lambda\right)$ is strongly precompact in ${\rm L}^1_{loc}(\R^d)$ (i.e. it lies in a compact subset of ${\rm L}^1_{loc}(\R^d)$). Precise assumptions on the coefficients are given in a)-e) below.

Equation \eqref{app-1} has two main components. The transport part
$$
\Div_\mx( f(\mx,\lambda) u_n(\mx,\lambda))
$$and the diffusive part
$$
\Div_\mx \bigl(\Div_\mx \left( a(\lambda) u_n(\mx,\lambda) \right)\bigr)=\Div_\mx \bigl( a(\lambda) \nabla_{\mx} u_n(\mx,\lambda) \bigr) \;,
$$ 
where $u_n(\mx,\lambda)$ is unknown, $a(\lambda)\in \R^{d\times d}$ is the diffusion matrix, $f(\mx,\lambda)$ is the flux, $\mx \in \R^d$ is the space (and time) 
variable and $\lambda$ is called the velocity variable, but it can be considered as a parameter.
For the sake of generality and simplicity of the exposition, we compressed 
the space-time variable into a single variable $\mx$, while still our main intention 
is to study evolution equations (see Remark \ref{rem:app-1-parbolic}). 
In the literature, velocity variable $\lambda$ is often denoted by $v$. 
The form of the remaining source terms in \eqref{app-1} is motivated by 
the kinetic formulation for degenerate parabolic equations, as can be seen in Section 5 and Section 6. 

The transport component $\Div_\mx( f(\mx,\lambda) u_n(\mx,\lambda))$ is a generalisation of the usual kinetic transport term $\langle v\,|\, \nabla_{\mx} h(\mx,v)\rangle$, i.e. the equation
\begin{equation*}
\label{agosh}
\pa_t h+\langle v\,|\, \nabla_{\mx} h\rangle=\Div_\mx \pa_v^\mkappa g, \ \ (t,\mx)\in (0,\infty)\times \R^d, \ \ v\in \R^d, \ \ g\in \mathrm{L}^2(\R^d\times \R^d) 
	\,, \ \mkappa\in\N^d \,,
\end{equation*} 
for which the velocity averaging results was proved in \cite{Ago}. Independently of \cite{Ago}, the corresponding results were discovered in \cite{nGPS} and further extended in \cite{nGLPS}. The mentioned results were given in the $\mathrm{L}^2$-setting. In \cite{nDpLM}, one can find the first ${\rm L}^p$, $p>1$, velocity averaging result obtained using the approach of multiplier operators  (see e.g.~\cite{Gra}). The optimal result in the sense of the ${\rm L}^p$-integrability of $(u_n)$ has been achieved in \cite{nB, Per}, while an ${\rm L}^2$ velocity averaging result for pseudo-differential operators can be found in \cite{nGG}.

Such a type of result appeared to be very useful and it was a substantial part of the proof of existence of the weak solution to the Boltzmann
equation \cite{DpL} as well as the regularity of admissible solutions to scalar conservation laws \cite{LPT}. In \cite{LPT}, one can also find the first result concerning the velocity averaging for the transport equations with the flux of the form $f=f(\lambda)$, $f \in \mathrm{C}(\R;\R^d)$, under the non-degeneracy conditions which essentially mean that for any $\mxi\in\mathbb{R}^{d}\setminus\{0\}$, the mapping
\begin{equation*}
\label{n-d}
\lambda \mapsto \langle f(\lambda)\,|\,\mxi\rangle
\end{equation*} is possibly zero only on a negligible set.

As for the non-hyperbolic situation ($a\ne 0$), the velocity averaging results for ultra-parabolic equations are proven in \cite{LM2},
while for the degenerate parabolic equations, i.e.~the ones in which $a$ changes rank for different $\lambda$, by our best knowledge, the only results
can be found in \cite{GH, LPT, TT} for the homogeneous flux $f$ and diffusion matrix $a$ (i.e.~both independent of $\mx$) in the $\mathrm{L}^p$-settig for any $p>1$. 
Let us note here again that in \cite{LPT} details of the proof are not provided (see \cite[Theorem C]{LPT}) since the authors conjectured that the proof could be accomplished by following the method from \cite{nGLPS}. 
However, such a method is actually applied in \cite{TT} and the authors needed an additional assumption (see \cite[(2.20)]{TT}) to finalise the arguments (see more precise explanation below).

Let us now briefly explain a main idea of the technique from \cite{GH, LPT, TT}. Since both flux and diffusion matrix are independent of $\mx$, this enables a separation of coefficients and unknown functions by means of the Fourier transform. Indeed, if $f$ is independent of $\mx$, by applying the Fourier transform to \eqref{app-1} with respect to $\mx$ one sees that
\begin{equation}\label{eq:informal_calc_fourier}
\hat{u}_n(\mxi,\lambda)=\frac{\sigma^2 |\mxi|^2 \hat{u}_n+i \langle\hat{P}_n\,|\,\mxi
	\rangle
	+ \pa_\lambda \hat{G}_n }{\sigma^2 |\mxi|^2+i \langle f\,|\, \mxi\rangle+
	2\pi\langle a \,\mxi \,|\,\mxi\rangle} \;,
\end{equation}
where we denoted by $\mxi$ the dual variable (the definition of the Fourier
transform used here can be found in Notation below). 
In the case when $a\equiv 0$, informally speaking, from \eqref{eq:informal_calc_fourier}:

\begin{itemize}

\item by controlling the term $\hat{u}_n$ on the right-hand side of the latter expression by the constant $\sigma$;

\item by integrating by parts with respect to $\lambda$ to remove the derivative from the functions $\hat{G}_n$;

\item by employing the non-degeneracy conditions \eqref{n-d};

\end{itemize} 
one can draw appropriate conclusions on the sequence $(u_n)$.

The generalisation on the situation when $a\neq 0$ is not straightforward. First we need to assume that 
\begin{equation}
\label{par-non-deg}
(\forall \mxi \in \R^d\setminus\{0\}) \qquad 
	\operatorname{meas} \bigl\{ \lambda\in K\subset\subset \R:\, 
	\langle f(\lambda)\,|\, \mxi\rangle 
	=\langle a(\lambda)\mxi \,|\, \mxi\rangle=0 \bigr\}=0 \,,
\end{equation} which are the non-degeneracy conditions corresponding to \eqref{app-1} (with the flux independent of $\mx$). 
 
However, since the integration by parts with respect to $\lambda$ 
(which is the second step in the procedure above) 
affects the non-negativity of the matrix $a$, it seems that additional assumptions 
on $a$ are needed in order to conclude about the strong compactness of the velocity averages.

In particular, such a result can be found in \cite{TT}, which is aimed to the regularity properties of the velocity averages (more precisely, $\mathrm{W}^{s,r}$-regularity, $s>0$, $r\geq 1$). In the essence of the proofs is the method described above (separation of the solution $u$ from coefficients) together with the so-called {truncation property} \cite[Definition 2.1]{TT} (see \cite[Lemma 2.3]{TT}) under a variant of assumption \eqref{par-non-deg} and an assumption on behaviour of the $\lambda$-derivative of the symbol ${\cal L}(\mxi,\lambda)=i \langle f(\lambda)\,|\,\mxi\rangle+2\pi\langle a(\lambda) \mxi\,|\,\mxi\rangle$ of equation \eqref{app-1} on \emph{layers} in $\mxi$-space defined by the Littlewood-Palley decomposition. In \cite{GH}, the results are repeated in the stochastic setting.  

We also mention results from \cite{gess, GST} where one can find velocity averaging results for degenerate parabolic equations obtained as a kinetic reformulation of the porous media equation.

Before stating our main result, let us first fix the notation used in the paper. 
\smallskip

\noindent \textbf{Notation.}
Throughout the paper we denote by $\langle\cdot\,|\,\cdot\rangle$ the complex Euclidean scalar product on $\mathbb{C}^d$, 
which we take to be antilinear in the second argument. However, in our situations we shall mainly work on $\Rd$. 
By $|\cdot|$ we denote the corresponding norm of vectors, while the same notation is used for the 2-norm
for matrices.  For a matrix $A$, by $A^T$ we denote its transpose. For the complex conjugate of a complex number $z$ 
we use $\bar z$.

By $\mx=(x_1,x_2,\dots,x_d)$ we write points (vectors) in $\Rd$,
while by $\mxi=(\xi_1,\xi_2,\dots,\xi_d)$ we denote the dual variables in the sense of the Fourier transform (if $t$ occurs, then we use
$\tau$ for the dual variable). The Fourier transform we define by $\hat{u}(\mxi) = {\cal F}u(\mxi) = \int_{\Rd} e^{-2\pi i \langle\mxi\,|\,\mx\rangle}u(\mx)\,d\mx$, and its inverse by $(u)^\vee(\mx) = \bar{\cal F}u(\mx) = \int_{\Rd} e^{2\pi i \langle\mxi\,|\,\mx\rangle}u(\mxi)\,d\mx$, while the Fourier multiplier operator by ${\cal A}_\psi u= (\psi\hat u)^\vee$. 
If ${\cal A}_\psi$ is bounded on ${\rm L}^p(\Rd)$ we call it the ${\rm L}^p$-Fourier multiplier operator and $\psi$ the 
${\rm L}^p$-Fourier multiplier. We will often have that $\psi$ depends (besides $\mxi$) on $\lambda$ which is always considered as a parameter.

For a Lebesgue measurable subset $A\subseteq\Rd$ we denote by $\mathop{Cl}A$, 
$\mathop{Int}A$, $cA$, $\operatorname{meas}(A)$ and $\chi_A$ 
the closure of $A$, the interior of $A$, the complement of $A$, 
the Lebesgue measure of $A$, and the characteristic function over $A$, respectively. 
The open (closed) ball in $\Rd$ centered at point $\mx$ with radius $r > 0$ we will denote by $B(\mx, r)$
($B[\mx, r]$), the unit sphere in $\Rd$ by $\Sdmj$, and in Section \ref{sec:application} we will use the shorthand $\R^{d+1}_+:=\R^+\times\Rd$.
The signum function is denoted by ${\rm sgn}$.

For a multi-index $\malpha=(\alpha_1,\alpha_2,\dots,\alpha_d)\in\N_0^d$ we denote by $|\malpha|=\alpha_1+\alpha_2+\cdots+\alpha_d$ its
length and by $\pa^\malpha=\pa_{x_1}^{\alpha_1}\pa_{x_2}^{\alpha_2}\cdots\pa_{x_d}^{\alpha_d}$ partial derivatives. 

By ${\rm L}^p_{loc}(\Omega)$, $\Omega\subseteq\Rd$ open and $p\in[1,\infty]$, we denote the Fr\'echet space of functions that 
are contained in ${\rm L}^p(\Omega')$ for any compactly contained set $\Omega'$ in $\Omega$ ($\Omega'\subset\subset\Omega$), and 
analogously for Sobolev spaces ${\rm W}^{s,p}_{loc}(\Omega)$, $s\in\R$, $p\in[1,\infty]$. 
${\rm C}_c(X)$ stands for the space of compactly supported continuous functions on a locally compact space $X$.
If $X$ is compact then ${\rm C}_c(X)={\rm C}(X)$. For the space of Lipschitz functions we use ${\rm C}^{0,1}(X)$.
Any dual product is denoted by $\langle\cdot,\cdot\rangle$, which we take to 
be linear in both arguments. ${\cal L}(X)$ stands for the space of bounded linear operators on a normed space $X$. 

When applicable, functions defined on a subset of $\Rd$ shall often be identified by their extensions by zero to the whole space.
\smallskip

In order to introduce the main results of the paper, we need the following assumptions on $(u_n)$ and the coefficients appearing in \eqref{app-1}:
\smallskip

\noindent\textbf{Assumptions}
\begin{itemize}
\item [a)] $(u_n)$ is uniformly compactly supported on open $\Omega\times S\subset\subset \R^{d}_\mx\times\R_\lambda$, $d\geq 2$, and weakly (weakly-$\star$ for $q=\infty$) converges to zero in 
${\rm L}^q(\R^d\times \R)$ for some $q\in(2,\infty]$;

\item [b)] $a=\sigma^T\sigma$, where $\sigma\in{\rm C}^{0,1}(S; \R^{d\times d})$;
\item [c)] $f\in {\rm L}^p(\Omega\times S; \Rd)$ for some $p>\frac{q}{q-1}$ 
($p>1$ if $q=\infty$), and
for any compact $K\subseteq S$ it holds
\begin{equation}\label{lp-111}
\esssup\limits_{\mx\in \Omega} \sup\limits_{\mxi \in \mathrm{S}^{d-1}} 
	\operatorname{meas}\bigl\{\lambda \in K :\, \langle f(\mx,\lambda)\,|\,\mxi
	\rangle=\langle a(\lambda)\mxi \,|\,\mxi \rangle=0 \bigr\}=0 \;; 
\end{equation}
\item [d)] $G_n \to 0$ strongly in ${\rm L}_{loc}^{r_0}(\R_\lambda; {\rm W}_{loc}^{-1/2,r_0}(\R^d_\mx))$ for some 
	$r_0\in(1,\infty)$;
\item [e)] $P_n=(P_1^n,\dots,P_d^n)\to 0$ strongly in ${\rm L}_{loc}^{p_0}(\R^d_\mx\times\R_\lambda;\Rd)$ for some $p_0\in(1,\infty)$.
\end{itemize}

Our main result is the following velocity averaging result for \eqref{app-1}.

\begin{theorem}\label{velocity averaging}
	Let $d\geq 2$ and let $(u_n)$ satisfies (a) and the sequence of equations \eqref{app-1} whose coefficients satisfy conditions (b)-(e).
	
	Then there exists a subsequence $(u_{n'})$ such that for any $\rho\in {\rm C}_c(S)$,
	\begin{equation}
	\label{cnv*}
	\int_S \rho(\lambda) u_{n'}(\mx,\lambda) \,d\lambda \,\longrightarrow\, 0 \quad\hbox{strongly in} \ {\rm L}^1_{loc}(\Rd) \;.
		\end{equation}
\end{theorem}

The theorem above generalises the compactness results of \cite{LPT, TT} to 
the case of the flux discontinuous with respect to the space variable,
while the diffusion matrix remains homogeneous, i.e.~dependent only on $\lambda$.
Moreover, the non-degeneracy condition \eqref{lp-111} can be seen as 
a natural generalisation of \eqref{par-non-deg} to the heterogeneous setting.
The heterogeneity prevents us of using the above explained method based on the 
Fourier transform, thus in our proof we follow the approach of \cite{Ger, LM2},
which is elaborated below. 
However, we are not able to obtain the result for $(u_n)$ bounded in ${\rm L}^p$ 
if $p \leq 2$, as achieved in \cite{LPT, TT}.

Assumption (b) from the above 
can be relaxed (see Remark \ref{rem:weaker_assump}), 
and the proofs remain essentially the same. The following corollary holds.

\begin{corollary}
\label{cor1}
If we replace in the assumptions of Theorem \ref{velocity averaging}
condition (b) by a weaker assumption:
\begin{enumerate}
\item[$\tilde{\rm b}$)] the mapping $\lambda \mapsto a(\lambda)=\sigma^{T}(\lambda)\sigma(\lambda) \in \R^{d\times d}$ is such that for almost every $\lambda_0\in S\subset\subset \R$ there exists $\eps>0$ such that $\sigma\in \mathrm{C}^{0,1}((\lambda_0-\eps,\lambda_0+\eps);\R^{d\times d})$;
\end{enumerate}
the statement of Theorem \ref{velocity averaging} still holds.
\end{corollary}

Moreover, under stronger assumptions on $(G_n)$ we have the 
following result. 

\begin{corollary}\label{velocity averaging_Gn_in_Lr}
If we replace in the assumptions of Theorem \ref{velocity averaging}
conditions (b) and (d) by
\begin{enumerate}
\item[b')] $a\in {\rm C}^{0,1}(S; \R^{d\times d})$ is such that, 
for every $\lambda\in S$, $a(\lambda)$ is a symmetric and positive 
semi-definite matrix;
\item [d')] $G_n \to 0$ strongly in ${\rm L}_{loc}^{r_0}(\R^d_\mx\times\R_\lambda)$ for some 
$r_0\in(1,\infty)$;
\end{enumerate}
the statement of Theorem \ref{velocity averaging} holds. 
\end{corollary}

Comparing the result of Corollary \ref{cor1},
when 
applied to the homogeneous setting ($f=f(\lambda)$), 
to \cite[Theorem C]{LPT}, one can see that the former does not 
reveal completely the latter, 
where only smoothness and positive semi-definiteness of $a$ is required, i.e.~(b') instead of ($\tilde{\rm b}$).
Nevertheless, conditions (b) and ($\tilde{\rm b}$) 
still cover many interesting cases of 
the degenerate diffusion effects. Let us illustrate this on the following example.

\begin{example}
\label{ex3}
\begin{enumerate}
\item[a)]
It is clear that all matrix functions of the form 
$a(\lambda)=Q(\lambda)^T\Lambda(\lambda) Q(\lambda)$
satisfy condition (b), 
where, for any $\lambda\in S$, $Q(\lambda)$ is orthogonal
and $\Lambda(\lambda)$ is positive semi-definite and diagonal, 
and $Q,\sqrt{\Lambda}\in\mathrm{C}^{0,1}(S;\R^{d\times d})$.
Indeed, in this case we can take $\sigma(\lambda)=\sqrt{\Lambda(\lambda)}Q(\lambda)$.

For instance,
\begin{equation}\label{eq:matrix_2d}
a(\lambda) =
\left(
\frac{1}{\sqrt{\lambda^2+1}}
\begin{bmatrix}
\lambda& 1\\
1& -\lambda
\end{bmatrix}
\right)
\begin{bmatrix}
0& 0\\
0& \lambda^2+1
\end{bmatrix}
\left(
\frac{1}{\sqrt{\lambda^2+1}}
\begin{bmatrix}
\lambda& 1\\
1& -\lambda
\end{bmatrix}
\right)
=\begin{bmatrix}
1&-\lambda \\
-\lambda& \lambda^2
\end{bmatrix} \;,
\end{equation}
is of the above form. 
Therefore, situations where the kernel of $a(\lambda)$ depends on $\lambda$ are
allowed, which supersedes the results of \cite{LM2, pan_jms} for ultra-parabolic
equations. 
\item[b)] Assumption (b) trivially implies (b'), while
\begin{equation}\label{eq:matrix_2d_notb}
a(\lambda) = 
\begin{pmatrix}
0&0 \\
0&|\lambda|
\end{pmatrix}
\end{equation}
is a simple example which illustrates that the converse does not hold. 
Indeed, $a$ satisfies (b'), but $\begin{pmatrix}0&0 \\ 0&\sqrt{|\lambda|}\end{pmatrix}$ 
is not Lipschitz continuous around zero, which implies that 
a matrix $\sigma$ such that condition (b) is satisfied does not exist.

Since this matrix is singular on the set of zero Lebesgue measure, 
we can apply Corollary \ref{cor1} and still obtain the result. 
However, one can easily generalise \eqref{eq:matrix_2d_notb} to the case where 
the singular set of $\sqrt{a}$ is of positive measure. 
For example, take a Cantor set $C$ on $[0,1]$ of a non-zero measure (so-called 
fat Cantor set) and on each connected component 
$(\alpha,\beta)$ of $[0,1]\setminus C$ define $g(\lambda)=\frac{\beta-\alpha}{2}-|\lambda-\frac{\alpha+\beta}{2}|$
(another possibility could be $g(\lambda)= (\beta-\lambda)(\lambda-\alpha)$).
Then
\begin{equation}\label{eq:matrix_2d_notb_b}
\begin{pmatrix}
0&0 \\
0& g(\lambda)
\end{pmatrix}
\end{equation}
satisfies (b'), but does not satisfy neither (b) nor ($\tilde{{b}}$).
Thus, for this matrix only Corollary \ref{velocity averaging_Gn_in_Lr} is applicable 
among our results. 
\end{enumerate}
\end{example}

\begin{remark}
In Example \ref{ex3}(a) we have seen that smoothness of eigenvectors of $a(\lambda)$
(i.e.~smoothness of $Q(\lambda)$) could help in fulfilling condition (b). 
Let us recall some known results in this direction (\cite[II.6.1-3]{kato}):
\begin{enumerate}
	\item If $ a(\lambda)$ is symmetric and analytic then both eigenvectors and eigenvalues are analytic functions;
	\item If $a(\lambda)$ is symmetric and $\mathrm{C}^1$, then eigenvalues are $\mathrm{C}^1$-functions,
	while eigenvectors need not to be even continuous.
	\end{enumerate} 
Of course, not even item (1), with addition of positive semi-definiteness, is sufficient 
to ensure (b) since we require, in principle, Lipschitz continuity of $\sqrt{a}$. 
\end{remark}

\begin{remark}\label{rem:app-1-parbolic}
Since we are particularly interested in the parabolic case, we refer to \eqref{app-1} as a degenerate parabolic 
equation with discontinuous flux, although it is not necessarily of the parabolic type. More precisely, in the application
to the Cauchy problem for nonlinear degenerate parabolic equation with discontinuous flux (see Section \ref{sec:application})
we shall have
$f(t,\mx,\lambda)=\begin{bmatrix}1\\ \tilde f(t,\mx,\lambda)\end{bmatrix}$ and $a(\lambda)=\begin{bmatrix}0&0\\ 0&\tilde a(\lambda)\end{bmatrix}$.
In order to have that $a$ and $f$ satisfy assumptions (b) and (c) it 
is sufficient to have $\tilde f\in{\rm L}^p(\Omega_{t,\mx}\times S;\Rd)$,
$\tilde a=\tilde{\sigma}^T\tilde{\sigma}$, where 
$\tilde\sigma\in{\rm C}^{0,1}(S;\R^{d\times d})$, and for any 
$K\subset\subset S$
\begin{equation*}
\esssup\limits_{(t,\mx)\in \Omega} \sup\limits_{(\tau,\mxi) \in \mathrm{S}^{d}} 
	\operatorname{meas}\bigl\{\lambda \in K :\, \tau+\langle\tilde f(t,\mx,\lambda)
	\,|\, \mxi\rangle=\langle \tilde a(\lambda)\mxi \,|\,\mxi \rangle=0 \bigr\}=0 \;. 
\end{equation*}
\end{remark}

We shall now briefly explain principles of our approach. 

Since we cannot separate the unknown function $u_n$ from the coefficients in \eqref{app-1}, here we use variants of micro-local defect measures (or H-measures) introduced in now seminal papers by P.~G\'erard \cite{Ger} and L.~Tartar \cite{Tar}.
Besides the velocity averaging results  \cite{Ger, LM2}, the H-measures and similar tools found
applications on existence of traces and solutions to nonlinear evolution equations \cite{1, EM2022, HKMP, pan_jhde}, generalisation of
compensated compactness results to equations with variable coefficients \cite{Ger, Tar}, applications in the control theory \cite{DLR, LZ},   
explicit formulae and bounds in homogenisation  \cite{ALjmaa, tar_book}, etc. 

Moreover, it initiated variety of different generalisations to the original micro-local defect measures which we call here micro-local defect functionals.
We mention parabolic and ultra-parabolic variants of the H-measures \cite{ALjfa, pan_aihp}, H-measures as duals of Bochner spaces \cite{LM2},
H-distributions \cite{AEM, AM, LM3, MM}, micro-local compactness forms \cite{rindler}, one-scale H-measures \cite{AEL, Tar2} etc.

Let us recall the first variant of H-measures \cite{Tar} (introduced at the same time as the micro-local defect measures \cite{Ger}).

\begin{theorem}
\label{tbasic1} If $(u_n)$ is a sequence in ${\rm L}_{loc}^2(\Omega;\R^r)$, $\Omega\subseteq\R^{d}$, such
that $u_n\rightharpoonup 0$ in ${\rm L}_{loc}^2(\Omega;\R^r)$, then there exist a subsequence $(u_{n'})\subset (u_{n})$ and a positive
complex Radon measure $\mu=\{\mu^{jk}\}_{j,k=1,...,r}$ on $\Omega\times {\rm S}^{d-1}$ such that for any
$\varphi_1,\varphi_2\in {\rm C}_c(\Omega)$ and $\psi\in {\rm C}({\rm S}^{d-1})$, it holds
\begin{equation*}
\begin{split}
\lim\limits_{n'\to \infty}& \int_{\Omega}(\varphi_1 u^j_{n'})(\mx) \overline{ {\cal A}_{\bar\psi\left(\frac{\cdot}{|\cdot|}\right)}(\varphi_2 u^k_{n'})(\mx)} \,d\mx=
\langle\mu^{jk},\varphi_1\bar{\varphi}_2\psi \rangle\\
&=\int_{\Omega\times {\rm S}^{d-1}}\varphi_1(\mx)\overline{\varphi_2(\mx)}\psi(\mxi)d\mu^{jk}(\mx,\mxi) \;,
\end{split}
\end{equation*} 
where ${{\cal A}_{\bar\psi\left(\frac{\cdot}{|\cdot|}\right)}}$ is the Fourier multiplier operator with the symbol $\bar\psi(\mxi/|\mxi|)$.
\end{theorem} 

The measure $\mu$ is called the H-measures and, with respect to the dual variable $\mxi$, it is defined on the sphere (since $\mxi/|\mxi| \in {\rm S}^{d-1}$).
\smallskip

It has been proved (see \cite{ALjmaa}) that applying H-measures on differential relations where the ratio of the highest orders of derivatives in each variable is not the same might lead to
unsatisfactory results. This is due to the projection $\mxi\mapsto \mxi/|\mxi|$, since scalings in all variables are the same. 
We can change the scaling and put, for example, $\frac{\mxi}{|(\xi_1,\dots,\xi_k)|+|(\xi_{k+1},\dots,\xi_d)|^2}$ instead of $\mxi/|\mxi|$, but such H-measure will
be able to \emph{see} e.g.~first order derivatives with respect to $(x_1,\dots,x_k)$ and second order derivatives with respect to $(x_{k+1},\dots,x_d)$ (for a parabolic variant, see \cite{ALjfa}).
In other words, no change of the highest order of the equation is permitted. For instance, this means that the matrix $a(\lambda)$ 
in equation \eqref{app-1} must have the rank and the kernel 
(locally) independent of $\lambda$
(see also \cite{pan_aihp}) otherwise, we cannot use the existing theory of the micro-local defect functionals (except in special situations \cite{HKMP}).

This represents a significant confinement since many challenging mathematical questions, especially from a view-point of modelling,
involve equations that change type. In particular, we have in mind degenerate parabolic equations which describe wide range of phenomena 
containing the combined effects of nonlinear convection and degenerate diffusion and which have the form

\begin{align}\label{d-p-1} 
\pa_t u+\Div_{\mx} \mff(t,\mx,u)&=D_\mx^2 \cdot A(u) \,,
\end{align}
where the matrix $A$ is such that the mapping $\R \ni \lambda \mapsto  \langle A(\lambda) \mxi \,|\,\mxi\rangle$ is non-decreasing,
i.e.~that the diffusion matrix $a(\lambda)=A'(\lambda)$ is merely positive semi-definite.
To this end, let us mention \cite{Zuazua}, where one of the first results on the case of degenerate parabolic equations was given (to be more precise, an ultra-parabolic equation was considered there).

Let us remark that in the case when the coefficients in \eqref{d-p-1} are regular, the theory of existence and uniqueness for
appropriate Cauchy problems is well-established (see e.g. \cite{CP, CK, GKM}).
Nevertheless, concrete applications such as flow in porous media  very often occur in highly heterogeneous environment causing rather rough
coefficients in \eqref{d-p-1} (e.g. during $CO_2$ sequestration process \cite{JMN}). Furthermore, even in a simplified situation in which the diffusion is neglected such as
a road traffic with variable number of lanes \cite{BKT2009-multilanes}, the Buckley-Leverett equation in a layered porous medium \cite{AndrCances, Kaas},
and sedimentation applications \cite{9, BurgerEtAl, Die1, Die2}, the flux is as a rule discontinuous.

However, due to obvious technical obstacles, most of the previous literature was dedicated either to homogeneous degenerate parabolic equations or to equations where the flux and
diffusion are regular functions (e.g. \cite{BK, Car, CdB, CP, VH}). We mention  \cite{HKMP} where \eqref{d-p-1} was considered with the flux $\mff(t,\mx,\lambda)$ (merely) continuous with respect to $\lambda$ and belonging to ${\rm L}^p$, $p>2$, with respect to $\mx \in \R^d$. Here, we are able to improve the result by relaxing the assumption of $p$ to $p>1$. More precise explanation can be found in Section 5.

A similar situation is regarding existence of traces. The trace of a function $u$ at $t=0$ is a function $u_0 \in \textrm{L}^1_{loc}(\R^d)$ such that
$$
u(t,\cdot) \to u_0(\cdot) \ \ {\rm as} \ \ t\to 0 \ \ {\rm in} \ \ \mathrm{L}^1_{loc}(\R^d) \,.
$$
One can find several results in the hyperbolic setting 
\cite{pan_jhde, pan_jhdeB, vass}, while degenerate parabolic equations of
the form \eqref{d-p-1} were successfully studied only in some special
situations. For example, in \cite{23} one has the scalar diffusion matrix
$a(\lambda)=\tilde{a}(\lambda)I$, where $I$ is the unitary matrix and
$\tilde{a}\in \mathrm{C}^1(\R)$ is a non-negative function. 
The case of ultra-parabolic matrices was considered in \cite{1}.
In both cases, assumptions were imposed so that the essential problem of $\lambda$-changing degeneracy directions
does not appear. 
To be more descriptive, we note that the matrix $a$ given 
in Example \ref{ex3}(a) is not covered by the results from \cite{1, 23}.
In the current contribution we shall thus provide a result regarding existence of strong traces in the case when the diffusion matrix degenerates in directions which depend on $\lambda$ (thus covering the case of the aforementioned matrix). 
Moreover, we allow that the flux depends explicitly on $\mx$ and it can even be discontinuous,
which is a novelty when compared to a recent result \cite{EM2022}.
However, here we are not able to avoid non-degeneracy assumptions.

We overcome the above mentioned technical problems caused by degenerate diffusion matrices by considering multiplier operators with symbols of the form
\begin{equation}\label{eq:form_of_symbols}
\psi\left(\frac{\mxi}{|\mxi|+\langle a(\lambda) \mxi \,|\, \mxi \rangle}\right) \;, \quad \psi\in {\rm C}(\R^d) \;,
\end{equation}
where the matrix $a$ represents the diffusion matrix in the degenerate parabolic equation \eqref{app-1}.

\smallskip

The paper is organised as follows.

In Section \ref{sec:multipliers} we study symbols of the form \eqref{eq:form_of_symbols} which shall be often used for the Fourier multiplier operators, 
and show two important results concerning their continuity (see Lemma \ref{multiplierlemma1}). 
Section \ref{sec:Hm} is devoted to the construction of adaptive micro-local defect functionals.

In Section 4, we use the results of Section \ref{sec:multipliers} and Section \ref{sec:Hm} to prove the main result of the paper, Theorem \ref{velocity averaging}.

In Section 5, as an application of the velocity averaging result, we show existence of a weak solution to the Cauchy problem of the degenerate advection-diffusion equation with discontinuous flux.
The strategy of the proof is to reduce the degenerate parabolic equation \eqref{d-p-1} to its kinetic counterpart of the form (below, $f=\partial_\lambda\mff$ and $a=A'$):
\begin{equation*}
\begin{split}
\pa_t h(t,\mx,\lambda) &+ \Div_\mx (f(t,\mx,\lambda) h(t,\mx,\lambda) )\\&=\Div_\mx \bigl(\Div_\mx \left( a(\lambda) h(t,\mx,\lambda) \right)\bigr)+ \pa_{\lambda} G(t,\mx,\lambda) + \Div_\mx P(t,\mx,\lambda) \,,
\end{split}
\end{equation*} and then to use the velocity averaging results. 

In Section 6, we provide another application of the velocity averaging result by proving that any bounded
quasi-solution to \eqref{d-p-1} (see Definition \ref{def25}) admits the strong trace at $t=0$ under the non-degeneracy condition:
\begin{equation*}
 \sup\limits_{\mxi \in \mathrm{S}^{d-1}} 
	\operatorname{meas}\bigl\{\lambda \in K :\, \langle a(\lambda)\mxi \,|\,\mxi \rangle=0 \bigr\}=0 \,. 
\end{equation*}

\section{Results on Fourier integral operators}\label{sec:multipliers}

Let $a:S\to \R^{d\times d}$, $S\subseteq \R$, be a Borel measurable matrix function such that for a.e.~$\lambda\in S$ matrix $a(\lambda)$ is 
symmetric and positive semi-definite, i.e.~$a(\lambda)^T=a(\lambda)$ and $\langle a(\lambda)\mxi \,|\,\mxi\rangle \geq 0$, $\mxi\in\Rd$. 
Further on, we define
\begin{equation}
\label{Pi}
\pi_P(\mxi,\lambda) := \frac{\mxi}{|\mxi|+\langle a(\lambda)\mxi \,|\,\mxi \rangle} \ , \quad (\mxi,\lambda)\in\Rd\!\setminus\!\{0\}\times S \;.
\end{equation} 
As $a(\lambda)$ is positive semi-definite, we have $\pi_P(\R^{d}\!\setminus\!\{0\}\times S)\subseteq B[0,1]\setminus\{0\}$, 
where $B[0,1]$ denotes the unit closed ball in $\Rd$.
Moreover, it is not difficult to show that (for a.e.~$\lambda\in S$)
\begin{equation}\label{eq:slika_od_pi}
	\mathop{Cl}\pi_P(\Rd\!\setminus\!\{0\},\lambda) = 
	\left\{\begin{array}{lcr}
	\Sdmj &:& a(\lambda)=0 \\
	B[0,1] &:& a(\lambda)\ne 0
	\end{array}\right. \;,
\end{equation}
where $\mathop{Cl} A$ denotes the closure of $A\subseteq\Rd$.

If $a(\lambda)=0$, $\pi_P(\cdot,\lambda)$ is the projection of $\Rd\setminus\{0\}$ to the unit sphere along the rays through the origin.
In general $\pi_P(\cdot,\lambda)$ is not a projection since $a(\lambda)\neq0$ implies
$\pi_P(\pi_P(\cdot,\lambda),\lambda)\neq\pi_P(\cdot,\lambda)$. However, for simplicity, in the text we shall often address $\pi_P$ as a projection.

In this paper, we are interested in symbols of Fourier multiplier operators
of the form 
$$
\mxi\mapsto\bar\psi(\pi_P(\mxi,\lambda),\lambda) \,,
$$
where $\lambda\in S$ is fixed, $\psi\in{\rm L}^\infty(S;{\rm C}(B[0,1]))$, and $\pi_P$ is as above. 
Here $\bar z$ denotes the complex conjugate of complex number $z$.

Of course, $\psi\in{\rm L}^\infty(S;{\rm C}(B[0,1]))$ is sufficient to have 
that the Fourier multiplier operator is bounded on $\mathrm{L}^2(\Rd)$, with the
norm independent on $\lambda$. However, we shall need such a result on an arbitrary
$\mathrm{L}^p$, for which we need some additional regularity of $\psi$ with respect
to $\mxi$. 
More precisely, we shall first obtain that for a.e.~$\lambda$ and for any $p\in\oi1\infty$ the operator ${\cal A}_{\bar \psi(\pi_P(\cdot,\lambda),\lambda)}$
is bounded on ${\rm L}^p(\Rd)$, with the norm independent of $\lambda$ (Lemma \ref{multiplierlemma1}). Finally, we show that commutators of the Fourier multiplier operators 
and operators of multiplication map weakly converging sequences to strongly converging in a certain sense (Corollary \ref{cor:komutacija}).

In order to prove the ${\rm L}^p$ boundedness, we use the following corollary of the Marcinkiewicz multiplier theorem \cite[Corollary 5.2.5]{Gra}:

\begin{theorem}
	\label{m1} Suppose that  $\psi\in \pC{d}{\R^d\setminus\cup_{j=1}^d\{\xi_j= 0\}}$ is a bounded function such that for some
	constant $C>0$ it holds
	\begin{equation}
	\label{c-mar} |\mxi^{\malpha} \partial^{\malpha}
	\psi(\mxi)|\leq C,\ \
	\mxi\in \R^d\backslash \cup_{j=1}^d\{\xi_j= 0\}
	\end{equation}
	for  every multi-index
	$\malpha=(\alpha_1,\dots,\alpha_d) \in {\N}_0^d$
	such that
	$|\malpha|=\alpha_1+\alpha_2+\dots+\alpha_d
	\leq d$. Then $\psi$ is an ${\rm L}^p$-multiplier for
	any $p\in \oi 1\infty$, and the operator norm of ${\cal A}_\psi$ equals 
	$C_{d,p} C$, where
	$C_{d,p}$ depends only on $p$ and $d$.
\end{theorem}

Before proceeding with the verification of the assumptions of the previous theorem, let us recall some 
well known results from matrix analysis and at the same time fix our notations. 

As $a(\lambda)$ is a positive semi-definite symmetric matrix of order $d$, there exist orthogonal matrix $Q(\lambda)$
and diagonal matrix $\Lambda(\lambda) = {\rm diag}(\kappa_1(\lambda),\kappa_2(\lambda),\dots,\kappa_d(\lambda))$, 
containing (non-negative) eigenvalues of $a(\lambda)$, such that the following eigendecomposition holds: 
\begin{equation}\label{eq:diagonalisation}
	a(\lambda)=Q(\lambda)^T\Lambda(\lambda) Q(\lambda) \;.
\end{equation}
If 
\begin{equation}\label{eq:sigma}
a(\lambda)=\sigma(\lambda)^T \sigma(\lambda)
\end{equation}
holds for $a$ given by \eqref{eq:diagonalisation}, then $\sigma$
is necessarily of the form
\begin{equation}\label{eq:sigma_b}
\sigma(\lambda)= \widetilde{Q}(\lambda)\sqrt{\Lambda(\lambda)}Q(\lambda) \,,
\end{equation}
where $\sqrt{\Lambda(\lambda)} = {\rm diag}(\sqrt{\kappa_1(\lambda)},\sqrt{\kappa_2(\lambda)},\dots,\sqrt{\kappa_d(\lambda)})$ and $\widetilde{Q}(\lambda)$ is an orthogonal matrix.
For example, when $\tilde Q(\lambda)$ is the identity matrix, then we 
just have $\sigma(\lambda)=\sqrt{\Lambda(\lambda)}Q(\lambda)$.

It is important to notice that
$$
\pi_P\bigl(Q(\lambda)^T\mxi,\lambda\bigr) = \frac{Q(\lambda)^T\mxi}{|Q(\lambda)^T\mxi|+\langle Q(\lambda)a(\lambda)Q(\lambda)^T \mxi \,|\, \mxi \rangle} 
	= \frac{Q(\lambda)^T\mxi}{|\mxi|+ \sum_{j=1}^d \kappa_j(\lambda)\xi_j^2 }\,,
$$
where we have used that $Q(\lambda)^T$ preserves the length of vectors. Hence, with the orthogonal change of variables we will 
manage to reduce the problem to the case of diagonal matrix $a$.

\begin{lemma}
	\label{multiplierlemma1}
	Let $a:S\to \R^{d\times d}$, $S\subseteq \R$ open, be a Borel measurable matrix function such that for a.e.~$\lambda\in S$ matrix $a(\lambda)$ is 
	symmetric and positive semi-definite, and let $\psi\in{\rm L}^\infty(S;{\rm C}^d(B[0,1]))$. 
	
	Then for a.e.~$\lambda\in S$ and any $p\in\oi1\infty$, function $\bar\psi(\pi_P(\cdot,\lambda),\lambda)$, where $\pi_P$
	is given by \eqref{Pi}, is an ${\rm L}^p$-Fourier multiplier and the ${\rm L}^p$-norm of the corresponding Fourier multiplier operator 
	is independent of $\lambda$.
\end{lemma}
\begin{proof}
	Since the space of ${\rm L}^p$-Fourier multipliers 
	is invariant under orthogonal change of variables 
	\cite[Proposition 2.5.14]{Gra}
	(see also Lemma \ref{vaznalema} below) and the corresponding 
	norms coincide, applying $\mxi\mapsto Q(\lambda)^T \mxi$, 
	where $Q(\lambda)$ is given in \eqref{eq:diagonalisation}, 
	it is sufficient to study 
	$\mxi\mapsto \psi\bigl(\pi_P(Q(\lambda)^T\,\cdot,\lambda),\lambda\bigr)$.
	We shall apply the Marcinkiewicz multiplier theorem 
	(Theorem \ref{m1}) on this function.
	
	Since $\pi_P(\R^{d}\!\setminus\!\{0\}\times S)\subseteq B[0,1]$ 
	and $\psi\in{\rm L}^\infty(S;{\rm C}^d(B[0,1]))$, the functions 
	$$
	(\mxi,\lambda)\mapsto(\partial^\malpha_\mxi\bar\psi)\bigl(\pi_P(Q(\lambda)^T\mxi,\lambda),\lambda\bigr)
	$$ 
	are bounded on $\Rd\!\setminus\!\{0\}\times S$  for all $\malpha\in\mathbb{N}_0^d$, $|\malpha|\leq d$.
	Therefore, by the generalised chain rule formula (known as the Fa\'a di Bruno formula; see e.g.~\cite{FaadiBruno}) it is enough to infer that \eqref{c-mar} is satisfied 
	for each component of $\pi_P(Q(\lambda)^T\,\cdot\,,\lambda)$,
	with constant $C$ independent of $\lambda$.  
	
	Furthermore, since the Riesz transform of order 1
	satisfies \eqref{c-mar}, $Q(\lambda)$ is orthogonal and
	$$
	\bigl(\pi_P(Q(\lambda)^T\mxi,\lambda)\bigr)_j=
	\frac{(Q(\lambda)^T\xi)_j}{|\mxi|+\langle \Lambda(\lambda)\mxi \,|\,\mxi \rangle}=\frac{(Q(\lambda)^T\xi)_j}{|\mxi|}\frac{|\mxi|}{|\mxi|+\langle \Lambda(\lambda)\mxi \,|\,\mxi \rangle} \,,
	$$
	by the Leibniz rule it is sufficient to check \eqref{c-mar}
	for 
	$$
	\mxi \mapsto \frac{|\mxi|}{|\mxi|+\langle \Lambda(\lambda)\mxi \,|\,\mxi \rangle} \,.
	$$ 
	The claim follows by Lemma \ref{TVRDNJA} below.
	\end{proof}

The proof of the following lemma we leave for the Appendix. 

\begin{lemma}
	\label{TVRDNJA}
	For any $\mkappa=(\kappa_1,\kappa_2,\dots,\kappa_d) \in [0,\infty)^d$, 
	$m\in\{1,2,\dots,d\}$, $s\in [0,\infty)$, and $p\in\oi1\infty$, 
	functions $f^s$ and $g^s$ ($s$ is an exponent), where $f,g:\Rd\to\R$ are given by 
	$$
	f(\mxi)=\frac{|\mxi|}{|\mxi|+\sum\limits_{j=1}^d \kappa_j \xi_j^2}
		\qquad\hbox{and}\qquad 
		g(\mxi)=\frac{\kappa_m\xi_m^2}{|\mxi|+\sum\limits_{j=1}^d \kappa_j \xi_j^2} \ ,
	$$
	are ${\rm L}^p$-Fourier multipliers and the norm of the corresponding 
	Fourier multiplier operators depends only on $d$, $s$ and $p$,
	i.e.~it is independent of $\mkappa$. 
	Moreover, $f^s$ and $g^s$ satisfy the Marcinkiewicz  condition \eqref{c-mar}
	with a constant $C$ independent of $\mkappa$.
\end{lemma}

Our next goal is to study the Fourier multiplier operator associated to the symbol $\mxi \mapsto \pa_\lambda \frac{1}{|\mxi|+\langle a(\lambda)\mxi \,|\,\mxi \rangle}$.
For a smooth $a$, we have
\begin{equation}\label{eq:symbol_with_derivative_in_lambda}
\pa_\lambda \frac{1}{|\mxi|+\langle a(\lambda)\mxi \,|\,\mxi \rangle} 
	= \frac{-\langle a'(\lambda)\mxi \,|\, \mxi\rangle}{\bigl(|\mxi|+\langle a(\lambda)\mxi \,|\,\mxi \rangle\bigr)^2}
	= \psi(\pi_P(\mxi,\lambda),\lambda) \,,
\end{equation}
where $\psi(\mxi,\lambda)=-\langle a'(\lambda)\mxi \,|\, \mxi\rangle$. Thus, by Lemma \ref{multiplierlemma1},
the Fourier multiplier operator is ${\rm L}^p$-bounded, $p\in\oi1\infty$, uniformly in $\lambda$, if $a'$ exists (almost everywhere) and it is bounded. 
However, we need that this operator has a smoothing property.

Let us additionally assume that
$\sigma$ given by \eqref{eq:sigma} is Lipschitz continuous. 
Since 
$$
\langle a(\lambda)\mxi\,|\,\mxi\rangle = |\sigma(\lambda)\mxi|^2
	= \sum_{j=1}^d \bigl(\sigma(\lambda)\mxi\bigr)_j^2 \;,
$$
we have (for almost every $\lambda\in S$)
\begin{equation*}
\langle a'(\lambda)\mxi\,|\,\mxi\rangle = \frac{d}{d\lambda} 
	\langle a(\lambda)\mxi\,|\,\mxi\rangle
	= 2\sum_{j=1}^d \bigl(\sigma(\lambda)\mxi\bigr)_j 
	\bigl(\sigma'(\lambda)\mxi\bigr)_j \;.
\end{equation*}
Thus, symbol \eqref{eq:symbol_with_derivative_in_lambda} can be rewritten 
as
\begin{equation}\label{eq:symbol_ader_sigma}
-2\sum_{j=1}^d \frac{1}{\sqrt{|\mxi|+\langle a(\lambda)\mxi\,|\,\mxi\rangle}}
	\frac{\bigl(\sigma(\lambda)\mxi\bigr)_j}
	{\sqrt{|\mxi|+\langle a(\lambda)\mxi\,|\,\mxi\rangle}}
	\frac{\bigl(\sigma'(\lambda)\mxi\bigr)_j}
	{|\mxi|+\langle a(\lambda)\mxi\,|\,\mxi\rangle} \;,
\end{equation}
and the term $\frac{1}{\sqrt{|\mxi|+\langle a(\lambda)\mxi\,|\,\mxi\rangle}}$
will provide a smoothing property of the half derivative.

\begin{lemma}\label{multiplierlemma2}
	In addition to the assumptions in Lemma \ref{multiplierlemma1}, assume that
	there exists a Lipschitz continuous matrix function $\sigma:S\to\R^{d\times d}$ 
	such that \eqref{eq:sigma} holds. 
	Then, for a.e.~$\lambda \in S$ and any $p\in\oi1\infty$ the operator 
	${\cal A}_{\pa_\lambda \frac{1}{|\mxi|+\langle a(\lambda)\mxi \,|\,\mxi \rangle}}:
	{\rm L}^p(\R^d) \to {\rm W}^{\frac{1}{2},p}(\R^d)$ is bounded uniformly with respect to $\lambda\in S$.
\end{lemma}
\begin{proof}
	Since $a$ is a Lipschitz map, $a'$ exists almost everywhere and it is bounded. 
	Thus, $\mxi\mapsto-\langle a'(\lambda)\mxi \,|\,\mxi\rangle$ satisfies assumptions of Lemma \ref{multiplierlemma1},
	and by \eqref{eq:symbol_with_derivative_in_lambda}, for any $p\in\oi1\infty$,
	the operator
	\begin{equation*}
	{\cal A}_{\pa_\lambda \frac{1}{|\mxi|+\langle a(\lambda)\mxi \,|\,\mxi \rangle}}: {\rm L}^p(\R^d)\to {\rm L}^p(\R^d)
	\end{equation*}
	is uniformly bounded in $\lambda$.
	
	To prove that ${\cal A}_{\pa_\lambda \frac{1}{|\mxi|+\langle a(\lambda)\mxi \,|\,\mxi \rangle}}$ possesses 
	a smoothing property, we need to prove that derivatives (with respect to $\mx$) of the operator are 
	${\rm L}^p\to {\rm L}^p$ bounded uniformly in $\lambda$:
	\begin{equation*}
	\pa_{x_k}^{\frac12}{\cal A}_{\pa_\lambda \frac{1}{|\mxi|+\langle a(\lambda)\mxi \,|\,\mxi \rangle}}: {\rm L}^p(\R^d)\to {\rm L}^p(\R^d), \ \ k=1,\dots,d \;,
	\end{equation*} 
	where by $\pa_{x_k}^{\frac12}$ we have denoted the operator 
	${\cal A}_{(2\pi i \xi_k)^{1/2}}$.
	The symbol of the latter operator is 
	$$
	\frac{-(2\pi i \xi_k)^\frac{1}{2} \langle a'(\lambda)\mxi \,|\,\mxi \rangle}{\left(|\mxi|+\langle a(\lambda)\mxi \,|\,\mxi \rangle\right)^2} \;,
	$$ 
	which by \eqref{eq:symbol_ader_sigma} can be rewritten as
	\begin{equation*}
	2\sum_{j=1}^d \frac{-(2\pi i\xi_k)^\frac12}{\sqrt{|\mxi|+\langle a(\lambda)\mxi\,|\,\mxi\rangle}}
	\frac{\bigl(\sigma(\lambda)\mxi\bigr)_j}
	{\sqrt{|\mxi|+\langle a(\lambda)\mxi\,|\,\mxi\rangle}}
	\frac{\bigl(\sigma'(\lambda)\mxi\bigr)_j}
	{|\mxi|+\langle a(\lambda)\mxi\,|\,\mxi\rangle} \;.
	\end{equation*}
	The space of ${\rm L}^p$-Fourier multipliers is an algebra
	\cite[Proposition 2.5.13]{Gra}, hence we can study each 
	factor separately.

	Since $\sigma'$ is bounded, by Lemma \ref{multiplierlemma1} 
	for a.e.~$\lambda\in S$ and any $p\in (1,\infty)$ 
	$$
	\mxi\mapsto\frac{\bigl(\sigma'(\lambda)\mxi\bigr)_j}
		{|\mxi|+\langle a(\lambda)\mxi\,|\,\mxi\rangle}
	$$
	is an $\mathrm{L}^p$-multiplier with the norm independent of $\lambda$.
	The same holds for
	$$
	\mxi\mapsto \frac{-(2\pi i\xi_k)^\frac12}{\sqrt{|\mxi|+\langle a(\lambda)\mxi\,|\,\mxi\rangle}}
	$$
	by Lemma \ref{TVRDNJA} (see the result for $\sqrt f$). 
	
	Furthermore, we have $\sigma(\lambda)=\widetilde{Q}(\lambda)\sqrt{\Lambda(\lambda)}Q(\lambda)$ 
	(see \eqref{eq:sigma_b}). Thus, as the space of ${\rm L}^p$-Fourier multipliers is invariant under 
	orthogonal change of variables \cite[Proposition 2.5.14]{Gra}
	(see also Lemma \ref{vaznalema} below) and the corresponding norms coincide, 
	applying $\mxi\mapsto Q(\lambda)^T \mxi$ and using $|\widetilde{Q}(\lambda)|=1$, it is left to study
	$$
	\mxi\mapsto 	\frac{\sqrt{\kappa_j(\lambda)}\xi_j}
	{\sqrt{|\mxi|+\sum_{l=1}^d \kappa_l(\lambda)\xi_l^2}} \;.
	$$
	Finally, by Lemma \ref{TVRDNJA} above (see the result for $\sqrt{g}$) we have that this mapping is 
	an $\mathrm{L}^p$ multiplier with the $\mathrm{L}^p$-norm of the corresponding
	Fourier multiplier operator independent of functions $\kappa_l$, and thus of
	$\lambda$. 
	\end{proof}

It is by now a classical result that if we have an ${\rm L}^p$-multiplier, then the composition with an orthogonal
matrix is also an ${\rm L}^p$-multiplier with the same operator norm (see Proposition 2.5.14 in \cite{Gra}).
Now, we will show something very similar when we have a regular change of variables. 
The result is well known but we include it here for completeness.

\begin{lemma}\label{vaznalema}
Let $\psi\in{\rm L}^\infty(\Rd)$. If there exists a regular real constant matrix $M$ and $p\in\oi 1\infty$
such that $\psi(M^{-1}\cdot)$ is an ${\rm L}^p$-multiplier, then $\psi$ is also an ${\rm L}^p$-multiplier
and $\|{\cal A}_\psi\|_{{\cal L}({\rm L}^p(\Rd))}=\|{\cal A}_{\psi(M^{-1}\cdot)}\|_{{\cal L}({\rm L}^p(\Rd))}$.
\end{lemma}
\begin{proof}
Let us denote by $A:=\|{\cal A}_{\psi(M^{-1}\cdot)}\|_{{\cal L}({\rm L}^p)}$ the operator norm, and 
by $J:=|\det{M}|>0$ the Jacobian. 

Take $\varphi\in{\rm C}^\infty_c(\Rd)$ and for an arbitrary $u\in{\rm C}^\infty_c(\Rd)$, consider the following:

\begin{align*}
\int_\Rd\varphi(\mx)\overline{{\cal A}_\psi(u)(\mx)}\; d\mx &= \int_\Rd \widehat{\varphi}(\mxi)\overline{\psi(\mxi)\widehat{u}(\mxi)} \; d\mxi\\
&=J^{-1}\int_\Rd \widehat{\varphi}(M^{-1}\meta)\overline{\psi(M^{-1}\meta)\widehat{u}(M^{-1}\meta)} \; d\meta \;,
\end{align*}
where we have used Plancherel's theorem in the first equality
and the regular change of variables $\meta=M\mxi$ in the second one.
Furthermore, we have

\begin{align*}
\widehat{\varphi}(M^{-1}\meta) &= \int_\Rd e^{-2\pi i\mx\cdot M^{-1}\meta}\varphi(\mx)\; d\mx
= \int_\Rd e^{-2\pi i M^{-T}\mx\cdot \meta} \varphi(\mx)\; d\mx\\
&= J\int_\Rd e^{-2\pi i\my\cdot \meta}\varphi(M^T\my)\; d\my = \widehat{\varphi(M^T\cdot)}(\meta)\, J \;,
\end{align*}
where we have used the change of variables $\my = M^{-T}\mx$ in the third equality.
After applying Plancherel's theorem once more, we get

\begin{align*}
\left| \int_\Rd \varphi(\mx) \overline{{\cal A}_\psi(u)(\mx)}\; d\mx\right| 
	&= J^{-1}\left|\int_\Rd \widehat{\varphi}(M^{-1}\meta)\overline{\psi(M^{-1}\meta)\widehat{u}(M^{-1}\meta)} \; d\meta\right|\\
&=J\left|\int_\Rd \widehat{\varphi(M^T\cdot)}(\meta) \overline{\psi(M^{-1}\meta)\widehat{u(M^{T}\cdot)}(\meta)} \; d\meta \right|\\
&=J\left|\int_\Rd \varphi(M^{T}\my) \overline{{\cal A}_{\psi(M^{-1}\cdot)}(u(M^{T}\cdot))(\my)} \; d\my \right|\\
&\leq J A \|\varphi(M^T\cdot)\|_{{\rm L}^{p'}(\Rd)} \|u(M^T\cdot)\|_{{\rm L}^p(\Rd)}\\
&\leq A \|\varphi\|_{{\rm L}^{p'}(\Rd)} \|u\|_{{\rm L}^p(\Rd)} \;,
\end{align*}
where we have used the H\"older inequality ($1/p+1/p'=1$) and the boundedness of ${\cal A}_{\psi(M^{-1}\cdot)}$ in 
${\rm L}^p(\Rd)$, while the last inequality follows by the fact that the composition with $M^T$ scales 
the ${\rm L}^p$ norm of the function by factor $1/\sqrt[p]{|\det{M}|}$.

From here we conclude that ${\cal A}_\psi(u)$ is a continuous linear functional defined on a dense subset of ${\rm L}^{p'}(\Rd)$.
Thus, by the density argument, we can uniquely extend it
to a linear functional on the whole ${\rm L}^{p'}(\Rd)$, implying that 
${\cal A}_\psi(u)\in{\rm L}^p(\Rd)$ with the following bound:
$$
\| {\cal A}_\psi(u) \|_{{\rm L}^p(\Rd)} \leq A \| u \|_{{\rm L}^p(\Rd)} \,.
$$
The lemma follows for an arbitrary $u\in{\rm L}^p(\Rd)$ once we again use the same density argument as above. 
\end{proof}

In the remaining part of the section we study commutators of Fourier multipliers and operator of multiplications. 

In this section we will need a variant of the First commutation lemma which is given in \cite[Lemma 1]{1stcommlemm}
(see also Remark 2 in the mentioned reference).

\begin{theorem}\label{izclankaokomutatoru}
	Let $(v_n)$ be a bounded, uniformly compactly supported sequence in ${\rm L}^\infty(\Rd)$, converging to 0 in the sense of
	distributions. Let $\psi\in{\rm C}^d(\Rd\!\setminus\!\{0\})\cap{\rm L}^\infty(\Rd)$ be an ${\rm L}^p$-multiplier, 
	$p\in\oi 1\infty$, which satisfies
	\begin{equation}\label{eq:compactness_cond}
	\lim_{|\mxi|\to\infty}\sup_{|{\bf h}|\leq1}\left| \psi(\mxi+{\bf h}) - \psi(\mxi) \right|=0 \;.
	\end{equation}
	
	Then for any $b\in{\rm L}^\infty(\Rd)$ and $r\in\oi1\infty$ the following holds:
	$$
	b{\cal A}_{\psi}(v_n)-{\cal A}_{\psi}(b v_n) \longrightarrow 0 \quad
	\hbox{strongly in} \quad {\rm L}^r_{loc}(\Rd) \;.
	$$
\end{theorem}

In the following lemma we show that symbols of the form \eqref{eq:form_of_symbols} satisfy condition \eqref{eq:compactness_cond}.

\begin{lemma}
Under assumptions of Lemma \ref{multiplierlemma1}, for a.e.~$\lambda\in S$, the function $\bar{\psi}(\pi_P(\cdot,\lambda),\lambda)$ satisfies \eqref{eq:compactness_cond}.
\end{lemma}
\begin{proof}
Since $\bar{\psi}(\cdot,\lambda)$ is uniformly continuous on $B[0,1]$, it is sufficient to prove 
that vector valued function $\pi_P(\cdot,\lambda)$ satisfies 
\eqref{eq:compactness_cond}. Moreover, since
\eqref{eq:compactness_cond} is invariant under orthogonal change of coordinates, 
it is sufficient to study $\pi_P(Q(\lambda)^T\cdot,\lambda)$, where orthogonal matrix $Q(\lambda)$
is given by \eqref{eq:diagonalisation}.

For an arbitrary $|\mathbf{h}|\leq 1$ let us estimate 
$|\pi_P\bigl(Q(\lambda)^T\mxi,\lambda\bigr)-\pi_P\bigl(Q(\lambda)^T(\mxi+\mathbf{h}),\lambda\bigr)|$. To make the calculations easier to read, we omit the 
fixed parameter $\lambda$. Thus, we have 

\begin{align*}
\Bigl|& \frac{Q^T\mxi}{|\mxi|+\langle \Lambda\mxi\,|\,\mxi\rangle}
	-\frac{Q^T(\mxi+\mathbf{h})}{|\mxi+\mathbf{h}|+\langle \Lambda
	(\mxi+\mathbf{h})\,|\,\mxi+\mathbf{h}\rangle}\Bigr| \\
&\qquad \leq \frac{|Q^T\mxi - Q^T(\mxi+\mathbf{h})|}{|\mxi|+\langle \Lambda\mxi\,|\,\mxi\rangle}
	+|Q^T(\mxi+\mathbf{h})|\Bigl|\frac{1}{|\mxi|+\langle \Lambda\mxi\,|\,\mxi\rangle}
	-\frac{1}{|\mxi+\mathbf{h}|+\langle \Lambda(\mxi+\mathbf{h})\,|\,\mxi+\mathbf{h}\rangle}\Bigr| \\
&\qquad \leq \frac{1}{|\mxi|} + |\mxi+\mathbf{h}|\frac{\bigl||\mxi+\mathbf{h}|
	-|\mxi|\bigr|+\bigl|\langle\Lambda(\mxi+\mathbf{h})\,|\,\mxi+\mathbf{h}\rangle
	-\langle\Lambda\mxi\,|\,\mxi\rangle\bigr|}{\bigl(|\mxi|+\langle \Lambda\mxi\,|\,\mxi\rangle\bigr)
	\bigl(|\mxi+\mathbf{h}|+\langle \Lambda(\mxi+\mathbf{h})\,|\,\mxi+\mathbf{h}\rangle\bigr)} \\
&\qquad \leq \frac{2}{|\mxi|} + \frac{\langle\Lambda\mathbf{h}\,|\,\mathbf{h}\rangle + 2|\langle\Lambda\mxi\,|\,\mathbf{h}\rangle|}{|\mxi|+\langle \Lambda\mxi\,|\,\mxi\rangle} \\
&\qquad \leq \frac{2+|\Lambda|}{|\mxi|} + \frac{2\sqrt{|\Lambda|}\sqrt{\langle\Lambda\mxi\,|\,\mxi\rangle}}{|\mxi|+\langle\Lambda\mxi\,|\,\mxi\rangle} \,,
\end{align*}
where in the last line we have used the Cauchy-Bunjakovskij-Schwartz inequality for 
semi-definite scalar product $(\mxi,\meta)\mapsto\langle\Lambda\mxi\,|\,\meta\rangle$.

Now, for $R>0$ and $|\mxi|>R$ we have
$$
\Bigl|\pi_P\bigl(Q(\lambda)^T\mxi,\lambda\bigr)-\pi_P\bigl(Q(\lambda)^T(\mxi+\mathbf{h}),\lambda\bigr)\Bigr| \leq \frac{2+|\Lambda(\lambda)|}{R} +
\frac{2\sqrt{|\Lambda(\lambda)|}}{\sqrt{R}} \,,
$$
implying the claim.
Indeed, if $\langle\Lambda(\lambda)\mxi\,|\,\mxi\rangle \leq R$, we have
$$
\frac{\sqrt{\langle\Lambda\mxi\,|\,\mxi\rangle}}{|\mxi|+\langle\Lambda\mxi\,|\,\mxi\rangle} \leq \frac{\sqrt R}{R} = \frac{1}{\sqrt R} \,,
$$
while for $\langle\Lambda(\lambda)\mxi\,|\,\mxi\rangle > R$ we get
$$
\frac{\sqrt{\langle\Lambda\mxi\,|\,\mxi\rangle}}{|\mxi|+\langle\Lambda\mxi\,|\,\mxi\rangle} < \frac{1}{\sqrt{\langle\Lambda(\lambda)\mxi\,|\,\mxi\rangle}}
 < \frac{1}{\sqrt R} \,.
$$
\end{proof}

By the previous lemma and Lemma \ref{multiplierlemma1}, all assumptions of Theorem \ref{izclankaokomutatoru}
are satisfied. Hence the following corollary holds.

\begin{corollary}\label{cor:komutacija}
	Let $(v_n)$ be a bounded, uniformly compactly supported sequence in ${\rm L}^\infty(\Rd)$, converging to 0 in the sense of
	distributions, and let $\psi$ and $a$ be as in Lemma \ref{multiplierlemma1}.
	
	Then for any $b\in{\rm L}^\infty(\Rd)$, $r\in\oi1\infty$ and a.e.~$\lambda\in S$ the following holds:
	$$
	b{\cal A}_{\bar\psi(\pi_P(\cdot,\lambda),\lambda)}(v_n)
		-{\cal A}_{\bar\psi(\pi_P(\cdot,\lambda),\lambda)}(b v_n) \longrightarrow 0 
		\quad \hbox{strongly in} \quad {\rm L}^r_{loc}(\Rd) \;.
	$$
\end{corollary}

\bigskip

\section{Adaptive micro-local defect functionals}\label{sec:Hm}

In what follows we will have an uniformly compactly supported sequence $(u_n(\mx,\lambda))$.
It means that there exists a bounded open subset $\Omega \times S \subseteq\R^{d+1}$ of finite 
Lebesgue measure such that supports of all functions $u_n$ are contained in it. 
Let us take one such $\Omega\times S$ and fix it.


Now, we need to introduce the space on which we shall define the appropriate micro-local defect functional. 
The space will be adapted to the considered equation \eqref{app-1}. For $p\in(1,\infty)$, we define 
\begin{equation*}
\widetilde{W}^p_\Pi(\Omega,S)=\Bigg\{\sum\limits_{j=1}^k \varphi_j(\mx)\psi_j(\mxi,\lambda): \, k\in\N \,, 
	\ \varphi_j\in {\rm L}^p(\Omega) \,, \, \psi_j\in {\rm C}_c(B[0,1]\times S) \,, \, j=1,\dots,k \Bigg\} \;,
\end{equation*} 
where $B[0,1]$ is the unit closed ball in $\R^{d}$.
We denote for $\Psi=\Psi(\mx,\mxi,\lambda)\in\widetilde{W}^p_\Pi(\Omega,S)$
\begin{align*}
\|\Psi\|_{W^p_{\Pi}}=\left( \int_{\Omega} \left[ \sup_{\mxi\in \R^{d}\setminus\{0\}}
\left( \int_{S} |\Psi\bigl(\mx,\pi_P(\mxi,\lambda),\lambda\bigr)|^2\; d\lambda \right)^{1/2} \right]^{p} d\mx\right)^{1/p} \;,
\end{align*}
where $\pi_P$ is given by \eqref{Pi}. Due to \eqref{eq:slika_od_pi}, 
for $a\not\equiv 0$, this map represents a norm on $\widetilde{W}^p_\Pi(\Omega,S)$.
On the other hand, for $a\equiv 0$ it is only a seminorm, so one needs to consider 
the quotient space by its kernel, or, equivalently, just replace $B[0,1]$ by $\Sdmj$ in the 
definition of $\widetilde{W}^p_\Pi(\Omega,S)$. 

Finally, we introduce the space $W^p_{\Pi}(\Omega,S)$ as the completion of $\widetilde{W}^p_\Pi(\Omega,S)$ 
with respect to the norm $\|\cdot\|_{W^p_{\Pi}}$. It is easy to see that $W^p_{\Pi}(\Omega,S)$ coincides
with the Bochner space ${\rm L}^p(\Omega; X)$, where $X$ is the completion of ${\rm C}(B[0,1];{\rm L}^2(S))$
equipped with the norm
$$
{\rm C}(B[0,1];{\rm L}^2(S)) \ni \psi \mapsto \sup_{\mxi\in\Rd\setminus\{0\}} 
	\left(\int_{S} |\psi\bigl(\pi_P(\mxi,\lambda),\lambda\bigr)|^2 \,d\lambda\right)^{1/2} \;.
$$
Moreover, the space $\mathrm{C}_0(\Omega\times B[0,1]\times S)$,
equipped by the standard (supremum) topology, 
is for any $p\in (1,\infty)$ dense in $W^p_{\Pi}(\Omega,S)$, so 
continuous linear functionals on $W^p_{\Pi}(\Omega,S)$ are in fact 
bounded Radon measures on $\Omega\times B[0,1]\times S$. 

In the following theorem we construct one such functional which will 
play an important role in the proof of the velocity averaging result.
The construction is based on the Banach-Alaoglu-Bourbaki theorem,
which applies on $W^p_{\Pi}(\Omega,S)$ as it is clearly a (separable) Banach space.

\begin{theorem}\label{bilinearboundedness}
Let $(u_n(\mx,\lambda))$ be bounded in ${\rm L}^q(\R^d\times\R)$, for some $q> 2$, and 
uniformly compactly supported on $\Omega \times S \subset\subset \R^d\times \R$.
Let $(v_n(\mx))$ be an uniformly compactly supported sequence on $\Omega$ weakly-$\star$ converging to zero in
${\rm L}^\infty(\R^d)$, and let $a: S\to \R^{d\times d}$ be a Borel measurable matrix function such that for 
a.e.~$\lambda\in S$ matrix $a(\lambda)$ is symmetric and positive semi-definite.

Then for $p= \frac{2q}{q-2}$ there exists a subsequence (not relabelled) 
and a continuous functional $\mu$ on $W^{p}_\Pi(\Omega,S)$ such that for every
$\varphi\in {\rm L}^{p}(\Omega)$ and $\psi\in{\rm C}_c(B[0,1]\times S)$ it holds
\begin{equation}
\label{mu-repr-h}
\mu(\varphi \psi) = \lim_{n\to\infty}\int_{\Omega\times S}\varphi(\mx) u_n(\mx,\lambda)\overline{{\cal A}_{\bar\psi(\pi_P(\cdot,\lambda),\lambda)}(v_n)(\mx)}\; d\mx \,d\lambda \,.
\end{equation} 
Furthermore, the bound of the functional $\mu$ is $C_{u,q,2}C_{v,2}$,
where $C_{u,q,2}$ is the ${\rm L}^q_\mx({\rm L}^2_\lambda)$-bound of $(u_n)$ and $C_{v,2}$ is the ${\rm L}^2$-bound of $(v_n)$.
\end{theorem}

\begin{proof} 	
First, notice that the mappings
$$
\varphi(\mx)\psi(\mxi,\lambda) \mapsto
\int_{\Omega\times S}\varphi(\mx) u_n(\mx,\lambda)\overline{{\cal A}_{\bar\psi(\pi_P(\cdot,\lambda),\lambda)}
	(v_n)(\mx)}\; d\mx \,d\lambda
$$
define a sequence of linear mappings $(\mu_n)$ defined on $\widetilde{W}^p_\Pi(\Omega,S)$.
We shall prove that the sequence $(\mu_n)$ is bounded on $\widetilde{W}^p_\Pi(\Omega,S)$ with respect to the norm $\|\cdot\|_{W^p_{\Pi}}$.
Since $\widetilde{W}^p_\Pi(\Omega,S)$ is dense in ${W}^p_\Pi(\Omega,S)$ this will imply that $(\mu_n)$ is a bounded sequence of linear 
functionals on ${W}^p_\Pi(\Omega)$.
According to the Banach-Alaoglu-Bourbaki theorem, we conclude that $(\mu_n)$ is weakly-$\star$ precompact and a subsequential limit 
$\mu$ will satisfy conditions of the theorem.

Now, notice that any function belonging to $\widetilde{W}^p_\Pi(\Omega,S)$ can be approximated by sums of the form
\begin{equation*}
\sum\limits_{j=1}^N \chi_j(\mx) \psi_j(\mxi,\lambda) \;,
\end{equation*}
where $N\in\N$, $\chi_j(\mx)$, $j=1,\dots,N$, are characteristic measurable functions with disjoint supports, 
and $\psi_j\in{\rm C}^d_c(B[0,1]\times S)$.
Thus, it is enough to derive bounds for $\mu_n$ on functions of the above form.

By the properties of the commutator given in Corollary \ref{cor:komutacija} we have for a.e.~$\lambda\in S$ and 
any $j$
$$
\lim\limits_{n\to \infty} \int\limits_{\Omega} \chi_j(\mx) u_n(\mx,\lambda)
	\overline{{\cal A}_{\bar\psi_j(\pi_P(\cdot,\lambda),\lambda)}\bigl((1-\chi_j)v_n\bigr)(\mx)} \,d\mx = 0 \;,
$$
where we have used $\chi_j(1-\chi_j)=0$. Thus, as the norm of ${\cal A}_{\bar\psi_j(\pi_P(\cdot,\lambda),\lambda)}$ 
is independent of $\lambda$ (Lemma \ref{multiplierlemma1}), by the Lebesgue dominated convergence theorem (with respect to $\lambda$), we get
\begin{align*}
\limsup\limits_{n\to \infty} & \,\Bigl|\int\limits_{\Omega\times S} \sum\limits_{j=1}^N \chi_j(\mx) u_n(\mx,\lambda)
	\overline{{\cal A}_{\bar\psi_j(\pi_P(\cdot,\lambda),\lambda)}(v_n)(\mx)} \,d\mx \,d\lambda\Bigr| \\
&=\limsup\limits_{n\to \infty}\Big|\int\limits_{\Omega\times S} \sum\limits_{j=1}^N \chi_j(\mx) u_n(\mx,\lambda)
\overline{{\cal A}_{\bar\psi_j(\pi_P(\cdot,\lambda),\lambda)}(\chi_j \, v_n)(\mx)} \,d\mx \,d\lambda\Big| \;.
\end{align*}

Applying the Plancherel formula, the Fubini theorem, and the Cauchy-Bunjakovskij-Schwartz (C-B-S) inequality in $\lambda$, the above term 
is estimated by
\begin{align}
& \limsup\limits_{n\to \infty} \int\limits_{\R^{d}}\sum\limits_{j=1}^N \int\limits_S \left|\bigl(\widehat{\chi_j u_n(\cdot,\lambda)}\bigr)(\mxi)
	\psi_j(\pi_P(\mxi,\lambda),\lambda) \right| d\lambda\;  \left|\widehat{\chi_j v_n}(\mxi)\right| d\mxi \label{eq:sharper_estimate_calc}\\
&\quad \leq \limsup\limits_{n\to \infty} \int\limits_{\R^{d}} \sum\limits_{j=1}^N  \left( \int\limits_{S} \left|\bigl(\widehat{\chi_j u_n(\cdot,\lambda)}\bigr)(\mxi) 
	\right|^2 d\lambda \right)^{1/2}\!\!\!
	\left( \int\limits_S |\psi_j(\pi_P(\mxi,\lambda),\lambda)|^2 d\lambda \right)^{1/2} \left|\widehat{\chi_j v_n}(\mxi)\right| d\mxi \nonumber \\
&\quad \leq \limsup\limits_{n\to \infty} \int\limits_{\R^{d}} \sum\limits_{j=1}^N  \, A_j^{1/2}\left( \int\limits_{S} \left|\bigl(\widehat{\chi_j u_n(\cdot,\lambda)}\bigr)(\mxi) 
\right|^2 d\lambda \right)^{1/2} \left|\widehat{\chi_j v_n}(\mxi)\right| d\mxi \,, \nonumber
\end{align}
where 
\begin{equation}\label{eq:Aj}
A_j := \sup_{\mxi\in\Rd\setminus\{0\}}\int\limits_S |\psi_j(\pi_P(\mxi,\lambda),\lambda)|^2 d\lambda \,.
\end{equation}

We continue the estimate by applying first the discrete version of C-B-S inequality, and then its integral version with respect to $\mxi$,
obtaining that the above term is majorised by
\begin{align*}
\limsup\limits_{n\to \infty} & \int\limits_{\R^{d}} \Biggl(\sum\limits_{j=1}^N A_j
	\int\limits_{S} \left| \bigl(\widehat{\chi_j u_n(\cdot,\lambda)}\bigr)(\mxi) \right|^2 d\lambda \Biggr)^{1/2}
	\Biggl(\sum\limits_{j=1}^N \left|\widehat{ \chi_j v_n}(\mxi)\right|^2 \Biggr)^{1/2} d\mxi \\
& \leq \limsup\limits_{n\to \infty} \left(\,\int\limits_{\R^{d}} \sum\limits_{j=1}^N A_j
	\int\limits_{S} \left| \bigl(\widehat{\chi_j u_n(\cdot,\lambda)}\bigr)(\mxi) \right|^2 d\lambda d\mxi\right)^{1/2} 
	\left(\,\int\limits_{\Rd}\sum\limits_{j=1}^N \left|\widehat{ \chi_j v_n}(\mxi)\right|^2 d\mxi \right)^{1/2} \\
&= \limsup\limits_{n\to \infty} \left(\int\limits_{\Omega} \sum\limits_{j=1}^N A_j
	\int\limits_{S} \left|(\chi_j u_n)(\mx,\lambda) \right|^2 d\lambda d\mx\right)^{1/2} 
	\left(\int\limits_{\Omega}\sum\limits_{j=1}^N \left|(\chi_j v_n)(\mx)\right|^2 d\mx \right)^{1/2} \,,
\end{align*}
where the Plancherel formula is used in the last equality. 

As supports of $\chi_j$ are disjoint, we have
$$
\int\limits_{\Omega}\sum\limits_{j=1}^N \left|(\chi_j v_n)(\mx)\right|^2 d\mx \leq 
	\int\limits_{\Omega} |v_n(\mx)|^2 d\mx = \|v_n\|_{{\rm L}^2(\Omega)}^2 \,,
$$
while on the first factor we apply the H\"older inequality $(1/q+1/p=1/2)$ in $\mx$:
\begin{align*}
&\left(\int\limits_{\Omega} \sum\limits_{j=1}^N A_j\int\limits_{S} \left|(\chi_j u_n)(\mx,\lambda) \right|^2 d\lambda d\mx\right)^{1/2} \\
&\qquad =\left(\int\limits_{\Omega} \|u_n(\mx,\cdot)\|_{{\rm L}^2(S)}^2 \biggl(\sum\limits_{j=1}^N A_j \chi_j(\mx)\biggr) d\mx\right)^{1/2} \\
&\qquad \leq \|u_n\|_{{\rm L}^q(\Omega;{\rm L}^2(S))}
	\left(\int\limits_{\Omega} \left( \sum\limits_{j=1}^N  \chi_j(\mx) 
	\biggl( \sup_{\mxi\in\Rd\setminus\{0\}} \int\limits_S  |\psi_j(\pi_P(\mxi,\lambda),\lambda)|^2 
	d\lambda \biggr) \right)^{\frac{p}{2}} d\mx \right)^{\frac{1}{p}} \\
&\qquad = \|u_n\|_{{\rm L}^q(\Omega;{\rm L}^2(S))} 
	\left(\,\int\limits_{\Omega} \left(\sup_{\mxi\in\Rd\setminus\{0\}} \Biggl(\int\limits_S  \biggl|\sum\limits_{j=1}^N  \chi_j(\mx)\psi_j(\pi_P(\mxi,\lambda),\lambda)\biggr|^2 
	d\lambda \Biggr)^{1/2}\right)^{p} d\mx \right)^{1/p} \,,
\end{align*}
where we have used once more that $\chi_j^2=\chi_j$ and that 
$\chi_j$ have disjoint supports. 

Therefore, the final estimate obtained in the above calculations reads
\begin{align*}
\limsup\limits_{n\to \infty} \,\Bigl|\int\limits_{\Omega\times S} \sum\limits_{j=1}^N \chi_j(\mx) u_n(\mx,\lambda)
	\overline{{\cal A}_{\bar\psi_j(\pi_P(\cdot,\lambda),\lambda)}(v_n)(\mx)} \,d\mx \,d\lambda\Bigr|
\leq  C_{u,q,2} \, C_{v,2} \,\biggl\|  \sum\limits_{j=1}^N  \chi_j\psi_j	\biggr\|_{{W}^p_\Pi(\Omega,S)} \,,
\end{align*}
where $C_{u,q,2}=\limsup_n\|u_n\|_{{\rm L}^q(\Omega;{\rm L}^2(S))}$ and $C_{v,2}=\limsup_n\|v_n\|_{{\rm L}^2(\Omega)}$, 
implying the boundedness of the sequence $(\mu_n)$ in $(\widetilde{W}^{p}_\Pi(\Omega,S), \| \cdot \|_{{W}^{p}_\Pi})$.
Thus, the sequence is bounded in ${W}^{p}_\Pi(\Omega,S)$ as well and, since  ${W}^{p}_\Pi(\Omega,S)$ is a separable
Banach space, the Banach-Alaoglu-Bourbaki theorem provides the statement of the theorem.
\end{proof}


\begin{remark}\label{rem<M}
Let us note that in \eqref{eq:sharper_estimate_calc} one could consider integration with respect 
to $\mxi$ only over $|\mxi|>M$ for any fixed $M>0$. 
Indeed, $\widehat{\chi_j v_n}(\mxi)\to 0$ as $n\to \infty$ for every fixed $\mxi\in\Rd$
and $j\in\{1,2,\dots,N\}$
since $(v_n)$ is uniformly compactly supported and converges weakly-$\star$ to $0$ in 
${\rm L}^\infty(\R^d)$.
On the other hand,
$$
\int\limits_S \left|\bigl(\widehat{\chi_j u_n(\cdot,\lambda)}\bigr)(\mxi)
	\psi_j(\pi_P(\mxi,\lambda),\lambda) \right| d\lambda
	\leq \|\psi_j\|_{\mathrm{L}^\infty(B[0,1]\times S)}\|u_n\|_{\mathrm{L}^1(\Omega\times S)}<\infty \,.
$$
Thus, we can apply the Lebesgue dominated convergence theorem to conclude that the part 
of \eqref{eq:sharper_estimate_calc} in which the integration is over $|\mxi|\leq M$ 
converges to zero as $n\to \infty$.

This means that in \eqref{eq:Aj} the supremum could be taken only for 
$|\mxi|>M$, implying that $\mu$ from Theorem \ref{bilinearboundedness} satisfies 
a sharper estimate:
\begin{equation}\label{eq:mu_sharp_estimate}
|\langle \mu, \Psi\rangle| \leq \left( \int_{\Omega} \left[ \sup_{|\mxi|>M}
\left( \int_{S} |\Psi\bigl(\mx,\pi_P(\mxi,\lambda),\lambda\bigr)|^2\; d\lambda \right)^{1/2} \right]^{p} d\mx\right)^{1/p} \;, \quad \psi\in W^p_\Pi(\Omega,S) \,,
\end{equation}
for any fixed $M>0$. 
\end{remark}

We are actually able to show the following representation for the functional from the previous theorem for less regular functions
with respect to $\mx$ and $\lambda$.

\begin{corollary}\label{gen-repr}
Under the conditions of the previous theorem, let us consider a subsequence (not relabelled) that defines $\mu\in \Bigl( W^{\frac{2q}{q-2}}_{\Pi}(\Omega,S) \Bigr)'$ 
by \eqref{mu-repr-h}. Then for any $\varphi\in {\rm L}^{r}(\Omega\times S)$, $r>\frac{q}{q-1}$, and $\psi \in {\rm C}_c (S; {\rm C}_c^d(B[0,1]))$ it holds
\begin{equation}
\label{mu-repr}
\mu(\varphi \psi) = \lim_{n\to\infty}\int_{\Omega\times S}\varphi(\mx,\lambda) u_n(\mx,\lambda)\overline{{\cal A}_{\bar\psi(\pi_P(\cdot,\lambda),\lambda)}(v_n)(\mx)}\; d\mx d\lambda \;.
\end{equation}
\end{corollary}
\begin{proof}
In order to prove \eqref{mu-repr}, we shall use a fairly direct approximation argument (see also \cite[Theorem 2.2]{LM4}). 
To this end, we take a function $\varphi \in {\rm L}^r(\Omega \times S)$ and choose its approximation in ${\rm L}^{\frac{2q}{q-2}}(\Omega) \times {\rm C}_c(S) $ of the form
$$
\varphi_s(\mx,\lambda)=\sum\limits_{k=1}^s \phi_k(\mx) \chi_k(\lambda)\,, \ \ \phi_k\in {\rm L}^{\frac{2q}{q-2}}(\Omega), \, \chi_k\in {\rm C}_c(S;\R)\,,
$$ 
i.e.~$\lim\limits_{s\to \infty}\|\varphi-\varphi_s\|_{{\rm L}^r(\Omega\times S)}=0$. 
Then, we define an extension of $\mu$ by
\begin{equation*}
\mu(\varphi \psi):= \lim\limits_{s\to\infty} \mu(\varphi_s \psi) \,.
\end{equation*}

Since $\chi_k$ is real-valued and depends only on $\lambda$, by the definition
of $\mu$ it follows
\begin{align*}
\mu(\varphi_s \psi) =& \lim\limits_{n\to \infty}\sum\limits_{k=1}^s\int_{\Omega \times S} \phi_k(\mx)  u_n(\mx,\lambda)\overline{{\cal A}_{\bar\psi(\pi_P(\cdot,\lambda),\lambda)\chi_k(\lambda)}(v_n)(\mx)}\; d\mx d\lambda \\
=& \lim\limits_{n\to \infty}\sum\limits_{k=1}^s\int_{\Omega \times S} \phi_k(\mx)\chi_k(\lambda)  u_n(\mx,\lambda)\overline{{\cal A}_{\bar\psi(\pi_P(\cdot,\lambda),\lambda)}(v_n)(\mx)}\; d\mx d\lambda \\
=& \lim\limits_{n\to \infty} \int_{\Omega \times S} \varphi_s(\mx,\lambda)  u_n(\mx,\lambda)\overline{{\cal A}_{\bar\psi(\pi_P(\cdot,\lambda),\lambda)}(v_n)(\mx)}\; d\mx d\lambda \;. \\
\end{align*}
Thus, the above definition is equivalent to
\begin{equation}
\label{app-1-II}
\mu(\varphi \psi)=\lim\limits_{s\to\infty}\lim\limits_{n\to \infty} \int_{\Omega \times S} \varphi_s(\mx,\lambda)  
	u_n(\mx,\lambda)\overline{{\cal A}_{\bar\psi(\pi_P(\cdot,\lambda),\lambda)}(v_n)(\mx)}\; d\mx d\lambda \;.
\end{equation} 
This limit is well-defined as one can see from the Cauchy criterion. Indeed, for $s_1,s_2 \in \N$ by means of the H\" older inequality 
and the multiplier lemma (Lemma \ref{multiplierlemma1}) we have
\begin{equation*}
\big{|} \mu(\varphi_{s_2} \psi)-\mu(\varphi_{s_1} \psi) \big{|} \leq C \|\varphi_{s_2}-\varphi_{s_1}\|_{{\rm L}^r(\Omega\times S)} \;,
\end{equation*} 
and the constant $C$ is equal to 
$C_{\bar r,\psi}\operatorname{meas}(S)^\frac{1}{\bar r} \limsup_{n\to\infty} \|u_n\|_{\mathrm{L}^q(\Omega\times S)}\|v_n\|_{\mathrm{L}^{\bar r}(\Omega)}$,
where $\frac{1}{r}+\frac{1}{q}+\frac{1}{\bar r}=1$, and $C_{\bar r, \psi}$ is the $\mathrm{L}^{\bar r}$-bound of the Fourier multiplier operator ${\cal A}_{\bar\psi(\pi_P(\cdot,\lambda),\lambda) }$. 
Since $(\varphi_s)$ is a Cauchy sequence, the above difference can be made arbitrarily small for $s_1,s_2$ large enough, hence \eqref{app-1-II} is well defined.
The same analysis leads to 
\begin{equation*}
\lim\limits_{s\to\infty} \int_{\Omega \times S} \Bigl(\varphi(\mx,\lambda)-\varphi_s(\mx,\lambda) \Bigr)  
	u_n(\mx,\lambda)\overline{{\cal A}_{\bar\psi(\pi_P(\cdot,\lambda),\lambda)}(v_n)(\mx)}\; d\mx d\lambda=0 \,,
\end{equation*}
and the convergence is uniform with respect to $n$. Therefore, we can exchange the limits in \eqref{app-1-II}, which proves \eqref{mu-repr}.
\end{proof}

\begin{remark}\label{rem:xi_symbol}
The representation \eqref{mu-repr} holds even for $\tilde\psi(\mxi):=2\pi(1-|\mxi|)$,
which is merely continuous (at the origin $\tilde \psi$ is not smooth). 

Indeed, in the construction of the previous 
corollary we only needed that for a.e.~$\lambda$ and any $p\in(1,\infty)$
mapping $\mxi\mapsto \bar\psi\bigl(\pi_P(\mxi,\lambda),\lambda\bigr)$ is 
an $\mathrm{L}^p(\Rd)$-multiplier, with the norm 
independent of $\lambda$. 
By Lemma \ref{TVRDNJA} function $\tilde\psi$ 
satisfies this requirement. 
\end{remark}

Let us now introduce a \emph{localisation principle} for functionals
$\mu$ given by Theorem \ref{bilinearboundedness} (see also Corollary \ref{gen-repr}), 
which can serve as a way 
of proving that $\mu\equiv 0$. A similar result holds for
arbitrary continuous functionals on $W^p_\Pi(\Omega,S)$ as well.


\begin{lemma}\label{loc-p}
Under the conditions of Theorem \ref{bilinearboundedness}, 
let $\mu\in \bigl( W^{p}_{\Pi}(\Omega,S) \bigr)'$,
$p= \frac{2q}{q-2}$,
be a functional defined in Theorem \ref{bilinearboundedness}.

If a function $F\in {\rm L}^r(\Omega\times S)\otimes \mathrm{C}^d(B[0,1])$,
$r>\frac{q}{q-1}$, is such that for any compact $K\subseteq S$
it holds
\begin{equation}\label{lp-11-g}
\lim_{\eps\to 0^+}g(\eps)=0 \,,
\end{equation} where
\begin{equation*}
g(\eps):=
\int\limits_\Omega\sup\limits_{|\mxi|>1}
	\operatorname{meas}
	\Bigl\{\lambda\in K:\; \bigl|F \bigl(\mx,\lambda,\pi_P(\mxi,\lambda) \bigr)\bigr|<\eps 
	 \Bigr\} \,d\mx \,,
\end{equation*} 
and
\begin{equation}
\label{assum-11}
F \mu \equiv 0\,,
\end{equation} 
then
$$
\mu\equiv 0 \,.
$$
\end{lemma}

\begin{proof}
Let us take an arbitrary $\varepsilon>0$ and $\phi\in \mathrm{C}_c(\Omega \times B[0,1]\times S)$. 
Denote by $K$ the projection of the support of $\phi$ to the last variable $\lambda$. 
Applying \eqref{assum-11} to $\phi\frac{\overline{F}}{|F|^2+\varepsilon }$
(which is an admissible test function since it is bounded and compactly supported, hence 
with respect to the variables $\mx$ and $\lambda$ 
it can be approximated by a smooth function in 
any $\mathrm{L}^s(\Omega\times K)$, $s\in [1,\infty)$), we get
$$
0 = \left\langle\mu,\phi\frac{|F|^2}{|F|^2+\varepsilon } 
	\right\rangle = \langle\mu,\phi\rangle
	-\left\langle\mu,\phi\Bigl(\frac{|F|^2}{|F|^2+\varepsilon 
	}-1\Bigr)\right\rangle \,.
$$
Thus, it is sufficient to show that the second term on the right hand side goes to 0 as $\varepsilon \to 0$. From \eqref{eq:mu_sharp_estimate} (applied for $M=1$) we have
\begin{align*}
&\left|\left\langle\mu,\phi\Bigl(\frac{|F|^2}{|F|^2+\varepsilon 
	}-1\Bigr)\right\rangle \right| \\
&\quad \leq \|\phi\|_{\mathrm{L}^\infty} \left( \int_{\Omega} \left[ \sup_{|\mxi|>1}
\left( \int_{K} \Big| \frac{\varepsilon}{|F|^2+\varepsilon}\Big|^2\; d\lambda \right)^{1/2} \right]^{p} d\mx\right)^{1/p} \;.
\end{align*}
Thus, it is left to prove 
\begin{equation}\label{eq:loc_princ_to_prove}
\lim_{\eps\to 0^+} \left( \int_{\Omega} \left[ \sup_{|\mxi|>1}
\left( \int_{K} \Big| \frac{\varepsilon}{|F|^2+\varepsilon}\Big|^2\; d\lambda \right)^{1/2} \right]^{p} d\mx\right)^{1/p}  = 0 \,.
\end{equation}	
	
To this end, denote
$$
K^\theta(\mxi,\mx)
	:=\Bigl\{\lambda\in K:\; \bigl|F \bigl(\mx,\lambda,\pi_P(\mxi,\lambda) 
	\bigr)\bigr|
	<\theta  \Bigr\} \,.
$$ 
Let us separately analyse \eqref{eq:loc_princ_to_prove} on 
$K^{\sqrt[4]{\eps}}=K^{\sqrt[4]{\eps}}(\mxi,\mx)$
and its complement. 

By the assumption \eqref{lp-11-g} we have
\begin{align*}
\int_{\Omega} \biggl[\, & \sup_{|\mxi|>1}
\biggl( \int_{K^{\sqrt[4]{\eps}}} \Bigl| \frac{\varepsilon}{|F|^2+\varepsilon}\Bigr|^2\; d\lambda \biggr)^{1/2} \,\biggr]^{p} d\mx  \\
& \leq \int_\Omega \Bigl(\sup\limits_{|\mxi|>1}
	\operatorname{meas}\bigl(K^{\sqrt[4]{\eps}}(\mxi,\mx)\bigr)\Bigr)^{p/2} \,d\mx \\
& \leq \operatorname{meas}(K)^{\frac{p-2}{2}}
	\int_\Omega \sup\limits_{|\mxi|>1}
	\operatorname{meas}\bigl(K^{\sqrt[4]{\eps}}(\mxi,\mx)\bigr) \,d\mx
		= \operatorname{meas}(K)^{\frac{p-2}{2}}
		 g(\sqrt[4]{\eps}) \longrightarrow 0
\end{align*}
as $\eps\to 0$. Note that in the second inequality we have used that 
$p=\frac{2q}{q-2}>2$.

On the other hand, on the complement we get
\begin{align*}
\int_{\Omega} \biggl[\, & \sup_{|\mxi|>1}
\biggl( \int_{(K\setminus K^{\sqrt[4]{\eps}})} \Big| \frac{\varepsilon}{|F|^2+\varepsilon}\Big|^2\; d\lambda \biggr)^{1/2} \,\biggr]^{p} d\mx  \\
&\leq \operatorname{meas}(\Omega)\,	\Biggl\|\frac{\varepsilon }
	{\sqrt{\varepsilon} 
		+\varepsilon } 
	\Biggr\|_{\mathrm{L}^2(K\setminus K^{\sqrt[4]{\eps}})}^p \leq \operatorname{meas}(\Omega) \Big(\eps \operatorname{meas}(K)\Bigr)^{p/2} \;.
\end{align*}

Therefore, \eqref{eq:loc_princ_to_prove} holds, which, by previous observations, provides 
$\langle \mu,\phi\rangle=0$, finishing the proof.
\end{proof}

In the application of the previous lemma we shall have a specific form of 
the function $F$ for which the \emph{non-degeneracy assumption} 
\eqref{lp-11-g} simplifies. 

\begin{lemma}\label{lm:loc-p_symbol}
Let us assume that the conditions of Theorem \ref{bilinearboundedness} are fulfilled
and that in addition $a\in\mathrm{C}(S;\R^{d\times d})$.
Let us define
\begin{equation*}
F(\mx,\lambda,\mxi) := i\langle f(\mx,\lambda)\,|\,\mxi\rangle
	+2\pi (1-|\mxi|) \,,
\end{equation*}
where $f\in \mathrm{L}^r(\Omega\times S;\Rd)$, $r>\frac{q}{q-1}$.

For the function $F$ condition \eqref{lp-11-g} is equivalent to
\begin{equation}\label{eq:non-deg-simpl}
(\forall K\subset\subset S) \qquad \esssup_{\mx\in \Omega} \sup_{\mxi\in \mathrm{S}^{d-1}} \operatorname{meas} \Bigl\{\lambda\in  K : \;
\langle f(\mx,\lambda)\,|\, \mxi\rangle=\langle a(\lambda)\mxi\,|\,\mxi \rangle=0\Bigr\}=0 \;.
\end{equation}
\end{lemma}

\begin{proof}
%
Notice first that although function $F$ is only in 
${\rm L}^r(\Omega\times S)\otimes \mathrm{C}(B[0,1])$,
by Remark \ref{rem:xi_symbol} the statement of the previous lemma still applies.

Suppose that \eqref{eq:non-deg-simpl} does not hold, i.e.~there exists $c>0$, a compact set $K\subseteq S$ and $\Omega'\subseteq \Omega$ of positive measure
such that
$$
\sup_{\mxi\in \mathrm{S}^{d-1}} \operatorname{meas}
	\Bigl\{\lambda\in  K : \; \langle f(\mx,\lambda)\,|\, \mxi\rangle
	=\langle a(\lambda)\mxi\,|\,\mxi \rangle=0\Bigr\} >c>0 \ , \quad \mx\in\Omega' \,.
$$
Since for any $\eps>0$ it holds
$$
\biggl\{\lambda\in K:\; \Bigl|F \Bigl(\mx,\lambda,
{\pi_P(\mxi,\lambda)} \Bigr)\Bigr|<\eps 
\biggr\}\supseteq
\Bigl\{\lambda\in  K : \;
\langle f(\mx,\lambda)\,|\, \mxi\rangle=\langle a(\lambda)\mxi\,|\,\mxi \rangle=0\Bigr\} \,,
$$
we get
\begin{align*}
g(\eps) &\geq \int\limits_{\Omega} \sup_{|\mxi|>1} \operatorname{meas} 
	\Bigl\{\lambda\in  K : \;
	\langle f(\mx,\lambda)\,|\, \mxi\rangle=\langle a(\lambda)\mxi\,|\,\mxi \rangle=0\Bigr\} \,d\mx \\
&= \int\limits_{\Omega} \sup_{\mxi\in \mathrm{S}^{d-1}} \operatorname{meas} 
	\Bigl\{\lambda\in  K : \;
	\langle f(\mx,\lambda)\,|\, \mxi\rangle=\langle a(\lambda)\mxi\,|\,\mxi \rangle=0\Bigr\} \,d\mx 
\geq \operatorname{meas}(\Omega') \,c >0 \;,
\end{align*}
implying that \eqref{lp-11-g} does not hold as well.

To prove the opposite, assume that \eqref{eq:non-deg-simpl} holds 
but $\lim\limits_{\eps\to 0^+} g(\eps)=0 $ fails to hold. 
This means that there exist $c>0$ and a decreasing sequence of positive 
numbers $(\eps_n)$ converging to zero such that for any $n\in\N$ 
we have $g(\eps_n)>c$. 

Let us define ($\mx\in\Omega$, $|\mxi|>1$):
$$
K_n(\mxi,\mx) := \Bigl\{\lambda\in K : |\langle f(\mx,\lambda)\,|\,\mxi\rangle|
	+2\pi\langle a(\lambda)\mxi\,|\,\mxi\rangle 
	< \eps_n \bigl(|\mxi|+ \langle a(\lambda)\mxi\,|\,\mxi\rangle\bigr)\Bigr\} \,.
$$
For any $\lambda\in K_n(\mxi,\mx)$ and sufficiently large $n$ we have
\begin{align*}
&|\langle f(\mx,\lambda)\,|\,\mxi\rangle| + \langle a(\lambda)\mxi\,|\,\mxi\rangle
	\leq
	|\langle f(\mx,\lambda)\,|\,\mxi\rangle| + (2\pi-\eps_n) \langle a(\lambda)\mxi\,|\,\mxi\rangle
	<  \eps_n |\mxi|.
\end{align*} 
Dividing by $|\mxi|$ and using $|\mxi|>1$ we finally get
\begin{equation*}
K_n(\mxi,\mx)\subseteq \tilde K_n\Bigl(\frac{\mxi}{|\mxi|},\mx\Bigr) := 
	\biggl\{\lambda\in K : \Bigl|\Bigl\langle f(\mx,\lambda)\, \Big|\,
	\frac{\mxi}{|\mxi|} \Bigr\rangle\Bigr| 
	+\Bigl\langle a(\lambda)\frac{\mxi}{|\mxi|} \,\Big|\,
	\frac{\mxi}{|\mxi|} \Bigr\rangle 
< \eps_n  \biggr\} \,. 
\end{equation*}
Therefore, for any $n\in\N$ it holds
\begin{align*}
c < g(\eps_n) &= \int\limits_\Omega \sup\limits_{|\mxi|>1} 
	\operatorname{meas}\bigl(K_n(\mxi,\mx)\bigr)\,d\mx \\
&\leq \int\limits_\Omega \sup\limits_{|\mxi|>1} 
	\operatorname{meas}\biggl(\tilde K_n\Bigl(\frac{\mxi}{|\mxi|},\mx\Bigr)\biggr)\,d\mx
	= \int\limits_\Omega \sup\limits_{\mxi\in\Sdmj} 
	\operatorname{meas}\bigl(\tilde K_n(\mxi,\mx)\bigr)\,d\mx \,.
\end{align*}

By the definition of supremum, for any $n\in\N$ and $\mx\in\Omega$ there exists
$\mxi_n(\mx)\in\Sdmj$ such that 
\begin{equation*}
\int_\Omega \operatorname{meas}\Bigl(\tilde K_n\bigl(\mxi_n(\mx),\mx\bigr)\Bigr) \,d\mx 
	> \frac{c}{2} \;.
\end{equation*}
Let us define a sequence of (measurable) non-negative bounded functions 
on $\Omega$ by $\Phi_n(\mx):=\operatorname{meas}
\Bigl(\tilde K_n\bigl(\mxi_n(\mx),\mx\bigr)\Bigr)$.
Since $\Phi_n(\mx)\leq \operatorname{meas}(K)$ and $\Omega$ is bounded, 
by the so-called reverse Fatou lemma we have
\begin{equation}\label{eq:loc_Phi_n}
\frac{c}{2} \leq \limsup_{n\to\infty} \int_\Omega \Phi_n(\mx) \,d\mx 
	\leq \int_\Omega \Phi(\mx) \,d\mx \;,
\end{equation}
where we used $\Phi(\mx) := \limsup_{n\to\infty} \Phi_n(\mx)$. 

For a fixed $\mx\in\Omega$, there exists a subsequence such that
$\Phi(\mx)=\lim_{n'\to \infty}\Phi_{n'}(\mx)$. Since $\mxi_{n'}(\mx)\in\Sdmj$, we can pass
to another subsequence (not relabelled) 
to have that $(\mxi_{n'}(\mx))$ converges, and let us denote 
$\mxi(\mx)=\lim_{n'\to \infty}\mxi_{n'}(\mx)$. 

We shall prove that 
\begin{equation}\label{eq:loc_limit_Phi_n}
\lim_{n'}\Phi_{n'}(\mx) \leq \operatorname{meas}\Bigl\{\lambda\in  K : \; \langle f(\mx,\lambda)\,|\, \mxi(\mx)\rangle
=\langle a(\lambda)\mxi(\mx)\,|\,\mxi(\mx) \rangle=0\Bigr\} \ , \quad \hbox{a.e.} \ 
\mx\in\Omega \,,
\end{equation}
which by the uniqueness of the limit implies 
$$
\Phi(\mx) \leq \operatorname{meas}\Bigl\{\lambda\in  K : \; \langle f(\mx,\lambda)\,|\, \mxi(\mx)\rangle
	=\langle a(\lambda)\mxi(\mx)\,|\,\mxi(\mx) \rangle=0\Bigr\} \ , \quad \hbox{a.e.} \ 
	\mx\in\Omega \,.
$$ 
Thus, this together with \eqref{eq:loc_Phi_n} is in contradiction with
\eqref{eq:non-deg-simpl}.

Let us prove \eqref{eq:loc_limit_Phi_n} for an arbitrary, but fixed, $\mx\in\Omega$.
For $\lambda\in \tilde K_{n'}\bigl(\mxi_{n'}(\mx),\mx\bigr)$
and sufficiently large $n'$ we have
\begin{align*}
|\langle f(\mx,\lambda) & \,|\,\mxi(\mx)\rangle |
	+ \langle a(\lambda)\mxi(\mx)\,|\,\mxi(\mx)\rangle \\
&< \eps_{n'} +  \Bigl(\max_{s\in K}|a(s)||\mxi(\mx)-\mxi_{n'}(\mx)|+|f(\mx,\lambda)|
	+2|a(\lambda)|\Bigr) |\mxi(\mx)-\mxi_{n'}(\mx)| \\
&\leq \eps_{n'} +  \Bigl(1+|f(\mx,\lambda)|
+2|a(\lambda)|\Bigr) |\mxi(\mx)-\mxi_{n'}(\mx)| \,.
\end{align*}
Hence, $\tilde K_{n'}\bigl(\xi_{n'}(\mx),\mx\bigr)\subseteq \hat K_{n'}(\mx)$,
where
$$
\hat K_{n'}(\mx) := \biggl\{\lambda\in K : 
	\frac{|\langle f(\mx,\lambda)\,|\,\mxi(\mx)\rangle|
	+ \langle a(\lambda)\mxi(\mx)\,|\,\mxi(\mx)\rangle}
	{1+|f(\mx,\lambda)|+2|a(\lambda)|} < \eps_{n'}+|\mxi(\mx)-\mxi_{n'}(\mx)| \biggr\} \;,
$$
leading to
$$
\infty>\operatorname{meas}(K)\geq \operatorname{meas}(\hat K_{n'}(\mx))
	\geq \operatorname{meas}\bigl(\tilde K_{n'}\bigl(\xi_{n'}(\mx),\mx\bigr)\bigr)
	= \Phi_{n'}(\mx) \;, \quad n\in\N \,.
$$ 
Furthermore, we can pass to another subsequence (not relabelled) such that 
(for fixed $\mx$) $(|\mxi(\mx)-\mxi_{n'}(\mx)|)_{n'}$ is a decreasing sequence of real
numbers, implying that $(\hat K_{n'}(\mx))$ is a decreasing sequence of sets. 
Therefore,
$$
\begin{aligned}
\lim_{n'}\Phi_{n'}(\mx)
&\leq \lim_{n'} \operatorname{meas}
	\bigl(\hat K_{n'}(\mx)\bigr) \\
& =\operatorname{meas}\biggl\{\lambda\in  K : \;
\frac{|\langle f(\mx,\lambda)\,|\,\mxi(\mx)\rangle|
	+ \langle a(\lambda)\mxi(\mx)\,|\,\mxi(\mx)\rangle}
{1+|f(\mx,\lambda)|+2|a(\lambda)|}=0 \biggr\} \\
&= \operatorname{meas}\Bigl\{\lambda\in  K : \;
	|\langle f(\mx,\lambda)\,|\,\mxi(\mx)\rangle|
	= \langle a(\lambda)\mxi(\mx)\,|\,\mxi(\mx)\rangle=0\Bigr\} \,,
\end{aligned}
$$ 
obtaining \eqref{eq:loc_limit_Phi_n}, thus completing the proof. 
\end{proof}

\section{Proof of the main theorem}\label{sec:velocity_averaging}

In this section, we shall apply previously developed tools to prove a velocity averaging result for the sequence of equations given in the Introduction:
\begin{equation}
\begin{split}
\Div_\mx \bigl(f(\mx,\lambda) u_n(\mx,\lambda) \bigr) =& \,\Div_\mx \bigl(\Div_\mx \left( a(\lambda) u_n(\mx,\lambda) \right)\bigr) \\
&\qquad\quad + \pa_{\lambda} G_n(\mx,\lambda) + \Div_\mx P_n(\mx,\lambda) \qquad \hbox{in} \ {\cal D}'(\R^{d+1})\;,
\end{split} \tag{\ref{app-1}}
\end{equation}
where we assume that conditions (a)-(e) are fulfilled. 
\smallskip

%

\noindent {\bf Proof of Theorem \ref{velocity averaging}:}
Let $\Omega\times S\subseteq \R^{d+1}$ be a bounded open subset such that supports of all functions $u_n$ are contained in it.
Let us take a bounded sequence of functions $(v_n)$ uniformly compactly supported on $\Omega$ and
weakly-$\star$ converging to zero in ${\rm L}^\infty(\Omega)$, which we take at this moment to be arbitrary. At the end 
of the proof the precise choice will be made. 
Let us pass to a subsequence of both $(u_n)$ and $(v_n)$ (not relabelled) which defines a bounded linear functional 
$\mu\in\bigl(W^\frac{2q}{q-2}_\Pi(\Omega,S)\bigr)'$
according to Theorem \ref{bilinearboundedness}, and consider its extension given by Corollary \ref{gen-repr}.

For arbitrary $\varphi \in {\rm C}_c(\Omega)$ and 
$\psi\in{\rm C}^{d+1}_c(B[0,1]\times S)$ we set
\begin{equation}\label{tf-1}
\theta_n(\mx,\lambda) := \overline{{\cal A}_{\frac{\bar\psi\left(\pi_P(\cdot,\lambda),\lambda\right)}
	{|\cdot|+\langle a(\lambda) \cdot\,|\,\cdot\rangle}}(\varphi v_n)(\mx)} \;.
\end{equation}
Testing \eqref{app-1} by $\theta_n$, i.e.~multiplying it by $\theta_n$, integrating over $\Omega\times S$, and applying the integration by parts we get the following:
	
\begin{align}
0=2\pi\sum\limits_{j=1}^d  & \,\int\limits_{\Omega\times S} f_j(\mx,\lambda) u_n(\mx,\lambda) 
	\overline{{\cal A}_{\frac{i\xi_j\bar\psi(\pi_P\left(\mxi,\lambda),\lambda\right)}{|\mxi|+\langle a(\lambda)\mxi \,|\,\mxi\rangle}}
	(\varphi v_n)(\mx)} \,d\mx d\lambda \label{l1a} \\
&-2\pi\int\limits_{\Omega\times S}  u_n(\mx,\lambda) \overline{{\cal A}_{\frac{2\pi\langle a(\lambda) 
		\mxi \,|\, \mxi\rangle \bar\psi\left(\pi_P(\mxi,\lambda),\lambda\right)}{|\mxi|+\langle a(\lambda)\mxi \,|\,\mxi\rangle}}(\varphi v_n)(\mx)} \,d\mx d\lambda \label{l2a} \\
&-\int\limits_{\Omega\times S} G_n(\mx,\lambda) \overline{{\cal A}_{\pa_{\lambda}\frac{\bar\psi\left(\pi_P(\mxi,\lambda),\lambda\right)}{|\mxi|+
		\langle a(\lambda)\mxi \,|\,\mxi\rangle}}(\varphi v_n)(\mx)} \,d\mx d\lambda \label{l3b}\\
&-\sum\limits_{j=1}^d \,\int\limits_{\Omega\times S}  P^n_j(\mx,\lambda) \overline{{\cal A}_{\frac{2\pi i\xi_j\bar\psi\left(\pi_P(\mxi,\lambda),\lambda\right)}
		{|\mxi|+\langle a(\lambda)\mxi \,|\,\mxi\rangle}}(\varphi v_n)(\mx)} \,d\mx d\lambda \;,\label{l4a}
\end{align} where we have used $\partial_{x_j}{\cal A}_\psi={\cal A}_{2\pi i\xi_j\psi}$, which holds according to our definition of the Fourier transform. 
Line \eqref{l3b} is to be understood as
$$
\int_S \left\langle G_n(\cdot,\lambda),\, \overline{{\cal A}_{\pa_{\lambda}\frac{\bar\psi\left(\pi_P(\cdot,\lambda),\lambda\right)}{|\cdot|+\langle a(\lambda)\cdot\,|\,\cdot\rangle}}
	(\varphi v_n)} \right\rangle \,d\lambda \,,
$$ 
where $\langle\cdot,\cdot\rangle$ represents the dual product between $\mathrm{W}^{-\frac{1}{2},r}_{loc}(\Rd)$ and $\mathrm{W}^{\frac{1}{2},r}_{c}(\Rd)$.

Let us consider term by term in the above expression as $n$ goes to infinity along the chosen subsequence.

Symbols of the Fourier multiplier operators in \eqref{l1a} and \eqref{l2a}
can be rewritten as
\begin{align*}
\frac{i\xi_j\bar\psi\bigl(\pi_P(\mxi,\lambda),\lambda\bigr)}
	{|\mxi|+\langle a(\lambda)\mxi \,|\,\mxi\rangle} 
	&= (\overline{\psi_j\psi})\bigl(\pi_P(\mxi,\lambda),\lambda\bigr) \\
\frac{2\pi\langle a(\lambda)\mxi \,|\,\mxi\rangle\bar
	\psi\bigl(\pi_P(\mxi,\lambda),\lambda\bigr)}
	{|\mxi|+\langle a(\lambda)\mxi \,|\,\mxi\rangle}
	&= (\overline{\tilde \psi\psi})\bigl(\pi_P(\mxi,\lambda),\lambda\bigr) \;,
\end{align*}
where $\psi_j(\mxi):=-i\xi_j$ and $\tilde{\psi}(\mxi):=2\pi(1-|\mxi|)$.
Thus, by applying first Corollary \ref{cor:komutacija} in order to move 
$\varphi$ outside of the Fourier multiplier operators, and then
Corollary \ref{gen-repr} (see also Remark \ref{rem:xi_symbol}),
the limit of the sum of \eqref{l1a} and \eqref{l2a} is equal to
\begin{equation}\label{eq:velocityHm}
-2\pi \Bigl\langle \mu,  F(\mx,\mxi,\lambda)\varphi(\mx)\psi(\mxi,\lambda)
	\Bigr\rangle \,,
\end{equation}
where $F(\mx,\mxi,\lambda) = i\langle f(\mx,\lambda)\,|\,\mxi\rangle
+2\pi\bigl(1-|\mxi|\bigr)$.


Unlike the situation with \eqref{l1a} and \eqref{l2a}, term \eqref{l4a} 
is zero at the limit $n\to\infty$. 
Indeed, $P_j^n$ strongly converges in ${\rm L}^{p_0}_{loc}(\Rd\times \R)$, 
while ${\cal A}_{(\overline{\psi_j\psi})(\pi_P(\cdot,\cdot),\cdot)}(\varphi v_n)$ 
weakly converges to zero in ${\rm L}^{p_0'}(\Omega\times S)$ 
by Lemma \ref{multiplierlemma1}, and the integration is over 
relatively compact set $\Omega\times S$.

The symbol appearing in \eqref{l3b} we divide into two parts, namely
\begin{equation}\label{eq:Gn_symbol_I}
\frac{(\partial_\lambda\bar\psi)\left(\pi_P(\mxi,\lambda),\lambda\right)}{|\mxi|+
	\langle a(\lambda)\mxi \,|\,\mxi\rangle}
\end{equation}
and
\begin{equation}\label{eq:Gn_symbol_II}
\partial_\lambda\left(\frac{1}{|\mxi|+\langle a(\lambda)\mxi\,|\,\mxi\rangle}\right)
	\biggl(\bar{\psi}\bigl(\pi_P(\mxi,\lambda),\lambda\bigr)
	+\sum_{j=1}^{d}\bigl(\xi_j\partial_{\xi_j}\psi(\mxi,\lambda)\bigr)\circ
	\bigl(\pi_p(\mxi,\lambda\bigr)\bigr)\biggr) \,.
\end{equation}
Let us study first the part of \eqref{l3b} associated to \eqref{eq:Gn_symbol_I}.

By Lemma \ref{pomocna} (given in the Appendix), Lemma \ref{multiplierlemma1}
and the Lebesgue dominated convergence theorem (applied for the integration in 
$\lambda$) we have for any $r\in [1,\infty)$ 
$$
{\cal A}_{\frac{(\partial_\lambda\bar\psi)\left(\pi_P(\mxi,\lambda),\lambda\right)}
	{|\mxi|+\langle a(\lambda)\mxi,\mxi\rangle}}(\varphi v_n)  =
	{\cal A}_{\frac{1}{|\mxi|+\langle a(\lambda)\mxi,\mxi\rangle}}
	\Bigl({\cal A_{(\partial_\lambda\bar\psi)\left(\pi_P(\mxi,\lambda),\lambda\right)}}(\varphi v_n)\Bigr) 
	\rightharpoonup 0 \quad \hbox{weakly in} \ 
	{\rm L}^r(S;\mathrm{W}^{1,r}(\Omega))\,,
$$
where we have used that $\Omega\times S$ is relatively compact. 
This, together with the assumption of $(G_n)$, implies the convergence to zero
of the part of \eqref{l3b} associated to \eqref{eq:Gn_symbol_I}.

The Fourier multiplier operator associated to the second factor of 
\eqref{eq:Gn_symbol_II} is bounded on $\mathrm{L}^r(\Rd)$ for any 
$r\in(1,\infty)$ uniformly in $\lambda$ (Lemma \ref{multiplierlemma1}),
while by Lemma \ref{multiplierlemma2} the first factor defines a bounded
operator (uniformly in $\lambda$) from $\mathrm{L}^r(\Rd)$ to
$\mathrm{W}^{\frac{1}{2},r}(\Rd)$ for any $r\in(1,\infty)$. 
Thus, the overall conclusion is that for any $r\in(1,\infty)$
$$
{\cal A}_{\partial_\lambda\left(\frac{1}
	{|\mxi|+\langle a(\lambda)\mxi\,|\,\mxi\rangle}\right)
	\left(\bar{\psi}(\pi_P(\mxi,\lambda),\lambda)
	+\sum_{j=1}^{d}\bigl(\xi_j\partial_{\xi_j}\psi(\mxi,\lambda)\bigr)\circ
	(\pi_p(\mxi,\lambda))\right) }(\varphi v_n) \rightharpoonup 0 \,,
$$
weakly in ${\rm L}^r(S;\mathrm{W}^{\frac{1}{2},r}(\Omega))$, implying the convergence to zero of the part of \eqref{l3b} 
associated to \eqref{eq:Gn_symbol_II}.

Collecting the previous considerations, we get from \eqref{l1a}-\eqref{l4a} after letting $n\to \infty$ that \eqref{eq:velocityHm} is 
the only non-trivial term, thus we have:
\begin{equation*}
\Bigl\langle \mu,  F(\mx,\mxi,\lambda)\varphi(\mx)\psi(\mxi,\lambda)\Bigr\rangle=0 \,.
\end{equation*}
Since $F$ satisfies the non-degeneracy assumption \eqref{eq:non-deg-simpl}
(see Lemma \ref{lm:loc-p_symbol}), by Lemma \ref{loc-p} we conclude from above
that 
$$
\mu\equiv 0 \;.
$$

	Let us assume that $u_n$ is real valued (if not, we just apply this procedure to the real and imaginary parts of $u_n$ separately).
	Let us take arbitrary real valued functions $\varphi\in{\rm C}_c(\Omega)$ and $\rho\in{\rm C}_c(S)$.
	As $\big({\rm sgn}(\int_{S}\rho(\lambda)u_n(\mx,\lambda) \,d\lambda)\big)$ is bounded in ${\rm L}^\infty(\Omega)$, it has a weakly-$\star$
	converging subsequence, whose limit we denote by $V\in{\rm L}^\infty(\Omega)$. We pass to that subsequence (not relabelled), and choose for
	$v_n$ in \eqref{tf-1}:
	$$
	v_n(\mx)=\varphi(\mx) \left({\rm sgn}\left(\int_{S}\rho(\lambda)u_n(\mx,\lambda) \,d\lambda\right) -V(\mx) \right) \,.
	$$
	As the subsequence defines the same functional $\mu$, by Theorem \ref{bilinearboundedness} we conclude
	\begin{align*}
	\lim\limits_{n\to \infty} & \int_{\Omega} \varphi(\mx)^2 \left| \int_{S}\rho(\lambda)u_n(\mx,\lambda)d\lambda \right| \,d\mx\\
	&=\lim\limits_{n\to \infty} \int_{\Omega} \varphi(\mx)^2 \left(\int_{S}\rho(\lambda)u_n(\mx,\lambda)d\lambda\right) \, 
		{\rm sgn}\left(\int_{S}\rho(\lambda)u_n(\mx,\lambda)d\lambda\right) \,d\mx \\
	&=\lim\limits_{n\to \infty} \int_{\Omega\times S} \varphi(\mx) \rho(\lambda)u_n(\mx,\lambda) v_n(\mx) \,d\mx \,d\lambda \\
	&= \langle \mu, \rho \varphi \otimes 1  \rangle=0 \,,
	\end{align*}
	where in the second equality we have used that $(u_n)$ converges weakly to zero. 
	Thus, the proof is over.
\endproof
\smallskip

\begin{remark}\label{rem:only_bdd}
Let us remark that the conditions of Theorem \ref{velocity averaging} 
(as well as of corollaries \ref{cor1} and \ref{velocity averaging_Gn_in_Lr}) 
can be relaxed
by assuming that $(u_n)$ is merely bounded in $\mathrm{L}^q(\Omega\times S)$,
while $(G_n)$ and $(P_n)$ strongly precompact in 
${\rm L}_{loc}^{r_0}(\R; {\rm W}_{loc}^{-1/2,r_0}(\R^d))$
and ${\rm L}_{loc}^{p_0}(\R^d\times\R;\Rd)$ respectively.
In that case there exists a subsequence $(u_{n'})$ such that for any 
$\rho\in\mathrm{C}_c(S)$ the sequence 
$(\int_S\rho(\lambda)u_{n'}(\mx,\lambda)\,d\lambda)$ converges toward 
$\int_S\rho(\lambda)u(\mx,\lambda)\,d\lambda$, where $u$ denotes the weak limit
of $(u_{n'})$. 
\end{remark}

We note that matrix function \eqref{eq:matrix_2d_notb} that we have provided in Example \ref{ex3}(b) 
is irregular only in one point. These situations we can handle by cutting off isolated singularities using appropriate 
$\lambda$-compactly supported functions. Let us briefly explain how to prove Corollary \ref{cor1}.
\smallskip

\noindent {\bf Proof of Corollary \ref{cor1}:} 
Let us denote by $S_0\subseteq S$ the set of points for which condition 
($\tilde{\mathrm{b}}$) is satisfied, i.e.~for any $\lambda\in S_0$ there 
exits $\eps(\lambda)>0$ such that $a=\sigma^T\sigma$ and 
$\sigma\in \mathrm{C}^{0,1}\bigl((\lambda-\eps(\lambda),\lambda+\eps(\lambda));
\R^{d\times d}\bigr)$. 
By ($\tilde{\mathrm{b}}$)
it is clear that $S_0$ is an open set of full Lebesgue measure, 
i.e.~$\operatorname{meas}(S_0)=\operatorname{meas}(S)$.

Let us take a countable set $\tilde{S}\subseteq S_0$ such that
\begin{equation}\label{eq:S0_cover}
\bigcup_{\lambda\in \tilde{S}}\bigl(\lambda-\eps(\lambda), 
	\lambda+\eps(\lambda)\bigr)  =S_0 \,.
\end{equation}
One way to obtain such $\tilde{S}$ is by considering an increasing sequence of 
compact sets $(S_k)_{k\in\N}$ which exhausts open set $S_0$ and such that 
$S_k\subseteq \mathop{Int}S_{k+1}$. Then for any $k\in\N$
we can reduce open cover 
$\bigl\{\bigl(\lambda-\eps(\lambda), \lambda+\eps(\lambda)\bigr) : \lambda\in S_0\bigr\}$
to a finite subcover of $S_k$, and let us denote by $\tilde{S}_k$ a finite set of 
corresponding $\lambda$'s. Now we can take $\tilde{S}=\bigcup_{k=1}^\infty \tilde{S}_k$.

Fix $\lambda_0\in \tilde{S}$. Then, in the proof of Theorem \ref{velocity averaging} given above, we simply take $\rho\in \mathrm{C}_c(\lambda_0-\eps,\lambda_0+\eps)$, where $\eps=\eps(\lambda_0)$, to obtain \eqref{cnv*}. 
Since $\tilde{S}$ is countable, we can take a subsequence $(u_{n'})$ 
in \eqref{cnv*} independent of $\lambda_0\in \tilde{S}$. 
 
From here, the statement of the corollary immediately follows.
Indeed, let us first take an arbitrary $ \rho \in \mathrm{C}_c(S_0)$.
Applying a partition of unity subordinate to \eqref{eq:S0_cover} and the previous 
conclusions, it is easy to see that \eqref{cnv*} still holds for the same sequence
$(u_{n'})$.  
Now, let us consider a general situation where $\rho\in\mathrm{C}_c(S)$.
Let us denote by $\vartheta_k\in\mathrm{C}(S;[0,1])$ a cut off function 
such that $\vartheta_k\equiv 1$ on $S_k$ and $\vartheta_k\equiv 0$ on $S\setminus S_{k+1}$.
Since $\lim_{k\to\infty}\operatorname{meas}(S_k)=\operatorname{meas}(S)$
and $(u_n)$ is bounded in $\mathrm{L}^q(\Omega\times S)$, $q>2$, the sequence 
$(\int_S \rho(\lambda)(1-\vartheta_k(\lambda))u_{n'}(\mx,\lambda)\,d\lambda)_k$
converges to zero strongly in $\mathrm{L}^1_{loc}(\Rd)$, where the convergence 
is uniform with respect to $n'$. 
Finally, as $\rho\vartheta_k\in\mathrm{C}_c(S_0)$, the claim follows by the previous 
observations. 
%
%
 \endproof


\begin{remark}
\label{rem:weaker_assump}

With the exception of Lemma \ref{multiplierlemma2}, 
all other (multiplier) results used in the above proof of 
Theorem \ref{velocity averaging}
holds under a weaker assumption on $a$:
\begin{enumerate}
\item[b')] $a\in {\rm C}^{0,1}(S; \R^{d\times d})$ is such that, 
for every $\lambda\in S$, $a(\lambda)$ is a symmetric and positive 
semi-definite matrix.
\end{enumerate}
Thus, a sufficient assumption on the diffusion matrix $a$ 
under which the statement of Theorem \ref{velocity averaging}
holds is that $a$ satisfies (b') and that \eqref{l3b} tends to 
0, as $n\to\infty$.
Moreover, it is enough to have local estimates of the multipliers with respect to $\lambda$ (see condition ($\tilde{\mathrm{b}}$) in the statement of Corollary 
\ref{cor1}). 

To be more specific, a possible weakening of assumption (b) which 
preserves the statement of Theorem \ref{velocity averaging} is:
\begin{enumerate}
\item[b'')] $a$ satisfies (b') and for a.e.~$\lambda_0\in S$ there exists its neighborhood $U(\lambda_0)$ such that for any
$\lambda\in U(\lambda_0)$ 
$$
\mxi\mapsto \frac{|\mxi|^{1/2}\langle a'(\lambda)\mxi\,|\,\mxi\rangle}
	{\bigl(|\mxi|+\langle a(\lambda)\mxi\,|\,\mxi\rangle\bigr)^2}
$$
is an $\mathrm{L}^{r_0}$-multiplier ($r_0$ given in condition (d)), 
with the norm uniformly bounded with respect to $\lambda \in U(\lambda_0)$. 
\end{enumerate} 
By lemmata \ref{multiplierlemma1} and \ref{multiplierlemma2} it is 
clear, which was also used in the previous proof, that (b'') is 
implied by ($\tilde{\mathrm{b}}$) (see Corollary \ref{cor1}). 
However, (b') does not imply (b'') in general. 
For instance, just consider \eqref{eq:matrix_2d_notb_b} from the Introduction.
\end{remark}

Although at this moment we cannot obtain the result by imposing only (b') on the diffusion matrix $a$,
while keeping all the other assumptions intact 
(see the previous remark), under stronger assumptions on $(G_n)$, given in Corollary 
\ref{velocity averaging_Gn_in_Lr}, that can be done. 
A proof of Corollary \ref{velocity averaging_Gn_in_Lr} we briefly explain below.
\smallskip

\noindent {\bf Proof of Corollary \ref{velocity averaging_Gn_in_Lr}:}
If we replace (d) by
\begin{enumerate}
\item [d')] $G_n \to 0$ strongly in ${\rm L}_{loc}^{r_0}(\R^d_\mx\times\R_\lambda)$ for some 
$r_0\in(1,\infty)$;
\end{enumerate}
then it is sufficient that \eqref{eq:Gn_symbol_II} is an $\mathrm{L}^p$-multiplier,
$p\in(1,\infty)$, which is ensured by Lemma \ref{multiplierlemma1}.
Thus, the statement of Corollary \ref{velocity averaging_Gn_in_Lr} holds. 
\endproof

\smallskip

If $(u_n)$ has a compact support only with respect to $\lambda$, the statement of the previous theorem still holds. 
Indeed, one just needs to test \eqref{app-1} by $\tilde\varphi\theta_n$ instead of $\theta_n$ (given by \eqref{tf-1}) 
for an arbitrary $\tilde\varphi\in{\rm C}_c(\Rd)$. By repeating the rest of the analysis of the proof of Theorem 
\ref{velocity averaging} we obtain that the functional $\mu$ from Theorem \ref{bilinearboundedness} corresponding to
$(\tilde\varphi u_n)$ equals zero, implying the strong convergence to zero in ${\rm L}^1_{loc}(\Rd)$ of 
$(\int_\R \rho(\lambda)\tilde\varphi(\mx) u_n(\mx,\lambda) d\lambda)$ for any $\rho\in{\rm C}_c(\R)$.
Due to arbitrariness of $\tilde{\varphi}$ we get the claim which we formulate in the following corollary.

\begin{corollary}\label{velocity averaging-cor}
	Let $d\geq 2$ and let $u_n\rightharpoonup 0$ in ${\rm L}^q_{loc}(\Rd\times\R)$, for some $q>2$, is
	uniformly compactly supported with respect to $\lambda$ on $S\subset\subset \R$. 
	Let $(u_n)$ satisfy the sequence of equations \eqref{app-1} whose coefficients satisfy conditions 
	(b), (d), (e), and
	\begin{enumerate}
	\item [c')] There exists $p>\frac{q}{q-1}$ ($p>1$ if $q=\infty$) such that for any 
	compacts $K_1\subseteq \Rd$ and $K_2\subseteq S$ it holds
	$f\in {\rm L}^p( K_1\times S; \Rd)$ and
	\begin{equation*}
	\esssup\limits_{\mx\in K_1} \sup\limits_{\mxi \in \mathrm{S}^{d-1}} 
	\operatorname{meas}\bigl\{\lambda \in K_2 :\, \langle f(\mx,\lambda)\,|\,\mxi\rangle=
	\langle a(\lambda)\mxi \,|\,\mxi \rangle=0 \bigr\}=0 \;.
	\end{equation*}
	\end{enumerate}
	
	Then there exists a subsequence $(u_{n'})$ such that for any $\rho\in {\rm C}_c(\R)$,
	$$
	\int_\R \rho(\lambda) u_{n'}(\mx,\lambda) \,d\lambda \,\longrightarrow\, 0 \quad\hbox{strongly in} \ {\rm L}^1_{loc}(\Rd) \;.
	$$
\end{corollary}


\section{Degenerate parabolic equation with rough coefficients -- existence proof}\label{sec:application}

In this section, we prove existence of a weak solution to the Cauchy problem for the advection-diffusion
equation:

\begin{align}
\label{d-p} \pa_t u+\Div_{\mx} \mff(t,\mx,u)&=D^2 \cdot A(u) \;, \quad (t,\mx)\in\R^+\times\Rd =: \R_+^{d+1}\;, \\
\label{ic}
u|_{t=0}&=u_0 \in {\rm L}^1({\R}^d)\cap {\rm L}^\infty({\R}^d) \;.
\end{align}
where $D^2\cdot A(u) = \sum_{i,j} \partial_{x_ix_j}^2[A(u)]_{ij}$. The equation describes a flow governed by

\begin{itemize}

\item the convection effects (bulk motion of particles) which are
represented by the first order terms;

\item diffusion effects which are represented by the second order
term and the matrix
$A(\lambda)=[A_{ij}(\lambda)]_{i,j=1,\dots,d}$ (more
precisely its derivative with respect to $\lambda$; see
\eqref{stand} below) describes direction and intensity of the diffusion;

\end{itemize} The equation is degenerate in the sense that the derivative of the diffusion matrix $A'$ can be equal to zero in some direction.
Roughly speaking, if this is the case, i.e.~if for some vector $\mxi \in \R^d$ we have $\langle A'(\lambda)\mxi \,|\, \mxi\rangle =0$,
then diffusion effects do not exist at the point $\mx$ for the state $\lambda$ in the direction $\mxi$.

Recently, several existence results for \eqref{d-p} in the case when the coefficients are irregular were obtained. In \cite{LM2, pan_jms,
sazh} the authors considered ultra-parabolic equations, while in \cite{HKMP} a degenerate parabolic equation was considered and a similar result as in Theorem \ref{main-tvr} below is obtained. 

Roughly speaking, in \cite{HKMP}, the authors had the assumptions that the flux $\mff(t,\mx,\lambda)$ is merely continuous with respect to $\lambda$ and $\max\limits_{|u|<M} |\mff(\mx,u)|\in {\rm L}^2_{loc}(\R^d)$ for every $M>0$. 
However, we still generalise this result by assuming the following for the coefficients of \eqref{d-p} (keep in mind that we need only ${\rm L}^p$, $p>1$ assumptions on the flux unlike the ${\rm L}^2$ assumptions from \cite{HKMP}):

\begin{itemize}
\item[i)] There exist $\alpha,\beta\in\R$ such that the initial data $u_0\in{\rm L}^1({\R}^d)\cap {\rm L}^\infty({\R}^d)$ 
is bounded between $\alpha$ and $\beta$ and the flux equals zero at $\lambda=\alpha$ and $\lambda=\beta$:
\begin{equation}
\label{bnd-fl}
\alpha\leq u_0(\mx) \leq \beta \quad {\rm and} \quad \mff(t,\mx,\alpha)=\mff(t,\mx,\beta)=0 \ \ {\rm a.e.} \ \ (t,\mx)\in \R^{d+1}_+ \;.
\end{equation}

\item[ii)] The convective term $\mff=\mff(t,\mx,\lambda)$ belongs to ${\rm C}^1([\alpha,\beta];{\rm L}^p_{loc}(\R^{d+1}_+))$ for some $p>1$, 
$$
\sup_{\lambda\in [\alpha,\beta]}|f(\cdot,\cdot,\lambda)|\in L^p_{loc}(\R^{d+1}_+) \,,
$$
and
\begin{equation*}
\Div_{\mx} \mff \in
{\cal M}(\R^{d+1}_+\times [\alpha,\beta]) \,,
\end{equation*} 
where ${\cal M}(X)$ denotes the space of Radon measures on $X\subseteq\Rd$. 

\item[iii)] The matrix $A(\lambda)\in{\rm C}^{1,1}([\alpha,\beta];\R^{d\times d})$
is symmetric and non-decreasing with respect to $\lambda\in [\alpha,\beta]$, i.e.~the (diffusion) matrix $a(\lambda):=A'(\lambda)$ satisfies
\begin{equation*}
\langle a(\lambda)\mxi \,|\, \mxi \rangle \geq 0 \,,
\end{equation*} 
and $a$ satisfies (b) (or ($\tilde{\mathrm{b}}$)) from the Introduction.

\item[iv)] $f:=\pa_\lambda\mff$ and $a=A'$ satisfy non-degeneracy assumption:
for any compact $K\subseteq\R_+^{d+1}$ it holds
	\begin{equation*}
	\esssup\limits_{(t,\mx)\in K} \sup\limits_{(\tau,\mxi) \in \mathrm{S}^{d}} 
	\operatorname{meas}\bigl\{\lambda \in [\alpha,\beta] :\, \tau+\langle f(t,\mx,\lambda)\,|\,\mxi\rangle=
	\langle a(\lambda)\mxi \,|\,\mxi \rangle=0 \bigr\}=0 \;.
	\end{equation*}
\end{itemize}

\begin{remark}\label{remX} 
	We note that condition (i) essentially provides the maximum principle for entropy solutions of \eqref{d-p}-\eqref{ic}. More precisely, if the flux is $C^1$ with respect to all the variables and the initial data $u_0$ are bounded between two stationary (and thus both entropy and weak) solutions $\alpha$ and $\beta$ (which is why we assume \eqref{bnd-fl}), then the entropy solution to \eqref{d-p}-\eqref{ic} remains bounded between $\alpha$ and $\beta$. This is called the maximum principle (see e.g.~\cite{CK}).  For more details, see the proof of Theorem \ref{main-tvr}. 
\end{remark}

Remark that equation \eqref{d-p} can be rewritten in the standard (more usual) form as follows (cf.~\cite{CK}):
\begin{equation}
\label{stand} \pa_t u+\Div_\mx 
\mff(t,\mx,u)=\Div_\mx (a(u) \nabla_\mx u) \,,
\end{equation} 
where $a=A'$.

Thus, by proving existence of solutions to equation \eqref{d-p}, we shall 
prove existence of solutions to Cauchy problems for a degenerate parabolic equation in the standard form given above.
We note that in \cite{CK} one can find a well-posedness results for entropy solutions of \eqref{stand}-\eqref{ic} under regularity assumptions for the flux. We stress that the concept of entropy solutions providing well-posedness does not exist for the degenerate parabolic equation with discontinuous flux.

The main theorem of the section is the following.

\begin{theorem}
\label{main-tvr} Let $d\geq 1$, and let $u_0,\mff$ and $A$ satisfy conditions (i)-(iv) above. 
Then there exists a weak solution to \eqref{d-p} augmented with the initial condition \eqref{ic}.
\end{theorem}

\begin{proof} 
Consider the sequence of admissible solutions to the following regularised Cauchy problems
\begin{align}
\label{d-p-r}
\pa_t u_n+\Div_{\mx} \mff_n(t,\mx,u_n)=&D^2 \cdot
A(u_n) \,,  
\\
u_n|_{t=0}=&u_0 \in {\rm L}^1(\R^d)\cap {\rm L}^\infty(\R^d) \,,
\label{i-c-r}
\end{align} 
where $\mff_n$ is a smooth regularisation with respect to $(t,\mx)$ of $\mff$
obtained by the convolution of $\mff$ with a standard mollifier.
Thus,
\begin{equation}\label{eq:fn-alpha-beta}
\mff_n(t,\mx,\alpha)=\mff_n(t,\mx,\beta)=0 \;,
\end{equation}
and for any compact $K\subseteq\R^{d+1}_+$ it holds (see \cite[(4.10)]{LM4})
\begin{equation}
\label{f-conv}
\lim_{n\to\infty}\Bigl\| \sup_{\lambda\in[\alpha,\beta]}\bigl|\mff_n(\cdot,\cdot,\lambda)-\mff(\cdot,\cdot,\lambda)\bigr|\Bigr\|_{{\rm L}^p(K)} =
\lim_{n\to\infty}\Bigl\| \sup_{\lambda\in[\alpha,\beta]}\bigl|\partial_\lambda\mff_n(\cdot,\cdot,\lambda)-\partial_\lambda\mff(\cdot,\cdot,\lambda)\bigr|\Bigr\|_{{\rm L}^p(K)} =0 \,.
\end{equation}

It is well known that there exists a solution $u_n$ to such an equation
satisfying for a.e.~$\lambda\in \R$ the Kru\v zkov-type entropy equality holds
(see \cite{CK} where existence was shown under much more restrictive conditions), i.e.
\begin{align}
\label{e-c-n} \pa_t |u_n-\lambda|&+\Div_\mx\Bigl( {\rm
sgn}(u_n-\lambda)\bigl(\mff_n(t,\mx,u)-\mff_n(t,\mx,\lambda) \bigr)
\Bigr)\\
&-D^2\cdot\Bigl( {\rm sgn}(u_n-\lambda)
\bigl(A(u_n)-A(\lambda)\bigr)\Bigr) = -\zeta_n(t,\mx,\lambda) \,,
\nonumber
\end{align} 
where $\zeta_n\in {\rm C}(\R_\lambda;w\star-{\cal M}(\R^{d+1}_+))$ and ${\cal M}(\R^{d+1}_+)$ is the space of Radon measures. We also note that the sequence $(\zeta_n)$ is bounded in ${\cal M}(\R^{d+1}_+\times \R)$ (when we consider $\zeta_n$ as a Radon measure in all three variables).

By finding derivative with respect to $\lambda$ in \eqref{e-c-n}, we see that the (kinetic) function
\begin{equation*}
h_n(t,\mx,\lambda)={\rm sgn}(u_n(t,\mx)-\lambda)=
 -\pa_\lambda |u_n(t,\mx)-\lambda|
\end{equation*} is a weak solution to the following linear equation:
\begin{align}
\label{kinetic-n} \pa_t h_n + \Div_\mx \bigl( {f}_n(t,\mx,\lambda)
h_n\bigr)-D^2\cdot\bigl( a(\lambda) h_n\bigr) = \pa_\lambda
\zeta_n(t,\mx,\lambda)\, ,
\end{align}
where ${f}_n=\pa_\lambda\mff_n$ and $a=A'$.

According to (i) and \eqref{eq:fn-alpha-beta}, we know that $(u_n)$ satisfies the maximum principle, i.e.~every $u_n$, $n\in \N$, remains bounded between $\alpha$ and $\beta$ (see Remark \ref{remX}). Therefore, the sequence $(\zeta_n)$ is bounded 
in $\mathrm{C}(\R_\lambda; w\star-{\cal M}(\R^{d+1}_+))$. This implies that $(\zeta_n)$ is actually strongly precompact in ${\rm L}^r_{loc}(\R_\lambda; \mathrm{W}_{loc}^{-\frac{1}{2},q}(\R^{d+1}_+))$ 
for any $r\geq 1$ and $q\in \bigl[1,\frac{d+1}{d+1-\frac{1}{2}}\bigr)$ 
(one can prove this in the same manner as in \cite[Theorem 1.6]{Evans} using 
that $\mathrm{W}^{\frac{1}{2},s}(\R^{d+1})$ is compactly embedded into 
$\mathrm{C}(\mathop{Cl} K)$, $K\subset\subset\R^{d+1}$, for $\frac{s}{2}>d+1$), 
and denote by $\zeta$ the limit along an appropriate subsequence 
(not relabelled and implied in what follows).

Let us pass to a subsequence (not relabelled) such that $(h_n)$ converges weakly-$\star$ in ${\rm L}^\infty(\R^{d+1}_+\times\R)$ and denote its limit by $h$.
Note that by passing to the limit $n\to\infty$ in \eqref{kinetic-n}
we obtain 
\begin{align*}
\pa_t h + \Div_\mx \bigl( {f}(t,\mx,\lambda)
h\bigr)-D^2\cdot\bigl( a(\lambda) h\bigr) = \pa_\lambda
\zeta(t,\mx,\lambda)\, ,
\end{align*}
where we have used that $(f_n)$ converges strongly to $f=\partial_\lambda \mff$
(see \eqref{f-conv}).
Thus, by subtracting this equation from \eqref{kinetic-n} we get
\begin{equation}
\label{kinetic-n1}
\begin{split}
\pa_t w_n + \Div & \bigl( {f}(t,\mx,\lambda)
w_n\bigr)-D^2\cdot\bigl( a(\lambda) w_n\bigr)\\
& =\Div \Bigl( \bigl({f}(t,\mx,\lambda)-{f}_n(t,\mx,\lambda)\bigr)
h_n \Bigr)+ \pa_\lambda 
\gamma_n(t,\mx,\lambda) \;,
\end{split}
\end{equation} 
where $w_n=h_n-h$ and $\gamma_n=\zeta_n-\zeta$, and both sequences converge to zero (the first convergence is weak-$\star$ in ${\rm L}^\infty(\R^{d+1}_+\times\R)$, 
and the latter strong in ${\rm L}^{r_0}_{loc}(\R_\lambda; \mathrm{W}_{loc}^{-\frac{1}{2},r_0}(\R^{d+1}_+))$ for any $r_0\in(1,\frac{2d+2}{2d+1})$). 

Due to the boundedness of $(u_n)$, $(w_n)$ is clearly 
uniformly compactly supported with respect to $\lambda$
on $[\alpha,\beta]$ ($h_n$, and thus $h$ as well, are constant functions for $\lambda$
out of $[\alpha,\beta]$).
As we also have \eqref{f-conv}, we see that \eqref{kinetic-n1} clearly satisfies conditions of Corollary \ref{velocity averaging-cor} with $q=\infty$ (see also Remark \ref{rem:app-1-parbolic}).

Therefore, on a subsequence (not relabelled), $(\int_{\R} \rho(\lambda) w_n(t,\mx,\lambda) d\lambda)$ converges strongly to zero in ${\rm L}^1_{loc}(\R^{d+1}_+)$
for any $\rho\in{\rm C}_c(\R)$.
Due to density arguments, we can insert $\rho(\lambda)=\chi_{[\alpha,\beta]}$, obtaining
\begin{align}
2u_n(t,\mx)-\alpha-\beta &= \int_{\alpha}^{u_n(t,\mx)} d\lambda - \int_{u_n(t,\mx)}^{\beta} d\lambda \nonumber \\
&= \int_{\alpha}^{\beta} {\rm sgn}(u_n(t,\mx)-\lambda ) \,d\lambda
\overset{n\to\infty}{\longrightarrow} \int_{\alpha}^{\beta} h(t,\mx,\lambda) \,d\lambda \,, \label{eq:strong_conv_un}
\end{align}
where the latter convergence is in ${\rm L}^1_{loc}(\R^{d+1}_+)$.
Therefore, $(u_n)$ strongly converges in ${\rm L}^1_{loc}(\R^{d+1}_+)$ toward
$$
u(t,\mx) := {\int_{\alpha}^{\beta} h(t,\mx,\lambda) d\lambda\,+{\alpha}+{\beta}\over 2} \,.
$$

The function $u$ is a weak solution to \eqref{d-p}-\eqref{ic}. Indeed, we have for any test function $\varphi \in C_c^2([0,\infty)\times \R^d)$ and a subsequence (not relabelled) $(u_n)$ converging a.e.~to $u$
(which exists due to \eqref{eq:strong_conv_un}):
\begin{align*}
&\iint\limits_{[0,\infty)\times \R^d} \left( u \pa_t \varphi + \mff(t,\mx,u)\cdot \nabla \varphi + A(u)\cdot D^2 \varphi \right) d\mx dt+\int\limits_{\R^d} u_0(\mx) \varphi(0,\mx) d\mx\\
& \qquad =\iint\limits_{[0,\infty)\times \R^d} \left( u_n \pa_t \varphi + \mff_n(t,\mx,u_n)\cdot \nabla \varphi + A(u_n)\cdot D^2 \varphi \right) d\mx dt+\int\limits_{\R^d} u_0(\mx) \varphi(0,\mx) d\mx\\
& \qquad \quad + \iint\limits_{[0,\infty)\times \R^d} \left( (u-u_n) \pa_t \varphi + \big(\mff(t,\mx,u)- \mff_n(t,\mx,u_n)\big) \cdot \nabla \varphi + \big( A(u)-A(u_n) \big)\cdot D^2 \varphi \right) d\mx dt\\
&\qquad =\iint\limits_{[0,\infty)\times \R^d} \left( (u-u_n) \pa_t \varphi + \big(\mff(t,\mx,u)- \mff_n(t,\mx,u_n)\big) \cdot \nabla \varphi + \big( A(u)-A(u_n) \big)\cdot D^2 \varphi \right) d\mx dt \,,
\end{align*} 
since $u_n$ is a weak solution to \eqref{d-p-r}-\eqref{i-c-r}. Moreover,  
\begin{align*}
\|f(\cdot,\cdot, & u(\cdot,\cdot)) - f_n (\cdot,\cdot,u_n(\cdot,\cdot))\|_{\mathrm{L}^p(\supp \varphi)} \\
&\leq \|f(\cdot,\cdot,u(\cdot,\cdot)) - f(\cdot,\cdot,u_n(\cdot,\cdot))\|_{\mathrm{L}^p(\supp \varphi)}
	+ \|f(\cdot,\cdot,u_n(\cdot,\cdot)) - f_n(\cdot,\cdot,u_n(\cdot,\cdot))\|_{\mathrm{L}^p(\supp \varphi)} \\
& \leq \|f(\cdot,\cdot,u(\cdot,\cdot)) - f(\cdot,\cdot,u_n(\cdot,\cdot))\|_{\mathrm{L}^p(\supp \varphi)}
+ \Bigl\|\sup_{\lambda\in[\alpha,\beta]}\bigl|f(\cdot,\cdot,\lambda) - f_n(\cdot,\cdot,\lambda)\bigr|\Bigr\|_{\mathrm{L}^p(\supp \varphi)} \,,
\end{align*}
and the right hand side tends to zero on the limit $n\to\infty$
by the Lebesgue dominated convergence theorem (note that
$(u_n)$ converges a.e.~to $u$ and that $f$ is continuous with respect to 
$\lambda$) and \eqref{f-conv}.
Therefore, after letting $n\to \infty$ in the above integral relation, 
we reach to the conclusion of the theorem.
\end{proof}

\section{Existence of traces for quasi-solutions to \eqref{d-p}}\label{sec:traces}

In this section we study existence of the strong trace (see Definition
\ref{traces}) for (entropy-like) solutions of \eqref{d-p}.
Since it is known (see \cite{pan_jhdeB} for $A=0$) that the notion of weak solutions
is not sufficient to ensure the existence of strong traces, we shall
study an entropy-like concept of solutions -- quasi-solutions.
Therefore, the obtained result (Theorem \ref{main-traces}) is derived for such solutions.  
Let us first recall the notion of quasi-solution for
\eqref{d-p} (see \cite{pan_jhde} for the hyperbolic conservation law).

\begin{definition}\label{def25}
A measurable function $u$ defined on $\R^+\times \R$ is called a quasi-solution to \eqref{d-p} if
$\mff_k(t,\mx,u), A_{kj}(u)\in {\rm L}^1_{loc}(\R^{d+1}_+)$, $k,j=1,\dots,d$, and for a.e.~$\lambda\in \R$ the Kružkov type entropy equality holds
\begin{align}
\label{e-c} 
\pa_t |u-\lambda|&+\Div\Bigl( {\rm
sgn}(u-\lambda)\bigl(\mff(t,\mx,u)-\mff(t,\mx,\lambda) \bigr)
\Bigr)\\
&-D^2\cdot\Bigl( {\rm sgn}(u-\lambda)
\bigl(A(u)-A(\lambda)\bigr)\Bigr) = -\zeta(t,\mx,\lambda),
\nonumber
\end{align} where $\zeta\in {\rm C}(\R_\lambda;w\star-{\cal
M}(\R^{d+1}_+))$ we call the quasi-entropy defect measure, and
here (and in this whole section) by ${\cal	M}(\R^{d+1}_+)$ we denote
the space of Radon measures on $\R^{d+1}_+$ which are locally finite up 
to the boundary $t=0$. 
\end{definition}

\begin{remark}We note that the quasi-solution is not necessarily a weak solution. However, for a regular flux $\mff$, the measure $\zeta(t,\mx,\lambda)$ can be rewritten in
the form $\zeta(t,\mx,\lambda)=\bar{\zeta}(t,\mx,\lambda)+{\rm
sgn}(u-\lambda)\Div_{\mx} \mff(t,\mx,\lambda)$,
for a measure $\bar{\zeta}$. If $\bar{\zeta}$ is non-negative, then the
quasi-solution $u$ is an entropy solution to \eqref{d-p} (i.e. a weak solution as well). For the uniqueness of such entropy solution, we additionally need the chain rule \cite{CK, CP}.
\end{remark}

\begin{remark}
	As before, we shall sometimes consider $\zeta\in C(\R_\lambda;w\star-{\cal
		M}(\R_+^{d+1}))$ as a measure in $\lambda$ as well.
	More precisely, when we write $d\zeta(t,\mx,\lambda)$ in fact we 
	would think of $d\zeta(\lambda)(t,\mx)d\lambda$, where $d\lambda$ 
	is the Lebesgue measure. 
	Since $\zeta$ is continuous in $\lambda$, the measure
	$\zeta\in {\cal M}(\R_+^{d+1}\times\R)$ is 
	locally finite up to the boundary $t=0$ as well. 
	Indeed, for any compact set $K_2\subseteq\R$, set $\{\zeta(\lambda) : \lambda\in K_2\}$
	is vaguely bounded in ${\cal M}(\R_+^{d+1})$, i.e.~with respect to the weak-$\star$
	topology.
	Thus, for any $T>0$ and any compact set $K_1\subseteq\Rd$ there exists 
	$c_{T,{K_1}}>0$ such that the total variation satisfies $\operatorname{Var}\bigl(\zeta(\lambda)\bigr)\leq c_{T,{K_1}}$
	(see e.g.~\cite[(13.4.2)]{JDieud}).
	Now, $\operatorname{Var}\zeta \leq c_{T,{K_1}} \operatorname{meas}(K_2)$,
	where $\operatorname{meas}(K_2)$ stands for the Lebesgue measure 
	of $K_2\subseteq\R$.
\end{remark}

From the notion of quasi-solution, the following kinetic
formulation can be proved (see also \cite[(4.4)]{TT}).

\begin{theorem}
\label{kinetic}

If function $u$ is a quasi-solution to \eqref{d-p}, then
the function
\begin{equation}
\label{equil} h(t,\mx,\lambda)={\rm sgn}(u(t,\mx)-\lambda)=
 -\pa_\lambda |u(t,\mx)-\lambda|
\end{equation} is a weak solution to the following linear equation:
\begin{align}
\label{k-1} \pa_t h + \Div \bigl( {f}(t,\mx,\lambda)
h\bigl)-D^2\cdot\bigl( a(\lambda) h\bigr) = \pa_\lambda
\zeta(t,\mx,\lambda)\, ,
\end{align}
where ${f}=\mff'_\lambda$ and $a=A'_\lambda$.
\end{theorem}

\begin{proof}
It is enough to find derivative of \eqref{e-c} with respect to $\lambda\in
\R$ to obtain \eqref{k-1}.
\end{proof}

The precise definition of the notion of strong traces is the following.

\begin{definition}
\label{traces}
Let $u\in \mathrm{L}^1_{loc}(\R^{d+1}_+)$. A locally integrable function $u_0$ defined on $\R^{d}$ is called the strong trace of $u$ at $t=0$ if $\operatorname{ess\,lim}_{t\to 0^+}u(t,\cdot)=u_0$ 
in ${\rm L}^1_{loc}(\Rd)$, i.e.~for some set $E\subseteq (0,\infty)$ of full Lebesgue 
measure and any relatively compact set $K\subset \subset \R^{d}$ it holds
\begin{equation*}
\lim\limits_{E\ni t\to 0} \|u(t,\cdot)-u_0\|_{{\rm L}^1(K)}=0 \,.
\end{equation*}
\end{definition}

The following theorem holds.

\begin{theorem}\label{main-traces} 
	Let $\mff\in {\rm C}^{0,1}(\R^{d+1}_+\times\R;\R^{d})$ and  let $A\in {\rm C}^{1,1}(\R;\R^{d\times d})$ be such that 
	there exists $\sigma\in {\rm C}^{0,1}(\R;\R^{d\times d})$ such that for any $\lambda\in\R$ 
	we have $a(\lambda):=A'(\lambda)=\sigma(\lambda)^T\sigma(\lambda)$.
	Moreover, assume that the non-degeneracy condition is satisfied: for any 
	compact $K\subseteq\R$,
	\begin{equation}\label{eq:non-deg-cond}
	\sup_{\mxi\in\mathrm{S}^{d-1}} \operatorname{meas} \Bigl\{\lambda\in  K : \;
	\langle a(\lambda)\mxi\,|\,\mxi\rangle=0\Bigr\}=0 \;,
	\end{equation}
	where $\mathrm{S}^{d-1}$ denotes the unit sphere in $\R^{d}$ centered at the origin. 
	
	Then any bounded quasi-solution $u\in {\rm L}^\infty(\R^{d+1}_+)$ to \eqref{d-p} admits the strong trace at $t=0$, i.e.~there exists $u_0\in \mathrm{L}^\infty(\R^{d})$ such that $$\operatorname{ess\,lim}_{t\to 0^+}u(t,\cdot)=u_0$$
	strongly in $\mathrm{L}^1_{loc}(\Rd)$.
\end{theorem}

The structure of proof can be presented as follows.

\begin{itemize}

\item We prove existence of weak traces.

\item We introduce the scaling $t=\frac{\hat{t}}{m},$ $x_1= y_1 + \frac{\hat{x}_1}{\sqrt{m}},\, x_2= y_2 + \frac{\hat{x}_2}{\sqrt{m}}, \dots, x_d= y_d + \frac{\hat{x}_d}{\sqrt{m}}$, where $\my\in\R^{d}$ is a fixed vector. 

\item We obtain a degenerate parabolic transport equation and apply the velocity averaging result.

\item From the previous item, we conclude the existence of the strong traces.

\end{itemize} 

In accordance with the described strategy, let us first show existence of the weak traces. 

\begin{proposition}
\label{w-t} 
Let $h\in \mathrm{L}^\infty(\R^{d+1}_+\times \R)$ be a
distributional solution to \eqref{k-1} and let us define
\begin{equation}\label{eq:set_E}
\begin{aligned}
E=\Bigl\{t\in \R^+: \; (t,\mx,\lambda)& \text{  is a
Lebesgue point of}\\ & h(t,\mx,\lambda) \text{  for a.e.} \, (\mx,\lambda) \in \R^{d}\times \R \Bigr\} \,.
\end{aligned}
\end{equation}

Then there exists $h_0\in \mathrm{L}^\infty(\R^{d+1})$, such that
$$
h(t,\cdot, \cdot) \rightharpoonup h_0\,,  \mbox{ \rm weakly-$\star$ in } \mathrm{L}^\infty(\R^{d+1}) \,,
	\mbox{ \rm as } t\to 0\,,\; t\in E\,.
	$$
\end{proposition}

\begin{proof}Note first that $E$ is of full measure. Since $h\in \mathrm{L}^\infty(\R^{d+1}_+\times \R)$, the family $\{
h(t,\cdot,\cdot)\}_{t\in E}$ is bounded in $\mathrm{L}^\infty(\R^{d+1})$.
             Due to the weak-$\star$ precompactness of $\mathrm{L}^\infty(\R^{d+1})$, there exists a sequence $\{t_m\}_{m\in \N}$ in $E$ such that $t_m\to 0$ as $m\to \infty$,
             and $h_0\in \mathrm{L}^\infty(\R^{d+1})$, such that
                         \begin{equation}\label{wc}
                                  h(t_m, \cdot,\cdot) \rightharpoonup h_0 \,, \mbox{ weakly-$\star$ in $\mathrm{L}^\infty(\R^{d+1})$, as } m\to \infty \,.
                         \end{equation}
             For $\phi\in \mathrm{C}_c^\infty(\R^{d})$, $\rho\in \mathrm{C}^1_c(\R)$, denote
             $$
               I(t):=\int_{\R^{d+1}}h(t,\mx,\lambda)\rho(\lambda)\phi(\mx)\,d\mx d\lambda \,, \quad t\in E\,.
             $$
             With this notation, \eqref{wc} means that
                         \begin{equation}\label{wc1}
                                \lim_{m\to \infty} I(t_m)= \int_{\R^{d+1}} h_0(\mx,\lambda)\rho(\lambda)\phi(\mx)\, d\mx d\lambda = :I(0) \,.
                         \end{equation}
             Now, fix $\tau\in E$ and notice that for the regularisation $I_\eps=I\star \omega_\eps$, where $\omega_\eps$ is the standard convolution kernel, it holds
             $$
             \lim\limits_{\eps \to 0} I_\eps(\tau)=I(\tau).
             $$
             
Then, fix $m_0\in \N$, such that $E\ni t_m\leq \tau$, for $m\geq m_0$, and remark that
 \begin{equation}
 \begin{split}\nonumber
 I(\tau)- I(t_m) &= \lim\limits_{\eps\to 0}\int_{t_m}^\tau I_\eps'(t)\,dt 
 \\
 &= \sum_{j=1}^d \int_{(t_m,\tau] \times \R^{d+1}} h(t,\mx,\lambda)f_j(t,\mx,\lambda) \rho(\lambda)\pa_{x_j}\phi(\mx)\, dt\, d\mx \,d\lambda \\
 &\qquad - \sum_{j,k=1}^d \int_{(t_m,\tau] \times \R^{d+1}} h(t,\mx,\lambda) a_{jk}(\lambda) \rho(\lambda) \pa_{x_j x_k}\phi(\mx)\, dt\, d\mx \,d\lambda \\
 &\qquad - \int_{(t_m,\tau] \times \R^{d+1}} \phi(\mx) \rho'(\lambda)\, d\zeta(t,\mx,\la) \,,
 \end{split}
 \end{equation}
 where we have used that $h$ is a distributional solution to \eqref{k-1}.
Hence, passing to the limit as $m\to\infty$, and having in mind \eqref{wc1} and the fact that
$\zeta$ is locally finite up to the boundary $t=0$, we obtain
  \begin{equation}\begin{split}\nonumber
  I(\tau)- I(0) &= \sum_{j=1}^d \int_{(0,\tau] \times \R^{d+1}} h(t,\mx,\lambda) f_j(t,\mx,\lambda) \rho(\lambda)\pa_{x_j}\phi(\mx) \, dt\, d\mx \,d\lambda \\
  &\qquad - \sum_{j,k=1}^d \int_{(0,\tau] \times \R^{d+1}} h(t,\mx,\lambda) a_{jk}(\lambda) \rho(\lambda)  \pa_{x_j x_k}\phi(\mx)\, dt\, d\mx \,d\lambda \\
  &\qquad - \int_{(0,\tau] \times \R^{d+1}} \phi(\mx) \rho'(\lambda) \, d\zeta (t,\mx,\la) \,.
  \end{split}\end{equation}
  The right hand side clearly tends to zero as $\tau\to 0$.
             Thus, for all $\phi\in \mathrm{C}_c^\infty(\R^{d+1})$ and $\rho\in \mathrm{C}^1_c(\R)$
             we have $\lim_{E\ni\tau\to 0}I(\tau)=I(0)$, i.e.
             $$
             \lim_{E\ni\tau\to 0 } \int_{\R^{d+1}} h(\tau,\mx,\lambda)\rho(\lambda) \phi(\mx)\, d\mx d\lambda = \int_{\R^{d+1}} h_0(\mx,\lambda)\rho(\lambda) \phi(\mx) \, d\mx d\la \,.
             $$
Having in mind that $h(\tau, \cdot,\cdot)$, $\tau\in E$, 
is bounded, and that $\mathrm{C}_c^\infty(\R^{d+1})$ is dense in
$\mathrm{L}^1(\R^{d+1})$, we complete the proof.
\end{proof}

\begin{remark}\label{rem:w-t_u}
If $u$ is a bounded quasi-solution to \eqref{d-p}, then in 
\cite[Corollary 2.2]{pan_jhde} was proved that it admits the weak trace.
The same conclusion can be derived from the previous proposition. 
 
Indeed, let $M>0$ be such that 
$$
|u(t,\mx)| \leq M \;, \quad \hbox{a.e} \ (t,\mx)\in\R^{d+1}_+ \;.
$$
Then, by the definition of $h$ (it is the sign function; see \eqref{equil})
$$
\int_{-M}^{M} h(t,\mx,\lambda)\,d\lambda 
	= \int_{-M}^{M} {\rm sgn}\bigl(u(t,\mx)-\lambda\bigr)\,d\lambda
	= \int_{-M}^{u(t,\mx)} d\lambda - \int_{u(t,\mx)}^{M} d\lambda
	= 2u(t,\mx) \,.
$$
Thus, the claim follows by Proposition \ref{w-t} by noting that the characteristic function 
$\chi_{[-M,M]}$ of the interval $[-M,M]$ is in $\mathrm{L}^1(\R)$.
More precisely, we have
\begin{equation}\label{eq:u_w-t}
u(t,\cdot) \overset{\star}{\rightharpoonup} u_0:=\frac{1}{2}\int_{-M}^{M} h_0(\cdot,\lambda)\,d\lambda
\end{equation}
weakly-$\star$ in $\mathrm{L}^\infty(\R^{d})$ as $t\to 0$, $t\in E$.

Moreover, for $\lambda\mapsto\lambda\chi_{[-M,M]}(\lambda)\in \mathrm{L}^1(\R)$ we have
$$
\int_{-M}^{M} \lambda h(t,\mx,\lambda)\,d\lambda 
	= \int_{-M}^{u(t,\mx)} \lambda \,d\lambda - \int_{u(t,\mx)}^{M} \lambda\,d\lambda
	= u(t,\mx)^2 -M^2 \,.
$$
Therefore, one can similarly conclude that
\begin{equation}\label{eq:uu_w-t}
u^2(t,\cdot) \overset{\star}{\rightharpoonup} u_1:=\int_{-M}^{M} \lambda h_0(\cdot,\lambda)
	\,d\lambda +M^2
\end{equation}
weakly-$\star$ in $\mathrm{L}^\infty(\R^{d})$ as $t\to 0$, $t\in E$.

If one can get that $u_1=u_0^2$, the strong trace can be obtained from the above weak convergences. 
Namely, for an arbitrary $\varphi\in \mathrm{C}_c(\Rd)$ by \eqref{eq:u_w-t}-\eqref{eq:uu_w-t}
we have
\begin{equation}\label{eq:u_strong_conv}
\begin{aligned}
\lim_{E\ni t\to 0}\int_\Rd \Bigl(u(t,\mx)- & u_0(\mx)\Bigr)^2 \varphi(\mx) \,d\mx \\
&= \lim_{E\ni t\to 0}\int_\Rd \Bigl(u(t,\mx)^2-2u(t,\mx)u_0(\mx)+u_0(\mx)^2\Bigr)\varphi(\mx) 
\,d\mx \\
&= \int_\Rd \Bigl(u_1(\mx)-u_0(\mx)^2\Bigr)\varphi(\mx) \,d\mx \,.
\end{aligned}
\end{equation}

A sufficient condition for $u_1=u_0^2$ in terms of a certain strong convergence 
of rescaled (sub)sequences is given in Proposition \ref{equivalence} below.
\end{remark}

Now we use the rescaling procedure (or the so-called blow-up method) in order to 
obtain a sufficient condition for the existence of the strong trace. 
More precisely, let us change the variables in \eqref{k-1} in the following way:
 $t=\frac{\hat{t}}{m},$ $x_1= y_1 + \frac{\hat{x}_1}{\sqrt{m}},\, x_2= y_2 + \frac{\hat{x}_2}{\sqrt{m}}, \dots, x_d= y_d + \frac{\hat{x}_d}{\sqrt{m}}$, i.e.
\begin{equation}
\label{blow-up}
(t,\mx,\la)= \Bigl(\frac{\hat t}{m}, \frac{\hat \mx}{\sqrt{m}}+\my, \la\Bigr) \,,
\end{equation} 
where $\my\in\R^{d}$ is a fixed vector. We get that a rescaled solution to \eqref{k-1}, denoted by 
$$
h_m(\hat{t},\hat{\mx},\my,\lambda):=h\Bigl(\frac{\hat{t}}{m},\frac{\hat{\mx}}{\sqrt{m}}+\my,\lambda\Bigr) \,
$$
satisfies       
\begin{align}
\label{lhm}
L_{h_m}:=&  \left(\pa_{\hat t} h_m +\frac{1}{\sqrt{m}}\sum_{k=1}^d \pa_{\hat x_k} 
	\left( f_k^m h_m \right) \right)-
	\sum_{k,j=1}^d \pa^2_{\hat x_j \hat x_k} (a_{jk}  h_m) = \frac{1}{m} 
	\pa_\lambda  \zeta_m^\my \,, \\
\label{ic-w}
h_m{\big|}_{t=0}=&h_0\bigl(\frac{\hat \mx}{\sqrt{m}}+\my,\lambda\bigr) \,,
\end{align}  where the initial conditions are understood in the weak sense, 
$h_0$ is the weak trace from Proposition \ref{w-t},
and $f_k^m(\hat{t},\hat{\mx},\lambda) 
= f_k\bigl(\frac{\hat{t}}{m},\frac{\hat{\mx}}{\sqrt{m}}+\my,\lambda\bigr)$.

Let us remark that the equality between $\zeta$ and $\zeta_m^\my$ is understood in the sense of distributions:
\begin{equation}
\label{**}
\langle \zeta_m^\my,\varphi \rangle=m^{d/2+1} 
\int_{\R^{d+1}_+\times\R}\varphi\bigl(m\, t,\sqrt{m}\, ({\mx}-{\my}),\lambda\bigr)\,d\zeta({t},{\mx},\lambda) \,.
\end{equation} 
If we prove that the sequence
\begin{equation}
\label{scal}      
\int_\R h\Bigl(\frac{\hat  t}{m},\frac{\hat \mx}{\sqrt{m}}+\my,\lambda\Bigr)\rho(\lambda)d\lambda
	\,, \ \ m\in \N \,,
\end{equation} converges strongly in $\mathrm{L}^1_{loc}(\R^{d+1}_+\times \R^{d})$ along a subsequence, we will obtain that function $u$ admits the trace in the sense of Definition \ref{traces}. More precisely, the following proposition holds.

\begin{proposition}
\label{equivalence}
Let $u$ be a bounded quasi-solution to \eqref{d-p} and let $h$ be given by \eqref{equil}.
Assume that for every $\rho \in \mathrm{C}^1_c(\R)$ the sequence given by \eqref{scal} converges toward 
$\int_\R h_0(\my,\lambda)\rho(\lambda)d\lambda$ in $\mathrm{L}^1_{loc}(\R^{d+1}_+\times \R_{\my}^d)$ along a subsequence, where $h_0$ is the weak trace of $h$ (see Proposition \ref{w-t}).

Then, the function $u$ admits the strong trace at $t=0$ and it is equal to 
$$
u_0(\mx):=\frac{1}{2}\int_{-M}^M h_0(\mx,\lambda) \,d\lambda \,,
$$
where $M=\|u\|_{\mathrm{L}^\infty(\R^{d+1}_+)}$. 
\end{proposition}
\begin{proof}
Since both $h$ and $h_0$ are bounded, using the density arguments we conclude that if the (sub)sequence from \eqref{scal} converges in $\mathrm{L}^1_{loc}(\R^{d+1}_+\times \R^{d})$ for any $\rho\in \mathrm{C}^1_c(\R)$, then it will also converge for any $\rho\in \mathrm{L}^1(\R)$. Let us take $\rho=\chi_{[-M,M]}$, where $\chi_{[-M,M]}$ is the characteristic function of the interval $[-M,M]$, and $M>0$ is such that
$$
|u(t,\mx)| \leq M\,, \quad \hbox{a.e.} \ (t,\mx)\in\R^{d+1}_+ \;.
$$

Thus, for any non-negative $\varphi\in \mathrm{C}_c(\R^{d+1}_+\times \R^{d})$, it holds (along the subsequence from the formulation of the proposition)
\begin{align*}
\lim\limits_{m\to \infty}\int_{\Rd\times\R^{d+1}_+}
	\varphi(\hat{t},\hat{\mx},\my)
	\biggl|\int_{-M}^M\biggl(
	h\Bigl(\frac{ \hat  t}{m},\frac{\hat \mx}{\sqrt{m}}+\my,\lambda\Bigr)
	-h_0(\my,\lambda)
	\biggr) \,d\lambda \biggr|d\my d\hat{\mx} d \hat{t}=0 \,.
\end{align*}  
Using the definition of the function $h$ (it is the sign function; see \eqref{equil})
we have
\begin{align*}
\int_{-M}^M h\Bigl(\frac{ \hat  t}{m},\frac{\hat \mx}{\sqrt{m}}+\my,\lambda\Bigr) \,d\lambda
	&=\int_{-M}^M {\rm sgn}\biggl(u\Bigl(\frac{\hat t}{m},
	\frac{\hat \mx}{\sqrt{m}}+\my\Bigr) - \lambda \biggr) d\lambda \\
&= \int_{-M}^{u(\frac{\hat t}{m},\frac{\hat \mx}{\sqrt{m}}+\my)} d\lambda
	- \int_{u(\frac{\hat t}{m},\frac{\hat \mx}{\sqrt{m}}+\my)}^{M} d\lambda \\
&= 2u\Bigl(\frac{\hat t}{m},\frac{\hat \mx}{\sqrt{m}}+\my\Bigr) \;.
\end{align*}
Taking this into account and the change of variables $\mz=\frac{\hat \mx}{\sqrt{m}}+\my$ 
(with respect to $\my$), the previous limit reads 
\begin{align*}
\lim\limits_{m\to\infty}\int_{\R^{d}\times\R^{d+1}_+}
	\varphi\Bigl(\frac{\hat{t}}{m},\hat{\mx},\mz-\frac{\hat \mx}{\sqrt{m}}\Bigr)
	\biggl|2 u\Bigl(\frac{\hat t}{m},\mz\Bigr)
	-\int_{-M}^M h_0\Bigl(\mz-\frac{\hat{\mx}}{\sqrt{m}},\lambda\Bigr) \,d\lambda
	\biggr| \,d\mz d\hat{\mx} d\hat{t}=0 \;.
\end{align*}
Furthermore, the limit still holds if we replace 
$\varphi(\hat{t},\hat{\mx},\mz-\frac{\hat \mx}{\sqrt{m}})$ by 
$\varphi(\hat{t},\hat{\mx},\mz)$ 
and $h_0(\mz-\frac{\hat{\mx}}{\sqrt{m}},\lambda)$ by 
$h_0(\mz,\lambda)$, i.e.
\begin{align}
\label{****}
\lim\limits_{m\to\infty}\int_{\R^{d}\times\R^{d+1}_+}
\varphi(\hat{t},\hat{\mx},\mz)
\biggl|2 u\Bigl(\frac{\hat t}{m},\mz\Bigr)
-\int_{-M}^M h_0(\mz,\lambda) \,d\lambda
\biggr| \,d\mz d\hat{\mx} d\hat{t}=0 \;.
\end{align}
Indeed, the first replacement is justified since 
$$
(\hat t, \hat \mx,\mz)\mapsto \biggl|2 u\Bigl(\frac{\hat t}{m},\mz\Bigr)
	-\int_{-M}^M h_0\Bigl(\mz-\frac{\hat{\mx}}{\sqrt{m}},\lambda\Bigr) \,d\lambda
	\biggr|
$$
is bounded and $\varphi$ is a continuous function with compact support, hence 
the convergence 
$$
\lim_{m\to\infty} \varphi\Bigl(\hat{t},\hat{\mx},\mz-\frac{\hat \mx}{\sqrt{m}}\Bigr)
	= \varphi(\hat{t},\hat{\mx},\mz)
$$
is uniform in $(\hat t, \hat \mx, \mz)$.
The second one follows by the convergence (implied by the continuity
of the average and the Lebesgue dominated convergence theorem)
$$
\lim_{m\to\infty} \int_{-M}^M \biggl|h_0\Bigl(\mz-\frac{\hat{\mx}}{\sqrt{m}},\lambda\Bigr) 
	- h_0(\mz,\lambda) \biggr| \,d\lambda = 0
$$ in $\mathrm{L}^1_{loc}(\Rd\times\R^{d+1}_+)$. 

Therefore, due to arbitrariness of $\varphi$ in \eqref{****}, we conclude 
$$
u\Bigl(\frac{\hat{t}}{m},\mz\Bigr)\to \frac{1}{2}\int_{-M}^M h_0(\mz,\lambda)\,d\lambda, \ \ m\to \infty
$$ 
in $\mathrm{L}^1_{loc}(\R^{d+1}_+)$ along the subsequence from the formulation of the proposition. This means that (for another subsequence not relabelled) 
there exists $\hat{E}\subseteq\R^+$ of full measure such that for any $\hat{t}\in\hat{E}$
we have
\begin{equation}
\label{sc1}
u\Bigl(\frac{\hat{t}}{m},\mz\Bigr)\to \frac{1}{2}\int_{-M}^M h_0(\mz,\lambda) \,d\lambda=u_0(\mz), \ \ m\to \infty
\end{equation}  in $\mathrm{L}^1_{loc}(\R^{d})$. 
It is easy to see that for $E$ given by \eqref{eq:set_E} set 
$E_\infty:=\bigcap\limits_{m\in\N} {m} E$ is of full measure
(since $E$ is of full measure). 
Thus, the intersection $\hat{E}\cap E_\infty$ is non-empty (in fact it is a set of full 
measure as well), so we can choose $\hat{t}\in\R^+$ such that \eqref{sc1} holds and 
$\frac{\hat{t}}{m}\in E$, $m\in\N$.

Now, choose $\rho(\lambda)=\lambda \,\chi_{[-M,M]}(\lambda)$ where $\chi_{[-M,M]}(\lambda)$ is the characteristic function of the interval $[-M,M]$. It holds according to Proposition \ref{w-t} (see also Remark \ref{rem:w-t_u})
$$
u^2(t,\mx) = \int_{-M}^M \lambda h(t,\mx,\lambda) \,d\lambda + M^2
	\overset{\star}{\rightharpoonup} \int_{-M}^M \lambda h_0(\mx,\lambda) \,d\lambda +M^2 
	=: u_1(\mx) 
$$ 
in $\mathrm{L}^\infty(\R^{d})$ as  $E\ni t\to 0$. 
Since the weak-$\star$ convergence in $\mathrm{L}^\infty(\Rd)$ implies the weak convergence in $\mathrm{L}^1_{loc}(\Rd)$, and since weak and strong limits coincide, 
from here and \eqref{sc1} we see that it must be $u_1=u_0^2$. 
Finally, by \eqref{eq:u_strong_conv} (see Remark \ref{rem:w-t_u})
we have
$$
u(t,\cdot) \rightarrow u_0 
$$
in $\mathrm{L}^2_{loc}(\Rd)$ as $t\to 0$, $t\in E$, which implies the convergence in $\mathrm{L}^1_{loc}$.
Hence, $u_0$ is the strong trace. \\
\end{proof}

Having the last proposition in mind, we clearly need the following theorem.

\begin{theorem}
\label{sc-ndeg}
Under assumption of Theorem \ref{main-traces}, let $h$ be given by \eqref{equil},
and let $h_0$ be the weak trace of $h$ (see Proposition \ref{w-t}).

Then, for any $\rho\in \mathrm{C}^1_c(\R)$, the sequence of functions
$$
(t,\mx,\my)\mapsto \int_{\R}  h\Bigl(\frac{t}{m},\frac{\mx}{\sqrt{m}}+\my,
	\lambda\Bigr) \rho(\lambda) \,d\lambda
$$ 
converges to $\int_\R h_0(\my,\lambda)\rho(\lambda)d\lambda$
in $\mathrm{L}^1_{loc}(\R^{d+1}_+\times \R^{d})$.
\end{theorem}
\begin{proof}
First, notice that for every $\my\in \R^d$, the sequence of function $(h_m)$ satisfies diffusive transport equation \eqref{lhm} which can be rewritten in the form
$$
\pa_{\hat t} h_m -	\sum_{k,j=1}^d \pa^2_{\hat x_j \hat x_k} (a_{jk}  h_m) =-\frac{1}{\sqrt{m}}\sum_{k=1}^d \pa_{\hat x_k} 
	\left( f_k^m h_m \right)+ \frac{1}{m} 
	\pa_\lambda  \zeta_m^\my \;.
$$ 
Clearly, $\frac{1}{\sqrt{m}} 
	\left( f_k^m h_m \right)$ converges strongly to zero in $\mathrm{L}^2_{loc}(\R^{d+1}_+\times\R)$. Moreover, $\frac{1}{m}  \zeta_m^\my$ converges to zero in ${\cal M}(\R^{d+1}_+\times \R)$ (this is proved in the same way as \cite[Lemma 2]{vass} or \cite[Lemma 3.2]{pan_jhdeB}). Therefore, keeping in mind conditions \eqref{eq:non-deg-cond}, we can apply Theorem \ref{velocity averaging} (see also Remark \ref{rem:only_bdd}) to conclude that for every $\rho\in \mathrm{C}^1_c(\R)$ there exists a subsequence of $(\int h_m(\cdot,\lambda)\rho(\lambda) d\lambda)$ (not relabelled) strongly converging in $\mathrm{L}^1_{loc}(\R^{d+1}_+)$ toward say $\int \tilde{h}(t,\mx,\my,\lambda) \rho(\lambda) d \lambda$. We note that $\tilde{h}(t,\mx,\my,\lambda)$ does not depend on $\rho$ since it is a weak limit of $(h_m)$ in $\mathrm{L}_{loc}^2(\R^{d+1}_+\times \R)$ along appropriate subsequence.  
	
On the other hand, the sequence $(h_m)$ satisfies \eqref{lhm}-\eqref{ic-w} and thus $\tilde{h}$ satisfies the Cauchy problem
\begin{align}
\label{lim-1}
&\pa_{\hat t} \tilde{h} -	\sum_{k,j=1}^d \pa^2_{\hat x_j \hat x_k} (a_{jk}  \tilde{h}) = 0\\
\label{lim-ic}
&\tilde{h}\Big|_{t=0}=h_0(\my,\lambda)\,,
\end{align} which implies $\tilde{h}\equiv h_0(\my,\lambda)$ since the solution of the latter Cauchy problem is unique. Moreover, from here it follows that the entire sequence $(\int h_m(\cdot,\lambda)\rho(\lambda) d\lambda)$ must converge toward $(\int h_0(\cdot,\lambda)\rho(\lambda) d\lambda)$ (since every converging subsequence must converge toward the solution to \eqref{lim-1}-\eqref{lim-ic}).
\end{proof}

We have thus proved that $(h_m)$ satisfies conditions of Proposition \ref{equivalence} and this in turn implies that the quasi-solution to \eqref{d-p} indeed admits existence of strong traces.
Thus, with this the proof of Theorem \ref{main-traces} is completed.

\section{Appendix}\label{sec:appendix}

In this section, we provide some auxiliary statements that we use in the proof of the velocity averaging result as well as the proof of Lemma \ref{TVRDNJA}. 

\begin{lemma}
\label{pomocna}
Assume $d\geq 2$ and that a matrix function $a$ satisfies conditions of Lemma \ref{multiplierlemma1}. Then

\begin{itemize}

\item[(i)] For a.e.~$\lambda \in S$ and $r\in (1,d)$ the Fourier multiplier operator ${\cal A}_{\frac{1}{|\cdot|+\langle a(\lambda)\cdot \,|\,\cdot \rangle}}$ 
is bounded operator from ${\rm L}^r(\R^d)$ to ${\rm L}^{r^*}(\R^d)$ uniformly with respect to $\lambda \in S$, where $r^*=\frac{dr}{d-r}$. 

\item[(ii)] For a.e.~$\lambda \in S$ and any $k\in\{1,2,\dots, d\}$, $r\in (1,\infty)$, the operator 
$\pa_{x_k}{\cal A}_{\frac{1}{|\cdot|+\langle a(\lambda)\cdot \,|\,\cdot \rangle}}$ is continuous operator from ${\rm L}^r(\R^d)$ to ${\rm L}^r(\R^d)$ 
uniformly with respect to $\lambda \in S$.

\end{itemize}

\end{lemma}
\begin{proof}
\begin{itemize}

\item[(i)] It is well known that the Riesz potential ${\cal A}_{\frac{1}{|\mxi|}}$ is continuous mapping from ${\rm L}^r(\R^d)$ to ${\rm L}^{r^*}(\R^d)$ 
for $r^*=\frac{dr}{d-r}$ \cite[Section 5]{stein}. Moreover, the operator ${\cal A}_{\frac{|\mxi|}{|\mxi|+\langle a(\lambda)\mxi \,|\,\mxi \rangle}}$ 
is continuous operator from ${\rm L}^r(\R^d)$ to ${\rm L}^r(\R^d)$ by Lemma \ref{TVRDNJA} (taking into account the invariance of multipliers under 
orthogonal transformations; see the proof of Lemma \ref{multiplierlemma1}). Now, the statement follows by:
$$
{\cal A}_{\frac{1}{|\mxi|+\langle a(\lambda)\mxi \,|\,\mxi \rangle}}={\cal A}_{\frac{1}{|\mxi|}}\circ {\cal A}_{\frac{|\mxi|}{|\mxi|+\langle a(\lambda)\mxi \,|\,\mxi \rangle}}.
$$  

\item[(ii)] The second part follows by applying Lemma \ref{multiplierlemma1} on $\psi(\mxi)=2\pi i \xi_k$ since

$$
\pa_{x_k}{\cal A}_{\frac{1}{|\mxi|+\langle a(\lambda)\mxi \,|\,\mxi \rangle}}={\cal A}_{\frac{2 \pi i \xi_k}{|\mxi|+\langle a(\lambda)\mxi \,|\,\mxi \rangle}} \,.
$$
\vskip-0,8cm
\end{itemize}
\end{proof}

Now, we are going to prove Lemma \ref{TVRDNJA}, but first we develop two 
auxiliary results that are given in the following two lemmata.

\begin{lemma}
\label{tvrdnja'} Let $\mkappa=(\kappa_1,\kappa_2,\dots,\kappa_d) \in [0,\infty)^d$ and define $f:\Rd\to\R$ by 
	$$
	f(\mxi)=\frac{|\mxi|}{|\mxi|+\sum\limits_{j=1}^d \kappa_j \xi_j^2} \;.
	$$
	For every multi-index $\malpha \in \N_0^d$ and $\mxi\in\Rd\setminus\{0\}$, it holds
\begin{equation*} 
(\pa^\malpha f)(\mxi)=\biggl(|\mxi|+\sum\limits_{j=1}^d \kappa_j \xi_j^2 \biggr)^{-|\malpha |-1} P_{\malpha}\Bigl(\mkappa,\mxi, \frac{1}{|\mxi|}\Bigr) \,,
\end{equation*}
where $P_{\malpha}(\mkappa,\mxi, \eta)$ is a polynomial consisting of the terms $C \mkappa^\mbeta \mxi^\mgamma \eta^l$ for a constant $C=C(\malpha,d)$, 
multi-indices $\mbeta, \mgamma \in \N^d_0$ and $ l\in \N_0$, such that
\begin{align}
\label{cond-app}
\alpha_j+\gamma_j\geq 2\beta_j\,, \quad j=1,\dots,d\,; \qquad |\malpha| \geq |\mbeta |\,; \qquad |\mgamma|=1+l+|\mbeta|\,.
\end{align} 
\end{lemma}
\begin{proof}
We prove the claim by the induction argument with respect to the order of derivative $n=|\malpha|$.

For $n=0$ ($|\malpha|=0$), we have $P_0=|\mxi|=\frac{|\mxi|^2}{|\mxi|}=\sum\limits_{j=1}^d \xi_j^2 \frac{1}{|\mxi|^2}$. It is easy to check that conditions \eqref{cond-app} are satisfied.

Assume now that the statement holds for some $n\in \N_0$ and let us prove that it holds for $n+1$. 
Let $\tilde{\malpha}\in \N_0^d$, $|\tilde{\malpha}|=n+1$, and let $s\in \{1,2,\dots,d \}$ and $\malpha \in \N_0^d$, $|\malpha|=n$, 
be such that $\tilde{\malpha}=\mathbf{e}_s+\malpha$, where $\mathbf{e}_s$ is the $s$-th vector of the canonical basis of $\Rd$. 
By the Schwarz rule and the assumption of the induction argument we have
\begin{align*}
&\pa^{\tilde{\malpha}} f(\mxi)=\pa_{\xi_s}\Biggl(\biggl(|\mxi|+\sum\limits_{j=1}^d \kappa_j \xi_j^2 \biggr)^{-|\malpha |-1} 
	P_{\malpha}\Bigl(\mkappa,\mxi, \frac{1}{|\mxi|}\Bigr) \Biggr) \\
&=(-|\malpha|-1)\biggl(|\mxi|+\sum\limits_{j=1}^d\kappa_j \xi_j^2 \biggr)^{-|\malpha|-2} \left(\frac{\xi_s}{|\mxi|}+2\kappa_s \xi_s \right) 
	P_{\malpha}\Bigl(\mkappa,\mxi, \frac{1}{|\mxi|}\Bigr) \\
&\quad +\biggl(|\mxi|+\sum\limits_{j=1}^d\kappa_j \xi_j^2\biggr)^{-|\malpha|+1} \!\!\!(\pa_{\xi_s} P_\malpha)
	\Bigl(\mkappa,\mxi,\frac{1}{|\mxi|}\Bigr)-\biggl(|\mxi|+\sum\limits_{j=1}^d\kappa_j \xi_j^2\biggr)^{-|\malpha|+1} \!\!\!(\pa_\eta P_\malpha)
	\Bigl(\mkappa,\mxi,\frac{1}{|\mxi|}\Bigr) \frac{\xi_s}{|\mxi|^3}\\
&=\biggl(|\mxi|+\sum\limits_{j=1}^d\kappa_j \xi_j^2 \biggr)^{-|\tilde{\malpha}|-1} \Biggl[ -|\tilde\malpha|\left(\frac{\xi_s}{|\mxi|}
	+2\kappa_s \xi_s \right)P_\malpha\Bigl(\mkappa,\mxi,\frac{1}{|\mxi|}\Bigr)\\
&\quad +\biggl( \frac{|\mxi|^2}{|\mxi|}+ \sum\limits_{j=1}^d\kappa_j \xi_j^2 \biggr) (\pa_{\xi_s}P_{\malpha})\Bigl(\mkappa,\mxi,\frac{1}{|\mxi|}\Bigr)
	- \biggl( \frac{\xi_s}{|\mxi|^2}+ \sum\limits_{j=1}^d\kappa_j \frac{\xi_s \xi_j^2}{|\mxi|^3} \biggr)
	(\pa_{\eta}P_{\malpha})\Bigl(\mkappa,\mxi,\frac{1}{|\mxi|}\Bigr) \Biggr] \\
&=:\biggl(|\mxi|+\sum\limits_{j=1}^d\kappa_j \xi_j^2 \biggr)^{-|\tilde{\malpha}|-1} P_{\tilde{\malpha}}\Bigl(\mkappa,\mxi,\frac{1}{|\mxi|}\Bigr) \;.
\end{align*}From here, a direct analysis of the six terms forming $P_{\tilde{\malpha}}(\mkappa,\mxi,\eta)$ provides \eqref{cond-app}. \end{proof}

Analogously one can prove the following result. 

\begin{lemma}
	\label{tvrdnja''} 
	Let $\mkappa=(\kappa_1,\kappa_2,\dots,\kappa_d) \in [0,\infty)^d$ 
	and $m\in\{1,2,\dots,d\}$, and define $g:\Rd\to\R$ by 
	$$
	g(\mxi)=\frac{\kappa_m\xi^2_m}{|\mxi|+\sum\limits_{j=1}^d \kappa_j \xi_j^2} \;.
	$$
	For every multi-index $\malpha \in \N_0^d$ and $\mxi\in\Rd\setminus\{0\}$, it holds
	\begin{equation*} 
	(\pa^\malpha g)(\mxi)=\biggl(|\mxi|+\sum\limits_{j=1}^d \kappa_j \xi_j^2 \biggr)^{-|\malpha |-1} P_{\malpha}\Bigl(\mkappa,\mxi, \frac{1}{|\mxi|}\Bigr) \,,
	\end{equation*}
	where $P_{\malpha}(\mkappa,\mxi, \eta)$ is a polynomial consisting of the terms $C \mkappa^\mbeta \mxi^\mgamma \eta^l$ for a constant $C=C(\malpha,d,m)$, 
	multi-indices $\mbeta, \mgamma \in \N^d_0$ and $ l\in \N_0$, such that
	\begin{align}
	\label{cond-app-g}
	\alpha_j+\gamma_j\geq 2\beta_j\,, \quad j=1,\dots,d\,; \qquad |\malpha|+1 \geq |\mbeta |\,; \qquad |\mgamma|=1+l+|\mbeta|\,; \qquad \beta_m\geq 1\,.
	\end{align} 
\end{lemma}
Now, we can prove Lemma \ref{TVRDNJA}.

\vspace{0.5cm}

\noindent {\bf Proof of Lemma \ref{TVRDNJA}:} Since the space of $\mathrm{L}^p$-Fourier
multipliers is an algebra, it is sufficient to prove the statement for $s\in[0,1]$. 
For $s=0$ the claim trivially holds, so let us first consider $s=1$. 

We use the Marcinkiewicz theorem (Theorem \ref{m1}). Let $\malpha \in \N^d_0$ and $\mxi\in\Rd\setminus\{0\}$. 
By the previous lemmata, for both functions $f$ and $g$ it is sufficient to estimate
$$
\mxi^\malpha \mkappa^\mbeta \mxi^\mgamma\frac{1}{|\mxi|^l} \biggl( |\mxi|+\sum\limits_{j=1}^d {\kappa_j\xi^2_j} \biggr)^{-|\malpha|-1} \;,
$$
where $\mbeta,\mgamma$ and $l$ satisfy \eqref{cond-app} and
\eqref{cond-app-g}, respectively. Thus, we have
\begin{align*}
\Biggl| \mxi^\malpha \mkappa^\mbeta \mxi^\mgamma\frac{1}{|\mxi|^l} & \biggl( |\mxi|+\sum\limits_{j=1}^d 
	{\kappa_j\xi^2_j} \biggr)^{-|\malpha|-1} \Biggr|
	= \prod_{j=1}^d (\kappa_j \xi_j^2)^{\beta_j} \,  \prod_{j=1}^d |\xi_j|^{\alpha_j+\gamma_j-2\beta_j} \, \frac{1}{|\mxi|^l} 
	\frac{1}{\Bigl(|\mxi|+\sum\limits_{j=1}^d {\kappa_j\xi^2_j}\Bigr)^{|\malpha|+1}} \\
&\leq \frac{|\mxi|^{|\malpha|+|\mgamma|-2|\mbeta|-l}}{\Bigl(|\mxi|+\sum\limits_{j=1}^d {\kappa_j\xi^2_j}\Bigr)^{|\malpha|+1-|\mbeta|}}
= \frac{|\mxi|^{|\malpha|+1-|\mbeta|}}{\Bigl(|\mxi|+\sum\limits_{j=1}^d {\kappa_j\xi^2_j}\Bigr)^{|\malpha|+1-|\mbeta|}} \leq 1 \;,
\end{align*} 
where in the first inequality we have used $|\xi_j|^{\alpha_j+\gamma_j-2\beta_j}\leq |\mxi|^{\alpha_j+\gamma_j-2\beta_j}$ as 
$\alpha_j+\gamma_j-2\beta_j\geq 0$ by \eqref{cond-app} and \eqref{cond-app-g}, 
while the last inequality is trivial since $|\malpha|+1-|\mbeta|\geq 0$
again by \eqref{cond-app} and \eqref{cond-app-g}. 
Therefore, by Theorem \ref{m1}, $f$ and $g$ are ${\rm L}^p$-multipliers for any $p\in\oi1\infty$, and the norm of the corresponding 
Fourier multiplier operators is independent of $\mkappa$.

For $s\in(0,1)$ the symbols are given by $h\circ \#$, 
where $h(x)=x^s$ and we us $\#$ to denote either $f$ or $g$, i.e.~$\#\in\{f,g\}$. 
By the Marcinkiewicz theorem and the generalised chain rule formula 
(known as the Fa\'a di Bruno formula; see e.g.~\cite{FaadiBruno})
it is sufficient to estimate
$$
\mxi^\malpha h^{(k)}(\#(\mxi)) \prod_{i=1}^k \partial^{\mdelta^i}\#(\mxi) \,,
$$
where $h^{(k)}$ represents the derivative of the $k$-th order, 
$k\in\{1,2,\dots,|\malpha|\}$ and $\mdelta^i\in\mathbb{N}^d_0\setminus\{0\}$
are such that $\sum_{i=1}^{k}\mdelta^i=\malpha$.  
By lemmata \ref{tvrdnja'} and \ref{tvrdnja''}, an arbitrary summand of 
$\partial^{\mdelta^i}\#(\mxi)$ is given by
(up to a constant factor)
$$
\mkappa^{\mbeta^i} \mxi^{\mgamma^i} \frac{1}{|\mxi|^{l_i}}
	\biggl(|\mxi|+\sum\limits_{j=1}^d \kappa_j \xi_j^2 \biggr)^{-|\mdelta^i|-1} \,,
$$
where $\mbeta^i, \mgamma^i, l_i$ satisfy either \eqref{cond-app} or 
\eqref{cond-app-g}, with $\mdelta^i$ in place of $\malpha$. Let us define
\begin{equation*}
\mbeta:=\sum_{i=1}^k \mbeta^i \ , \qquad \mgamma:=\sum_{i=1}^{k}\mgamma^i \ ,
	\qquad l:=\sum_{i=1}^{k} l_i \;.
\end{equation*}
Since the derivative of $h$ of 
the $k$-th order is equal to (up to a constant factor) $x^{s-k}$, we are finally
left to estimate
\begin{equation}\label{eq:lemma_Marc_Bruno}
\mxi^\malpha (\#(\mxi))^{s-k} \mkappa^{\mbeta} \mxi^{\mgamma} \frac{1}{|\mxi|^{l}}
	\biggl(|\mxi|+\sum\limits_{j=1}^d \kappa_j \xi_j^2 \biggr)^{-|\malpha|-k} \;,
\end{equation}
where we have used $\sum_{i=1}^{k}\mdelta^i=\malpha$.

Let us consider first $\#=f$. In this case, using \eqref{cond-app}, 
we have
\begin{align*}
\alpha_j+\gamma_j\geq 2\beta_j\,, \quad j=1,\dots,d\,; \qquad |\malpha| \geq |\mbeta |\,; \qquad |\mgamma|=k+l+|\mbeta|\,,
\end{align*} 
and \eqref{eq:lemma_Marc_Bruno} reads
$$
\mxi^{\malpha+\mgamma} \mkappa^\mbeta \frac{1}{|\mxi|^{l+k-s}}
	\biggl(|\mxi|+\sum\limits_{j=1}^d \kappa_j \xi_j^2 \biggr)^{-|\malpha|-s} \;.
$$
With the analogous approach as in the case $s=1$, one can get that the 
term above is estimated by
$$
\frac{|\mxi|^{|\malpha|+|\mgamma|-2|\mbeta|-l-k+s}}
	{\Bigl(|\mxi|+\sum_{j=1}^{d}\kappa_j\xi_j^2\Bigr)^{|\malpha|-|\mbeta|+s}}
	= \frac{|\mxi|^{|\malpha|-|\mbeta|+s}}
	{\Bigl(|\mxi|+\sum_{j=1}^{d}\kappa_j\xi_j^2\Bigr)^{|\malpha|-|\mbeta|+s}} \leq 1 \,,
$$
where we have used that $|\mgamma|=k+l+|\mbeta|$.

In the case $\#=g$ by \eqref{cond-app-g} we have
\begin{align*} 
\alpha_j+\gamma_j\geq 2\beta_j\,, \quad j=1,\dots,d\,; \qquad |\malpha| +k \geq |\mbeta |\,; \qquad |\mgamma|=k+l+|\mbeta|\,; \qquad \beta_m\geq k \,,
\end{align*} 
which we use in estimating \eqref{eq:lemma_Marc_Bruno} to get
\begin{align*}
\Biggl|\mxi^{\malpha+\mgamma} \mkappa^\mbeta
	\frac{1}{(\kappa_m\xi_m)^{k-s}} & \frac{1}{|\mxi|^{l}}
	\biggl(|\mxi|+\sum\limits_{j=1}^d \kappa_j \xi_j^2 \biggr)^{-|\malpha|-s}\Biggr| \\
&\leq (\kappa_m\xi_m^2)^{\beta_m-k+s}\prod_{\substack{j=1 \\ j\ne m}}^{d}
	(\kappa_j\xi_j^2)^{\beta_j}\frac{|\mxi|^{|\malpha|+|\mgamma|-2|\mbeta|-l}}
	{\Bigl(|\mxi|+\sum_{j=1}^{d}\kappa_j\xi_j^2\Bigr)^{|\malpha|+s}} \\
&\leq\frac{|\mxi|^{|\malpha|-|\mbeta|+k}}
	{\Bigl(|\mxi|+\sum_{j=1}^{d}\kappa_j\xi_j^2\Bigr)^{|\malpha|-|\mbeta|+k}}
	\leq 1 \;.
\end{align*}
Thus, the statement is proven. 

\endproof

\section{Acknowledgements}\label{sec:thanks}
This work was supported in part by the Croatian Science Foundation under projects 
IP-2018-01-2449 (MiTPDE) and UIP-2017-05-7249 (MANDphy), and by the projects P33594 and M 2669 Meitner-Programm  of the Austrian Science Fund FWF.


Permanent address of D.~Mitrovi\'c is University of Montenegro.

\end{document}